\documentclass[11pt,reqno]{amsart}
\usepackage{enumerate}
\usepackage[sorted,compressed-cites,sorted-cites,initials]{amsrefs}
\usepackage{amsfonts}
\usepackage{amsmath}
\usepackage{amsthm}
\usepackage{amssymb}
\usepackage{latexsym}
\usepackage{multicol}
\usepackage{verbatim}
\usepackage{amscd}
\usepackage{hyperref}
\usepackage[cmtip,arrow]{xy}
\usepackage{pb-diagram}
\usepackage{pb-xy}
\allowdisplaybreaks[1]

\advance\textwidth by 1.2in \advance\oddsidemargin by -.6in
\advance\evensidemargin by -.6in
\newtheorem*{cor}{Corollary}%[section]
\newtheorem*{lem}{Lemma}
\newtheorem*{prop}{Proposition}

\theoremstyle{definition}
\newtheorem*{defn}{Definition}
\newtheorem{thm}{Theorem}

\newtheorem*{rem}{Remark}

\newtheorem*{ex}{Example}

\newenvironment{pf}{\proof}{\endproof}

\numberwithin{equation}{section}

\DeclareMathOperator{\ad}{ad}

\let\bwdg\bigwedge
\def\bigwedge{{\textstyle\bwdg}}

\begin{document}

\newcommand{\thmref}[1]{Theorem~\ref{#1}}
\newcommand{\secref}[1]{Section~\ref{#1}}
\newcommand{\lemref}[1]{Lemma~\ref{#1}}
\newcommand{\propref}[1]{Proposition~\ref{#1}}
\newcommand{\corref}[1]{Corollary~\ref{#1}}
\newcommand{\remref}[1]{Remark~\ref{#1}}
\newcommand{\defref}[1]{Definition~\ref{#1}}
\newcommand{\er}[1]{(\ref{#1})}
\newcommand{\id}{\operatorname{id}}
\newcommand{\sgn}{\operatorname{sgn}}
\newcommand{\wt}{\operatorname{wt}}
\newcommand{\tensor}{\otimes}
\newcommand{\from}{\leftarrow}
\newcommand{\nc}{\newcommand}
\newcommand{\rnc}{\renewcommand}
\newcommand{\dist}{\operatorname{dist}}
\newcommand{\qbinom}[2]{\genfrac[]{0pt}0{#1}{#2}}
\nc{\cal}{\mathcal} \nc{\goth}{\mathfrak} \rnc{\bold}{\mathbf}
\renewcommand{\frak}{\mathfrak}
\newcommand{\supp}{\operatorname{supp}}
\newcommand{\Ann}{\operatorname{Ann}}
\renewcommand{\Bbb}{\mathbb}
\nc\bomega{{\mbox{\boldmath $\omega$}}} \nc\bpsi{{\mbox{\boldmath
$\Psi$}}}
 \nc\balpha{{\mbox{\boldmath $\alpha$}}}
 \nc\bpi{{\mbox{\boldmath $\pi$}}}

\newcommand{\lie}[1]{\mathfrak{#1}}
\makeatletter
\def\section{\def\@secnumfont{\mdseries}\@startsection{section}{1}%
  \z@{.7\linespacing\@plus\linespacing}{.5\linespacing}%
  {\normalfont\scshape\centering}}
\def\subsection{\def\@secnumfont{\bfseries}\@startsection{subsection}{2}%
  {\parindent}{.5\linespacing\@plus.7\linespacing}{-.5em}%
  {\normalfont\bfseries}}
\makeatother
\def\subl#1{\subsection{}\label{#1}}
 \nc{\Hom}{\operatorname{Hom}}
\nc{\End}{\operatorname{End}} \nc{\wh}[1]{\widehat{#1}}
\nc{\Ext}{\operatorname{Ext}} \nc{\ch}{\text{ch}}
\nc{\ev}{\operatorname{ev}} \nc{\Ob}{\operatorname{Ob}}
\nc{\soc}{\operatorname{soc}} \nc{\rad}{\operatorname{rad}}
\nc{\head}{\operatorname{head}}
\def\Im{\operatorname{Im}}
\def\gr{\operatorname{gr}}
\def\mult{\operatorname{mult}}

 \nc{\Cal}{\cal} \nc{\Xp}[1]{X^+(#1)} \nc{\Xm}[1]{X^-(#1)}
\nc{\on}{\operatorname} \nc{\Z}{{\bold Z}} \nc{\J}{{\cal J}}
\nc{\C}{{\bold C}} \nc{\Q}{{\bold Q}}
\renewcommand{\P}{{\cal P}}
\nc{\N}{{\Bbb N}} \nc\boa{\bold a} \nc\bob{\bold b} \nc\boc{\bold c}
\nc\bod{\bold d} \nc\boe{\bold e} \nc\bof{\bold f} \nc\bog{\bold g}
\nc\boh{\bold h} \nc\boi{\bold i} \nc\boj{\bold j} \nc\bok{\bold k}
\nc\bol{\bold l} \nc\bom{\bold m} \nc\bon{\bold n} \nc\boo{\bold o}
\nc\bop{\bold p} \nc\boq{\bold q} \nc\bor{\bold r} \nc\bos{\bold s}
\nc\bou{\bold u} \nc\bov{\bold v} \nc\bow{\bold w} \nc\boz{\bold z}
\nc\boy{\bold y} \nc\ba{\bold A} \nc\bb{\bold B} \nc\bc{\bold C}
\nc\bd{\bold D} \nc\be{\bold E} \nc\bg{\bold G} \nc\bh{\bold H}
\nc\bi{\bold I} \nc\bj{\bold J} \nc\bk{\bold K} \nc\bl{\bold L}
\nc\bm{\bold M} \nc\bn{\bold N} \nc\bo{\bold O} \nc\bp{\bold P}
\nc\bq{\bold Q} \nc\br{\bold R} \nc\bs{\bold S} \nc\bt{\bold T}
\nc\bu{\bold U} \nc\bv{\bold V} \nc\bw{\bold W} \nc\bz{\bold Z}
\nc\bx{\bold x}
\title[Quivers with relations from algebras of invariants]{Quivers with relations arising from Koszul algebras of $\lie g$-invariants.}
\author{Jacob Greenstein}
\thanks{This work was partially supported by the NSF grant DMS-0654421}
\address{Department of Mathematics, University of
California, Riverside, CA 92521.} 
\email{jacob.greenstein@ucr.edu}\maketitle
\begin{abstract} Let $\lie g$ be a complex simple Lie algebra and let $\Psi$ be an extremal set of positive roots. One associates with~$\Psi$ an infinite dimensional Koszul algebra $\bs_\Psi^{\lie g}$ which is a graded subalgebra of the locally finite part of $((\End \bv)^{op}\tensor S(\lie g))^{\lie g}$,
where $\bv$ is the direct sum of all simple finite dimensional $\lie g$-modules.
We describe the structure of the algebra~$\bs_\Psi^{\lie g}$
explicitly in terms of an infinite quiver with relations for $\lie g$ of types~$A$ and~$C$.
We also describe several infinite families of quivers and finite dimensional algebras arising from this construction.
\end{abstract}

\maketitle

\section*{Introduction}
One of the classical methods in the representation theory is to replace a category one wishes to study by an equivalent category of modules over an associative algebra.
This approach was extensively used in the study of the category~$\cal O$ (cf. for example~\cites{BGG,BGS,BKM,FKM,CPS,Soe,St}) and in many other situations and
led to the introduction of highest weight categories in~\cite{CPS}. The associative algebra in question is usually the endomorphism algebra of a generator or
a co-generator of the category. On the other hand, it is also known that endomorphism algebras often give rise to nice associative algebras
(for example, in the case of the category~$\cal O$ these algebras are Koszul). However, describing them in terms of generators and relations, or in 
terms of quivers with relations, is usually a rather involved task (cf. for example~\cite{St}).

In~\cite{CG} the category~$\cal G$ of graded finite dimensional modules  over the polynomial current algebra~$\lie g[t]=\lie g\tensor \bc[t]$ of
a finite dimensional complex Lie algebra~$\lie g$ was studied. That category can be perceived as a non-semisimple ``deformation'' of the semi-simple category 
of finite dimensional $\lie g$-modules. We proved that this category is highest weight in the sense of~\cite{CPS}. We also
studied a family of quivers arising from the endomorphism algebras of injective co-generators of certain Serre subcategories with finitely many simples in the
cases when they are hereditary. For example, all Dynkin quivers can be realised in this way. We also considered an example where the
endomorphism algebra was not hereditary and computed the relations in that algebra. However, it was already clear from that computation that describing quivers and
relations for these algebras in general would be rather difficult.

The situation becomes more manageable if we pass to the truncated current algebra $\lie g\tensor \bc[t]/(t^2)$ which is isomorphic to 
the semidirect sum of~$\lie g$ with its adjoint representation. The motivation for the study of graded representation of that algebra stems from the fact
that several interesting families of indecomposable objects in~$\cal G$ can be regarded as modules over $\lie g\ltimes \lie g$,
namely the classical limits of Kirillov-Reshetikhin modules for $\lie g$ of classical types (\cite{Kir}) or more generally, of the minimal affinisations (\cites{Ch1,Ch2}).
The category~$\cal G_2$ of graded modules over $\lie g\ltimes\lie g$ was studied in~\cite{CG1}. In particular, we studied families of Serre subcategories
of~$\cal G_2$ associated with sets of roots maximising some linear functional. We call these sets extremal since they correspond to
faces of the convex hull of roots of~$\lie g$. A study of these subcategories was motivated by the observation that, after~\cites{Ch2,CM}, 
there is a natural extremal set of positive roots associated with a Kirillov-Reshetikhin module.
Extremal sets have many interesting combinatorial properties and were studied in~\cite{CRD} (in particular, their complete list for~$\lie g$ of classical types
was provided). Given an extremal set~$\Psi$ contained in a fixed set of positive roots of~$\lie g$, one obtains a family of Serre subcategories which have enough projectives and for which the endomorphism algebra of a projective generator is Koszul. Then one constructs an infinite dimensional Koszul algebra~$\bs_\Psi^{\lie g}$
which is ``approximated'' by these finite dimensional Koszul algebras. The advantage of this infinite dimensional algebra is that it allows us to study all
these finite dimensional subalgebras simultaneously.

The aim of the present paper is to describe the structure of algebras~$\bs_\Psi^{\lie g}$. 
We show that they can be realised as path algebras of quite nice quivers with relations.
In some cases these quivers admit very explicit combinatorial presentations. We compute all relations in these algebras for~$\lie g$ of types~$A$ and~$C$.
Quite expectedly, that turns out to be rather difficult and uses monomial bases of the universal enveloping algebra of the lower triangular part of~$\lie g$. 
Due to very restrictive properties of extremal sets, in types~$A$ and~$C$ we can perform all computations using only the monomial bases in type~$A$ which are known very
explicitly (\cite{Lit}). On the other hand, it is quite remarkable that to study the relations in~$\bs_\Psi^{\lie g}$ we only need the most elementary properties
of the extremal sets described in~\cite{CG1}. It should also be noted that, although an extremal set
is conjugate under the action of the Weyl group to the set of roots of an abelian ideal in a suitable Borel subalgebra (cf.~\cite{CRD}), the algebras $\bs_\Psi^{\lie g}$ behave quite differently even for conjugate sets~$\Psi$. For example, depending on whether the highest root of~$\lie g$ is contained in~$\Psi$, all connected
subalgebras of~$\bs^{\lie g}_\Psi$ are infinite or finite dimensional. Another example is discussed in~\ref{G2-20}. 

The paper is organised as follows. In Section~\ref{MR} we briefly review the construction of the algebras~$\bs_\Psi^{\lie g}$ and present the main results. 
In Section~\ref{RELS} we develop the technique for computing relations, while in Section~\ref{FEX} we consider several relatively simple examples which illustrate how
these methods are applied. In Section~\ref{RECF} we construct a family of elements in the universal enveloping algebra of a Borel subalgebra of
$\lie g$ corresponding to parabolic subalgebras with the Levi factor of type~$A$ which play the central role in our computations. Finally, in
Sections~\ref{A} and~\ref{C} we undertake a systematic study of relations in the algebras $\bs_\Psi^{\lie g}$ for~$\lie g$ of types~$A$ and~$C$. We also describe several infinite families of quivers arising from the study of connected subalgebras of~$\bs_\Psi^{\lie g}$ when~$\Psi$ satisfies some ``regularity'' condition.

\subsection*{Acknowledgements}
The principal part of this paper was written while the author was visiting the Weizmann Institute of Science. It is a pleasure to thank
Anthony Joseph for his hospitality and support. The author thanks Arkady Berenstein,
Maria Gorelik, Vladimir Hinich, Bernhard Keller, Anna Melnikov and Shifra Reif for numerous interesting discussions. 
%\newpage

\section{Main results}\label{MR}
Throughout this paper we denote by~$\bz_+$ the set of non-negative integers and by~$\bc$ the field of complex numbers. We consider $\bz_+\cup\{+\infty\}$
as a totally ordered semigroup with $+\infty> n$ and $+\infty+n=+\infty$ for all~$n\in\bz_+$.
All algebras and vector spaces
are considered over~$\bc$. Tensor products and $\Hom$ spaces are taken over~$\bc$ unless specified otherwise. For an associative algebra~$A$,
$A^{op}$ denotes its opposite algebra. For a vector space~$V$, $V^*=\Hom(V,\bc)$. Given a Lie algebra~$\lie a$, we denote by~$ U(\lie a)$ its universal enveloping algebra and by~$ U(\lie a)_+$ the augmentation ideal in~$ U(\lie a)$. In particular,
if~$\lie a$ is abelian, $U(\lie a)$ is the symmetric algebra~$S(\lie a)$. Given an $\lie a$-module $V$ we denote by $V^{\lie a}$
the subspace of~$\lie a$-invariant elements in~$V$, that is~$V^{\lie a}=\{ v\in V\,:\, x v=0\, \forall x\in\lie a\}$.

\subsection{}\label{MR10}
Let~$\lie g$ be a finite dimensional simple complex Lie algebra and fix its Cartan subalgebra~$\lie h$. The Killing form of~$\lie g$
induces a non-degenerate bilinear form $(\cdot,\cdot)$ on~$\lie h^*$. Let~$P\subset \lie h^*$ be a weight lattice and
let~$R\subset P$ be the set
of roots of~$\lie g$ with respect to~$\lie h$. 
Choose the set of simple roots~$\alpha_i\in R$, $i\in I:=\{1,\dots,\dim\lie h\}$ and the corresponding fundamental
weights $\varpi_i\in P$.
Let~$P^+\subset P$  be the $\bz_+$-span of the~$\varpi_i$ and let $R^+$ be the intersection of~$R$ with the~$\bz_+$-span of the~$\alpha_i$.
Given~$\beta\in R$, set for all~$i\in I$
$$
\varepsilon_i(\beta)=\max\{ t\in\bz_+\,:\,\beta+t\alpha_i\in R\},\qquad 
\varphi_i(\beta)=\max\{t\in\bz_+\,:\, \beta-t\alpha_i\in R\}
$$
and define
$$
\varepsilon(\beta):=\sum_{i\in I} \varepsilon_i(\beta)\varpi_i,\qquad
\varphi(\beta):=\sum_{i\in I}\varphi_i(\beta)\varpi_i.
$$
Clearly, $\varepsilon(\beta),\varphi(\beta)\in P^+$. It is well-known that~$\varphi(\beta)=\varepsilon(\beta)+\beta$.
For~$\alpha\in R$ let~$\lie g_\alpha$ be the corresponding root subspace of~$\lie g$ and, given~$\Psi\subset R^+$, let~$\lie n^\pm_\Psi=\bigoplus_{\alpha\in \Psi}
\lie g_{\pm\alpha}$. In particular, we write $\lie n^\pm=\lie n^\pm_{R^+}$ and set~$\lie b=\lie h\oplus\lie n^-$.

We say that~$\Psi\subset R$ is {\em extremal} if there exists~$\xi\in P$ such that
$$
\Psi=\{ \alpha\in R\,:\, (\xi,\alpha)=\max_{\beta\in R} (\xi,\beta)\}.
$$
Geometrically, an extremal subset of~$R$ is the intersection with~$R$ of a face of the convex hull of~$R$ in the euclidean space spanned by~$R$.
Note that if~$\xi\in P^+$ then~$\Psi\subset R^+$.
\label{*}
We will need the following property of extremal sets.
\begin{lem}[\cite{CG1}*{Lemma~2.3}]
Let~$\Psi\subset R$ be extremal and suppose that
$$
\sum_{\alpha\in R} m_\alpha \alpha=\sum_{\beta\in\Psi} n_\beta \beta,\qquad m_\alpha,n_\beta\in\bz_+.
$$
Then
\begin{equation}
\sum_{\beta\in\Psi} n_\beta\le \sum_{\alpha\in R} m_\alpha
\end{equation}
with equality if and only if $m_\alpha=0$ for all~$\alpha\notin\Psi$.\qed
\end{lem}
We note the following
\begin{cor}
Let~$\Psi\subset R$ be extremal. Then~$\Psi+\Psi\cap (R\cup\{0\})=\emptyset$ and
\begin{equation*}
\alpha,\beta\in R,\,\alpha+\beta\in\Psi+\Psi\implies \alpha,\beta\in\Psi.\tag*{\qedsymbol}
\end{equation*}
\end{cor}
\begin{rem}
It is shown in~\cite{CRD} that this property characterises extremal sets. 
\end{rem}

\subsection{}\label{MR30}
Let~$\Psi\subset R^+$ be extremal. In~\cite{CG1} 
two infinite dimensional
Koszul algebras~$\bs_\Psi^{\lie g}$ and~$\be_\Psi^{\lie g}$ were constructed and it was shown that~$(\be_\Psi^{\lie g})^{op}$ is the quadratic dual of~$\bs_\Psi^{\lie g}$
and the left global dimension of~$\bs_\Psi^{\lie g}$ equals~$|\Psi|$. This construction
was motivated by the study of categories of graded representations of current algebras initiated in~\cite{CG}.

Given~$\lambda\in P^+$, let~$V(\lambda)$ be the unique, up to an isomorphism, simple finite dimensional $\lie g$-module of highest weight~$\lambda$.
Let
$$
\bv=\bigoplus_{\lambda\in P^+} V(\lambda),\qquad \bv^\circledast =\bigoplus_{\lambda\in P^+} V(\lambda)^*.
$$
Then $\bv^{\circledast}\tensor \bv$ with the product given by
$$
(f\tensor v)(g\tensor w)=g(v) f\tensor w,\qquad f,g \in\bv^\circledast,\, v,w\in \bv
$$
is isomorphic to a subalgebra of the associative algebra~$(\End_{\bc} \bv)^{op}$
and hence for any associative algebra~$A$, the space $\ba=A\tensor \bv^{\circledast}\tensor \bv$ has a natural structure of an associative algebra.
Moreover, if~$A=\bigoplus_{n\in\bz_+} A[n]$ is a $\bz_+$-graded associative algebra, we obtain a grading on~$\ba$ by assigning to the elements of~$\bv^\circledast\tensor\bv$
the grade zero, that is, $\ba[k]=A[k]\tensor \bv^\circledast\tensor \bv$.
In the rest of the paper, we identify the algebra $\ba$
with $\bv^{\circledast}\tensor A\tensor \bv$ under the natural isomorphism of~$\lie g$-modules and with the induced algebra structure given by $$(f\otimes a\otimes v)(g\otimes b\otimes w)= g(v)f\otimes ab\otimes w,\qquad a,b\in A,\, f,g\in\bv^{\circledast},\, v,w\in\bv .$$

We write~$\bt$ (respectively, $\bs$, $\be$) for $\ba$ with~$A$ the tensor algebra $T(\lie g)$ of~$\lie g$ (respectively, 
the symmetric algebra $S(\lie g)$ and the exterior algebra $\bigwedge \lie g$). In particular, in these cases $A$ is a $\lie g$-module with respect to the diagonal action, hence
$\ba$ is also a $\lie g$-module and the multiplication is a homomorphism of $\lie g$-modules. It follows that $\ba^{\lie g}$ is a subalgebra of~$\ba$.
From now on, we let~$\ba$ be one of the algebras $\bt$, $\bs$ or~$\be$. Given~$\lambda\in P^+$, the algebra $\ba^{\lie g}$ contains
a primitive idempotent~$1_\lambda$ corresponding to the canonical $\lie g$-invariant element  in~$V(\lambda)^*\tensor V(\lambda)$, or,
equivalently to the identity element in~$\End V(\lambda)$. Then we have
$$
\ba^{\lie g}=\bigoplus_{\lambda,\mu\in P^+} 1_\lambda\ba^{\lie g} 1_\mu,\qquad
1_\lambda \ba^{\lie g} 1_\mu=(V(\lambda)^*\tensor A\tensor V(\mu))^{\lie g}.
$$

\subsection{}\label{MR50}
Given~$\Psi\subset R^+$, define a relation $\le_\Psi$ on~$P$ by $\lambda\le_\Psi\mu$ if~$\mu-\lambda\in\bz_+\Psi$. It is straightforward to check that
$\le_\Psi$ is a partial order. In particular, $\le:=\le_{R^+}$ is the standard partial order on~$P$.
If~$\lambda\le_{\Psi} \mu$ and~$\lambda\not=\mu$ we write~$\lambda<_\Psi\mu$. 
 Note that for all~$\lambda\in P$ and for all~$\Psi\subset R^+$, the set~$\{\mu\in P^+\,:\, \mu\le_\Psi
\lambda\}$ is finite.
Define a function $d_\Psi: \{ (\lambda,\mu)\in P^+\times P^+\,:\, \lambda\le_\Psi\mu\}\to \bz_+$ by
$$
d_\Psi(\lambda,\mu)=\min\{ \sum_{\alpha\in\Psi} m_\alpha\,:\, \mu-\lambda=\sum_{\alpha\in\Psi} m_\alpha\alpha,\, m_\alpha\in\bz_+\}.
$$
Clearly, $d_\Psi(\lambda,\mu)=0$ if and only if~$\lambda=\mu$ and $d_\Psi(\lambda,\mu)+d_\Psi(\mu,\nu)\le d_\Psi(\lambda,\nu)$ 
for all $\lambda\le_\Psi\mu\le_\Psi\nu$. Furthermore, if~$\Psi$ is extremal, we have
$$
d_\Psi(\lambda,\mu)+d_\Psi(\mu,\nu)=d_\Psi(\lambda,\nu)
$$
and if $\mu$ covers~$\lambda$ then~$d_\Psi(\lambda,\mu)=1$. In particular, in this case $d_\Psi$ is the unique distance function for
the poset $(P^+,\le_\Psi)$.

Fix an extremal set~$\Psi\subset R^+$. Given~$F\subset P^+$, define
$$
\ba_\Psi^{\lie g}(F)=\bigoplus_{\lambda,\mu\in F\,:\, \lambda\le_\Psi\mu} 1_\lambda \ba^{\lie g}[d_\Psi(\lambda,\mu)] 1_\mu.
$$
It is not hard to check that~$\ba_\Psi^{\lie g}(F)$ is a subalgebra of~$\ba^{\lie g}$. Let~$\ba_\Psi^{\lie g}:=\ba_\Psi^{\lie g}(P^+)$.
Given~$\lambda,\mu\in P^+$, define the following subsets of~$P^+$
$$
\le_\Psi\lambda=\{ \nu\in P^+\,:\, \nu\le_\Psi\lambda\},\qquad \mu\le_\Psi=\{\mu\le_\Psi\nu\,:\nu\in P^+\}
$$
and~$[\mu,\lambda]_\Psi=(\le_\Psi\lambda)\cap (\mu\le_\Psi)$. We say that~$F\subset P^+$ is interval closed in the partial order~$\le_\Psi$ if~$\lambda,\mu\in F$ implies that~$[\lambda,\mu]_\Psi\subset F$.

The main properties of the algebras~$\bs^{\lie g}_\Psi$ established in~\cite{CG1} are summarised below.
\begin{thm}[\cite{CG1}*{Theorem~1}] Let~$\Psi$ be an extremal set of positive roots. 
 \begin{enumerate}[{\rm(i)}]
\item\label{kosalg.i} Let $\mu,\nu\in P^+$. The  subalgebras $\bs_\Psi^{\lie g}(\le_\Psi\nu)$,
$\bs_\Psi^{\lie g}(\mu\le_\Psi)$ and $\bs_\Psi^{\lie g} ([\mu,\nu]_\Psi)$ of
$\bs_\Psi^{\lie g}$ are Koszul and have global dimension at most $|\Psi|$.
The   bound is attained for some $\mu',\nu'\in P^+$ with $\mu'\le_\Psi\nu'$.

\item\label{kosalg.ii} The algebra $\bs_\Psi^{\lie g}$ is Koszul and has left global dimension $|\Psi|$.
\end{enumerate}
\end{thm}
\begin{rem}
The argument of~\cite{CG1}*{Proposition~4.5} actually proves that~$\bs_\Psi^{\lie g}(F)$ is Koszul for any~$F\subset P^+$
interval closed in the partial order~$\le_\Psi$.
\end{rem}

\subsection{}\label{MR65}\label{PRE0}
Being Koszul, the algebras $\bs^{\lie g}_\Psi$ 
are quadratic and so to describe all relations in them it is enough to describe the quadratic relations. 
A convenient language for that is provided by quivers. We mostly follow the conventions from~\cite{RinBook}. Let us briefly review the
quiver terminology which will be used in the sequel.

Recall that
a quiver~$\Delta$ is a pair~$\Delta=(\Delta_0,\Delta_1)$ where~$\Delta_0$ is the vertex set, $\Delta_1$ is the set of arrows.
{\em In this paper we only consider quivers without multiple arrows}, that is, for any pair $x,y\in\Delta_0$, there is at most one arrow $x\from y\in\Delta_1$ (in
other words, $\Delta_1$ identifies with a subset of~$\Delta_0\times\Delta_0$).
A path of length~$k$ in such a quiver is a sequence $x_0,\dots,x_k\in\Delta_0$ such that for all~$0\le i<k$, there is an arrow $x_i\from x_{i+1}\in\Delta_1$. 
Denote by~$\Delta(x,y)$ the set of all paths in~$\Delta$ from~$y$ to~$x$.
With every vertex~$x\in \Delta_0$ we associate a trivial path~$1_x$ of length~$0$.  

The opposite quiver~$\Delta^{op}$ of~$\Delta$ is the quiver with the same vertex set obtained by reversing all arrows. 
The underlying graph~$\bar \Delta$ of~$\Delta$ is obtained from~$\Delta$ by forgetting the orientation of the arrows. We say that~$\Delta$ is connected if~$\bar \Delta$ is connected. 
A connected quiver $\Delta$ is said to be of type $\mathbb X$ (respectively, of type~$\tilde{\mathbb X}$), where $X=ADE$ if $\bar\Delta$ is the Dynkin diagram (respectively, extended 
Dynkin diagram) of a simple finite dimensional Lie algebra of type~$X$.

A vertex $x$ is said to be a direct successor (respectively, predecessor) of~$y$ if there is an arrow $x\from y$ (respectively, $y\from x$) in~$\Delta_1$. The
set of direct successors (predecessors) of~$x\in \Delta_0$ is denoted by~$x^+$ (respectively, $x^-$).
A vertex~$x\in \Delta_0$ is called a source if~$x^-=\emptyset$ and a sink
if~$x^+=\emptyset$.

Given~$\Delta_0'\subset \Delta_0$, the {\em full subquiver} of~$\Delta$ defined by~$\Delta_0'$ is~$\Delta'=(\Delta'_0,\Delta'_1)$ where
$\Delta_1'$ is the set of all arrows in~$\Delta_1$ with starting and ending points in~$\Delta_0'$.
A subquiver~$\Delta'$ of~$\Delta$ is called {\em convex} if for any vertices $x,y\in \Delta_0'$ we have ${\Delta'}(x,y)={\Delta}(x,y)$, that is
a path in~$\Delta$ from~$y$ to~$x$ is completely contained in~$\Delta'$.
In particular, a convex subquiver is full. 
A connected component of~$\Delta$ is a subquiver~$\Delta'$ such that~$\overline{\Delta'}$ is a connected
component of~$\bar\Delta$. Then~$\Delta$ is a disjoint union of its connected components.
Given~$x\in\Delta_0$, we denote the connected component of~$\Delta$ containing~$x$ by~$\Delta[x]$.

A full embedding of quivers $\Delta\to \Delta'$ is a pair of injective maps $F_0:\Delta_0\to\Delta_0'$ and~$F_1:\Delta_1\to\Delta_1'$
which are compatible in a natural way and such that the $(F_0(\Delta_0),F_1(\Delta_1))$ is a full
subquiver of~$\Delta'$. If both maps are  bijective we say that~$\Delta$ is
isomorphic to~$\Delta'$.

Given a quiver~$\Delta=(\Delta_0,\Delta_1)$, let~$\bc\Delta$ be the complex vector space with the basis consisting of all finite paths in~$\Delta$. The product 
of two paths is set to be their composition when they are composable, and zero otherwise. This defines on~$\bc\Delta$ a structure of a $\bz_+$-graded associative algebra,
the grading being given by the length of paths. 
In particular, the $1_x$, $x\in\Delta_0$ are primitive orthogonal idempotents and~$\bc\Delta[0]$ is commutative and semi-simple.
Clearly, $\bc(\Delta^{op})\cong (\bc\Delta)^{op}$.
The group~$(\bc^\times)^{\Delta_1}$ acts naturally on~$\bc\Delta[1]$ and for all~$\boldsymbol{z}\in (\bc^\times)^{\Delta_1}$ the action of~$\boldsymbol{z}$ extends to an automorphism of~$\bc\Delta$ preserving the grading and
the~$1_x$, $x\in\Delta_0$. Clearly, an isomorphism of quivers induces an isomorphisms of the corresponding path algebras.

A relation on~$\Delta$ is a linear combination of paths from~$x$ to~$y$ for some~$x,y\in\Delta_0$. In particular, a relation of the form $p$, where $p$
is a path,
is called a zero relation, while a relation of the form $p-p'$ is called a commutativity relation.
Given a quiver~$\Delta$ and a set of relations~$\cal R$, we can form an algebra~$\bc(\Delta,\cal R)=\bc(\Delta,V):=\bc \Delta/\langle \cal R \rangle$, where~$V$ 
is the vector subspace of~$\bc \Delta$ spanned by~$\cal R$. This algebra is often referred to
as the path algebra of the quiver~$\Delta$ with relations~$\cal R$.

\subsection{}\label{MR70}
We define the infinite quiver~$\Delta_\Psi$ as 
\begin{align*}
&(\Delta_\Psi )_0=P^+\\
&(\Delta_\Psi)_1=\{ (\lambda,\mu)\in P^+\times P^+\,:\, \mu-\lambda=\beta\in\Psi,\, \lambda-\varepsilon(\beta)=\mu-\varphi(\beta)\in P^+\}.
\end{align*}
Thus, $(\Delta_\Psi)_1$ is a subset of the cover relation in $(P^+,\le_\Psi)$. 
It is immediate that if there is a path from~$\mu$ to~$\lambda$ in~$\Delta_\Psi$, then
$\lambda\le_\Psi\mu$ and the length of any such path is~$d_\Psi(\lambda,\mu)$. In particular, the quiver~$\Delta_\Psi$ has no loops or oriented cycles.
Since for all~$\lambda\in P^+$ the set~$\le_\Psi\lambda$ is finite,
it follows that every vertex in~$\Delta_\Psi$ is connected to a sink.
Given $F\subset P^+$, denote~$\Delta_\Psi(F)$ the full subquiver of~$\Delta_\Psi$ defined by~$F$. If~$F$ is interval closed in the partial order~$\le_\Psi$,
$\Delta_\Psi(F)$ is convex. 

\begin{prop}
Let~$F\subset P^+$ be interval closed in the partial order~$\le_\Psi$. 
There exists a natural isomorphism of~$\bz_+$-graded associative algebras $\bt^{\lie g}_\Psi(F)\to
\bc\Delta_\Psi(F)$. This isomorphism is unique up to an automorphism of~$\bc\Delta_\Psi(F)$ extending the natural action of~$(\bc^\times)^{(\Delta_\Psi(F))_1}$ on~$\bc\Delta_\Psi(F)[1]$.
\end{prop}

\subsection{}\label{MR75}
As proved in~\cite{CG1}*{Lemma~4.2}, for all $F\subset P^+$ which is interval closed in the partial order~$\le_\Psi$,
$\bs^{\lie g}_\Psi(F)$ is isomorphic to the quotient of~$\bt^{\lie g}_\Psi(F)$ by a quadratic ideal
and 
$$
\ker(\bt^{\lie g}_\Psi(F)\to \bs^{\lie g}_\Psi(F))=\bt^{\lie g}_\Psi(F)\cap \ker(\bt_\Psi^{\lie g}\to\bs_\Psi^{\lie g}).
$$
Fix an isomorphism~$\Phi:\bt^{\lie g}_\Psi\to\bc\Delta_\Psi$.
Then~\propref{MR70} allows us to identify the idempotents (respectively, some fixed generators of degree~$1$) of~$\bs^{\lie g}_\Psi$ with vertices (respectively, arrows) in the quiver~$\Delta_\Psi$. 
To describe the quadratic relations, we need
to consider, for all $\lambda,\mu\in P^+$ with $\lambda\le_\Psi \mu$ and~$d_\Psi(\lambda,\mu)=2$, that is, for all $\lambda,\lambda+\eta\in P^+$, $\eta\in\Psi+\Psi$,
the subquivers~$\Delta_\Psi([\lambda,\lambda+\eta]_\Psi)$ of~$\Delta_\Psi$ 
and the subalgebras $\bs_\Psi^{\lie g}([\lambda,\lambda+\eta]_\Psi)$
of~$\bs_\Psi^{\lie g}$. 

Denote by $\mathfrak R_\Psi(\lambda,\lambda+\eta)$ the image of the canonical map 
$\ker(\bt^{\lie g}_\Psi([\lambda,\lambda+\eta]_\Psi)\twoheadrightarrow \bs^{\lie g}_\Psi([\lambda,\lambda+\eta]_\Psi))$ in
$\bc\Delta_\Psi([\lambda,\lambda+\eta]_\Psi)$
under~$\Phi$. We set~$\mathfrak R_\Psi(\lambda,\lambda+\eta)=0$ if~$\lambda+\eta\notin P^+$.

\subsection{}\label{M75}
Let $\eta\in \Psi+\Psi$ and set
$$
m_\eta=\#\{ (\beta,\beta')\in\Psi\times\Psi\,:\, \beta+\beta'=\eta\}.
$$
Note that~$m_\eta=1$ implies that~$\eta\in2 \Psi$.
For all~$\lambda\in P^+$, let $t_{\lambda,\eta}=\dim 1_\lambda\bt^{\lie g}_\Psi 1_{\lambda+\eta}$ if~$\lambda+\eta\in P^+$
and set~$t_{\lambda,\eta}=0$ otherwise. Since by~\propref{MR70}, $t_{\lambda,\eta}$ equals the number of paths from~$\lambda+\eta$
to~$\lambda$ in~$\Delta_\Psi$, it is immediate that
$t_{\lambda,\eta}\le m_\eta$ for all~$\lambda\in P^+$. 
\begin{defn}
An extremal set~$\Psi\subset R^+$ is said to be {\em regular} if for all~$\eta\in\Psi+\Psi$ and for all~$\lambda\in P^+$, $t_{\lambda,\eta}>0\implies t_{\lambda,\eta}=m_\eta$.
\end{defn}

The quiver 
$\Delta_\Psi([\lambda,\lambda+\eta]_\Psi)$ identifies with the quiver
\begin{equation}\label{std.quiv}
\Gamma(t)=\dgHORIZPAD=2pt \dgVERTPAD=2pt\divide\dgARROWLENGTH by 3\multiply\dgARROWLENGTH by 2
\def\dgeverynode{\scriptstyle}
\begin{diagram}
\node{}\node{}\node{t+1}\arrow{sww}\arrow{sw}\arrow{se}\arrow{see}\\
\node{1}\arrow{see}\node{2}\arrow{se}\node{\displaystyle\cdots}\node{t-1}\arrow{sw}\node{t}\arrow{sww}\\
\node{}\node{}\node{0} 
\end{diagram}
\end{equation}
with~$t=t_{\lambda,\eta}$ paths $\bop_i=(0\from r\from t+1)$, $1\le r\le t$ of length~$2$. 

Let~$V$ be a $k$-dimensional subspace of~$\bc\{\bop_1,\dots,\bop_t\}$. 
We say that $V$ is generic if it is generic with respect to
any coordinate flag corresponding to the basis $\bop_1,\dots,\bop_t$, that is
for all~$1\le i_1<\cdots<i_r\le t$, $1\le r\le t$ we have 
$$\dim (V\cap \bc\{ \bop_{i_1},\dots,\bop_{i_r}\})=\begin{cases}
                                              0,& 1\le r< t-k,\\
					      r+k-t,&t-k\le r\le t.
                                             \end{cases}
$$
In particular, if~$t=1$, $V$ is generic if and only if~$\dim V=1$.
For instance, if~$t=2$, $\dim V=1$ and~$V$ then
$\bc(\Gamma(t),V)$ is of finite type. However, it has different isomorphism classes of indecomposable modules,
depending on $V$ being or not being generic.
If~$t=3$, $\dim V=1$ and~$V$ is generic
then~$\bc(\Gamma(t),V)$ is unique up to an isomorphism, is canonical (cf.~\cite{RinBook}*{\S3.7}),  
of tubular type~$\mathbb D_4$ and tame. If~$t=4$, $\dim V$ is generic and $\dim V=2$ when
we can assume, without loss of generality, that $V$ is spanned by $\bop_1+\bop_2+\bop_3$, $\bop_1+z\bop_2+\bop_4$ for some $z\in\bc^\times$. In particular,
we have a family of algebras parametrised by elements of~$\mathbb P^1$. The algebra~$\bc(\Gamma(t),V)$ is again canonical,
of tubular type~$\tilde{\mathbb D}_4$, and is tame (cf.~\cite{RinBook}). In these cases the module categories of~$\bc(\Gamma(t),V)$ are 
described completely (\cite{RinBook}).
If~$V$ is not generic, $\bc(\Gamma(t),V)$ it is still tame (cf.~\cite{Gei}).
If~$t>4$, it is easy to see, using~\cite{dlPen}*{Proposition~1.3},
that $\bc(\Gamma(t),V)$ is wild for all choices of~$V$ of dimension~$\lfloor t/2\rfloor$.

From now on, we identify $\xi\in\lie h^*$ with the canonical algebra homomorphism $S(\lie h)\to \bc$ extending~$\xi$.
Let
$$
\cal N_\eta=\{ \lambda\in P^+\,:\, \text{$t_{\lambda,\eta}>0$, $\mathfrak R_\Psi(\lambda,\lambda+\eta)$ is not generic}\}.
$$
We can now formulate our main result.
\begin{thm}\label{mainthm}
Suppose that~$\lie g$ is of type~$A$ or~$C$ and
let~$\Psi$ be an extremal set of positive roots, $|\Psi|>1$. 
\begin{enumerate}[{\rm(i)}]
\item\label{mainthm.i} The algebra $\bs_\Psi^{\lie g}$ is isomorphic to the quotient of the path algebra of the quiver~$\Delta_\Psi$
by the ideal generated by the spaces $\mathfrak R_\Psi(\lambda,\lambda+\eta)$, $\eta\in\Psi+\Psi$, $\lambda,\lambda+\eta\in P^+$.

Fix~$\eta\in \Psi+\Psi$. 
\item\label{mainthm.ii'} If $m_\eta=1$  then $\mathfrak R_\Psi(\lambda,\lambda+\eta)=0$ for all~$\lambda\in P^+$.
\item\label{mainthm.ii} If $t_{\lambda,\eta}>1$, then~$\dim\mathfrak R_\Psi(\lambda,\lambda+\eta)=\lfloor t_{\lambda,\eta}/2\rfloor>0$.
\item\label{mainthm.iii} Suppose that~$m_\eta>1$. 
Then $\cal N_\eta$ is contained
in a Zariski closed subset of~$\lie h^*$. Moreover $\cal N_\eta\cap \{\lambda\in P^+\,:\, t_{\lambda,\eta}=2,3\}=\emptyset$ and
if~$\Psi$ is regular, then either $\cal N_\eta=\emptyset$ or there exists 
a linear polynomial $H_\eta\in S(\lie h)$ such that
$$
\cal N_\eta=P^+\cap \{ \xi\in\lie h^*\,:\, \xi(H_\eta)=0\}.
$$
\end{enumerate}
\end{thm}
Analysis of other examples allows us to conjecture that
a similar result should hold for~$\lie g$ of all types.

The above theorem is established in Propositions~\ref{A95} and~\ref{A110} for~$\lie g$ of type~$A$ and in Propositions~\ref{C120}, \ref{C140},
\ref{C160} and~\ref{C180} for~$\lie g$ of type~$C$. In fact, we do not just establish the genericity of the spaces~$\mathfrak R_\Psi(\lambda,\lambda+\eta)$
but also compute the relations explicitly. Needless to say, as we write the relations as linear combinations of paths, the
specific coefficients we obtain depend on
a fixed isomorphism~$\Phi$, or equivalently, on the choice of generators of degree one in~$\bs_\Psi^{\lie g}$, which are unique up to non-zero scalars,
while the genericity of the spaces~$\mathfrak R_\Psi(\lambda,\lambda+\eta)$ is independent of that choice. 
We choose~$\Phi$ so that the relations for~$\eta\in\Psi+\Psi$ and~$\lambda\in P^+$
satisfying~$t_{\lambda,\eta}=m_\eta=2$, are the commutativity relations. 

Let us briefly explain how to compute relations in~$\be_\Psi^{\lie g}$ from those in~$\bs_\Psi^{\lie g}$. 
There is a natural map~$\langle\cdot,\cdot\rangle:\bc\Delta_\Psi\tensor \bc\Delta_\Psi^{op}\to \bc$,
such that~$\langle (\bc\Delta_\Psi)[k],(\bc\Delta^{op}_\Psi)[r]\rangle=0$, $k\not=r$, 
$$\langle 1_\lambda,1_\mu\rangle=\delta_{\lambda,\mu},\qquad 
\langle \lambda_1\from\cdots \from\lambda_k,\mu_k\from\cdots\from \mu_1\rangle=\delta_{\lambda_1,\mu_1}\cdots\delta_{\lambda_k,\mu_k}.
$$
It is not hard to see from~\cite{CG1}*{Proposition~5.3} that~$\be_\Psi^{\lie g}$ is isomorphic to the quotient of~$\bc\Delta_\Psi^{op}$
by the ideal generated by the spaces~$\mathfrak R_\Psi(\lambda,\lambda+\eta)^!=\{ x\in\bc\Delta_\Psi^{op}\,:\, \langle \mathfrak R_\Psi(\lambda,\lambda+\eta),x\rangle=0\}$.

\subsection{}\label{M95}
We conclude this section with a description of an infinite family of quivers arising from this construction.

Given~$\boldsymbol{x}=(x_1,\dots,x_r)\in\bz_+^r$, let~$|\boldsymbol{x}|=\sum_{j=1}^r x_j$. Set~$\boldsymbol{e}_i^{(r)}=(\delta_{i,j})_{1\le j\le r}
\in\bz_+^r$.
Given~$\boldsymbol{m}=(m_1,\dots,m_r)\in(\bz_+\cup\{+\infty\})^r$,
we define the quiver~$\Xi(\boldsymbol{m})$ as follows.
The vertices of~$\Xi(\boldsymbol{m})$ are
the lattice points in the $r$-dimensional rectangular parallelepiped
$[0,m_1]\times\cdots\times [0,m_r]$. 
Given~$\boldsymbol{x}=(x_1,\dots,x_r)\in\Xi(\boldsymbol{m})_0$, the arrows ending at~$\boldsymbol{x}$ are 
\begin{alignat*}{2}
&\boldsymbol{x}\from \boldsymbol{x}+2\boldsymbol{e}_j^{(r)},&
\qquad &x_j<m_j-1,\, 1\le j\le r
\\
\intertext{and}
&\boldsymbol{x}\from \boldsymbol{x}+\boldsymbol{e}_j^{(r)}+\boldsymbol{e}_k^{(r)},&
&x_i<m_i,\,x_j<m_j,\,1\le i<j\le r.
\end{alignat*}
Let $\Xi_a(\boldsymbol{m})$, $a=0,1$ be the full subquiver of~$\Xi(\boldsymbol{m})$ defined by the set 
$$\{\boldsymbol{x}\in\Xi(\boldsymbol{m})_0\,:\, |\boldsymbol{x}|=a\pmod 2\}.$$ 
It is immediate that $\Xi_a(\boldsymbol{m})$ is a convex subquiver of~$\Xi(\boldsymbol{m})$.

For instance, for $r=2$ and~$m_1=m_2=1$, $\Xi_0(\boldsymbol{m})$ is the quiver of type~$\mathbb A_2$ with the linear orientation and $\Xi_1(\boldsymbol{m})$
has two isolated vertices (in fact, this is the only case when $\Xi_1(\boldsymbol{m})$ is not connected).
For~$m_1=m_2=2$,
$\Xi_0(\boldsymbol{m})$ is
the quiver~\eqref{std.quiv} with~$t=3$, while $\Xi_{1}(\boldsymbol{m})$ is
$$
\divide\dgARROWLENGTH by 3
\multiply\dgARROWLENGTH by 2
\def\dgeverynode{\scriptscriptstyle}\def\dgeverylabel{\scriptscriptstyle}
\dgVERTPAD=2pt \dgHORIZPAD=2pt
\begin{diagram}
\node{(1,2)}\arrow{s}\arrow{e}\node{(1,0)}
\\
\node{(0,1)}\node{(2,1)}\arrow{w}\arrow{n}
\end{diagram}
$$
An example with~$m_1=6$, $m_2=5$ is shown below
$$
\divide\dgARROWLENGTH by 3
\def\dgeverynode{\scriptscriptstyle}\def\dgeverylabel{\scriptscriptstyle}
\dgVERTPAD=2pt \dgHORIZPAD=2pt
\Xi_0(\boldsymbol{m})=
\begin{diagram}
\node{}\node{(5,1)}\arrow{sw}\arrow[2]{s}\node{}\node{(5,3)}\arrow[2]{w}\arrow{sw}\arrow[2]{s}\node{}\node{(5,5)}\arrow[2]{w}\arrow{sw}\arrow[2]{s}\\
\node{(4,0)}\arrow[2]{s}\node{}\node{(4,2)}\arrow[2]{w}\arrow{sw}\arrow[2]{s}\node{}\node{(4,4)}\arrow[2]{w}\arrow{sw}\arrow[2]{s}\node{}\node{(4,6)}\arrow[2]{w}\arrow[2]{s}\arrow{sw}
\\
\node{}\node{(3,1)}\arrow{sw}\arrow[2]{s}\node{}\node{(3,3)}\arrow[2]{w}\arrow{sw}\arrow[2]{s}\node{}\node{(3,5)}\arrow[2]{w}\arrow{sw}\arrow[2]{s}\\
\node{(2,0)}\arrow[2]{s}\node{}\node{(2,2)}\arrow[2]{w}\arrow{sw}\arrow[2]{s}\node{}\node{(2,4)}\arrow[2]{w}\arrow{sw}\arrow[2]{s}\node{}\node{(2,6)}\arrow[2]{w}\arrow[2]{s}\arrow{sw}
\\
\node{}\node{(1,1)}\arrow{sw}\node{}\node{(1,3)}\arrow[2]{w}\arrow{sw}\node{}\node{(1,5)}\arrow[2]{w}\arrow{sw}
\\
\node{(0,0)}\node{}\node{(0,2)}\arrow[2]{w}\node{}\node{(0,4)}\arrow[2]{w}\node{}\node{(0,6)}\arrow[2]{w}
\end{diagram}
\qquad
\Xi_1(\boldsymbol{m})=
\begin{diagram}
\node{(5,0)}\arrow[2]{s}\node{}\node{(5,2)}\arrow[2]{w}\arrow[2]{s}\arrow{sw}\node{}\node{(5,4)}\arrow[2]{w}\arrow[2]{s}\arrow{sw}\node{}\node{(5,6)}
\arrow[2]{w}\arrow[2]{s}\arrow{sw}
\\
\node{}\node{(4,1)}\arrow[2]{s}\arrow{sw}\node{}\node{(4,3)}\arrow[2]{w}\arrow[2]{s}\arrow{sw}\node{}\node{(4,5)}\arrow[2]{w}\arrow[2]{s}\arrow{sw}\\
\node{(3,0)}\arrow[2]{s}\node{}\node{(3,2)}\arrow[2]{w}\arrow{sw}\arrow[2]{s}\node{}\node{(3,4)}\arrow[2]{w}\arrow[2]{s}\arrow{sw}\node{}
\node{(3,6)}\arrow[2]{w}\arrow[2]{s}\arrow{sw}\\
\node{}\node{(2,1)}\arrow{sw}\arrow[2]{s}\node{}\node{(2,3)}\arrow[2]{w}\arrow{sw}\arrow[2]{s}\node{}\node{(2,5)}\arrow[2]{w}\arrow[2]{s}\arrow{sw}\\
\node{(1,0)}\node{}\node{(1,2)}\arrow{sw}\arrow[2]{w}\node{}\node{(1,4)}\arrow{sw}\arrow[2]{w}\node{}\node{(1,6)}\arrow{sw}\arrow[2]{w}\\
\node{}\node{(0,1)}\node{}\node{(0,3)}\arrow[2]{w}\node{}\node{(0,5)}\arrow[2]{w} 
\end{diagram}
$$
Note that in this case~$\Xi_{1}(\boldsymbol{m})\cong \Xi_{0}(\boldsymbol{m})^{op}$ (cf.~\propref{C2}).
For~$r=3$ and $m_1=m_2=m_3=1$, $\Xi_0(\boldsymbol{m})$ (respectively, $\Xi_1(\boldsymbol{m})$) is the quiver of type~$\mathbb D_4$ where the triple node is the unique 
sink (source). 
Finally, $\Xi_0((2,1,1))$ is the quiver~\eqref{std.quiv} with~$t=4$ where~$(0,0,0)$ is the sink and~$(2,1,1)$ is the source, while
$\Xi_0((1,1,1,1))$ is the quiver~\eqref{std.quiv} with~$t=6$, where~$(0,0,0,0)$ is the sink and~$(1,1,1,1)$ is the source.
We prove (cf.~\propref{C2}) that the isomorphism classes of quivers~$\Xi_a(\boldsymbol{m})$ with~$r>1$ are par\-a\-met\-rised by partitions.

\begin{prop}
Suppose that~$\lie g$ is of type~$C$ and~$\Psi$ is regular. Let~$\lambda\in P^+$ and suppose that~$|\lambda^-\cup\lambda^+|>0$.
Then the connected component~$\Delta_\Psi[\lambda]$ of~$\Delta_\Psi$
is isomorphic to~$\Xi_a(\boldsymbol{m})$ for some  $\boldsymbol{m}\in(\bz_+\cup\{+\infty\})^r$, $r>0$ and~$a\in\{0,1\}$.
\end{prop}
In particular, our isomorphism $\bt_\Psi^{\lie g}\to\bc\Delta_\Psi$ induces an isomorphism of a subalgebra of~$\bt_\Psi^{\lie g}$ corresponding
to an interval closed set
onto~$\Xi_a(\boldsymbol{m})$. In particular, 
we can define a family of relations on $\Xi_a(\boldsymbol{m})$, depending on positive integer parameters, which yields an infinite family of
finite dimensional Koszul algebras.

\section{Relations in~\texorpdfstring{$\bs_\Psi^{\lie g}$}{S\_Psi g}}\label{RELS}

\subsection{}\label{PRE100}
Let~$V$ be a $\lie g$-module. Given~$\mu\in\lie h^*$, let
$$
V_\mu=\{v\in V\,:\, hv=\mu(h) v,\, h\in\lie h\}.
$$
If~$V$ is finite dimensional, then~$V=\bigoplus_{\mu\in P} V_\mu$.
Moreover, $V$ is  isomorphic
to a direct sum of simple finite dimensional modules~$V(\lambda)$, $\lambda\in P^+$. In particular, the adjoint representation~$\lie g$ is isomorphic to~$V(\theta)$
where~$\theta$ is the highest root of~$\lie g$.

Fix Chevalley generators~$e_i\in\lie g_{\alpha_i}$, $f_i\in\lie g_{-\alpha_i}$ and $h_i\in \lie h$, $i\in I$ of~$\lie g$.
The module~$V(\lambda)$ is generated by a highest weight vector~$v_\lambda\in V(\lambda)_\lambda$
satisfying
$$
\Ann_{ U(\lie g)} v_\lambda= U(\lie g)(\lie n^+ + \ker \lambda)+ \sum_{i\in I}  U(\lie g) f_i^{\lambda(h_i)+1} 
$$
For each~$\lambda\in P^+$, we fix~$v_\lambda$ once for all and then
we fix $\xi_{-\lambda}\in V(\lambda)^*_{-\lambda}$ such that~$\xi_{-\lambda}(v_\lambda)=1$. Then we have
$$
\Ann_{ U(\lie g)} \xi_{-\lambda}= U(\lie g)(\lie n^- + \ker(-\lambda))+\sum_{i\in I}  U(\lie g) e_i^{\lambda(h_i)+1}.
$$
In particular,
$$
\Ann_{ U(\lie n^-)} v_\lambda=\sum_{i\in I}  U(\lie n^-) f_i^{\lambda(h_i)+1},\qquad \Ann_{ U(\lie n^+)} \xi_{-\lambda}=\sum_{i\in I}  U(\lie n^+) e_i^{\lambda(h_i)+1}.
$$

Given~$\lambda\in P^+$ and a finite dimensional $\lie g$-module~$M$, let
$$
M^{\lambda}=\{m \in M\,:\, \Ann_{ U(\lie n^+)}m\supset \Ann_{ U(\lie n^+)} \xi_{-\lambda}\}.
$$
If~$N$ is a subspace of~$M$, let~$N^\mu=N\cap M^\mu$. We will need the following results (cf.~\cite{PRV}; we use some of them in the form in which
they are presented in~\cite{Jos}, there the corresponding statements are established in the case of integrable modules over quantised enveloping algebras
of Kac-Moody algebras).
\begin{prop}
Let~$\mu,\nu\in P^+$ and let~$M$ be a finite dimensional $\lie g$-module. 
\begin{enumerate}[{\rm(i)}]
\item\label{PRE100.0} $\Hom_{\lie g}(V(\lambda),M)\cong M^{\lie n^+}\cap M_\lambda$.
 \item\label{PRE100.i} $V(\mu)\tensor V(\nu)^*= U(\lie g)(v_\mu\tensor \xi_{-\nu})$.
\item\label{PRE100.ii} There exists canonical isomorphisms of vector spaces
\begin{align*}
\Hom_{\lie g}(V(\mu)\tensor V(\nu)^*,M)&\cong \Hom_{\lie g}(V(\mu),M\tensor V(\nu))\cong (V(\mu)^*\tensor M\tensor V(\nu))^{\lie g}\\
&\cong (V(\nu)^*\tensor M^*\tensor V(\mu))^{\lie g}
\cong \Hom_{\lie g}(V(\nu),M^*\tensor V(\mu)).
\end{align*}
\item\label{PRE100.iii} $M_{\mu-\nu}^\nu=\{ m\in M_{\mu-\nu}\,:\, \Ann_{ U(\lie n^-)} m\supset \Ann_{ U(\lie n^-)} v_\mu\}$.
\item\label{PRE100.iv} The linear map $M_{\mu-\nu}^\nu\to \Hom_{\lie g}(V(\mu)\tensor V(\nu)^*,M)$ given by $m\mapsto \chi_m$,
where
$$
\chi_m( a(v_\mu\tensor \xi_{-\nu}))=a m,\qquad a\in U(\lie g)
$$
is an isomorphism of vector spaces. In particular, all vector spaces in~\eqref{PRE100.ii} are isomorphic to~$M_{\mu-\nu}^\nu$.\qed
\end{enumerate}
\end{prop}

\subsection{}\label{PRE60}
Let~$K=\bigoplus_{x\in J} \bc e_x$ be a semi-simple commutative algebra with primitive pairwise orthogonal idempotents~$e_x$ and let~$V$ an $K$-bimodule. Assume
that~$\dim e_x V e_y<\infty$ for all~$x,y\in J$ and that~$V=\bigoplus_{x,y\in J} e_x V e_y$ (which is always the case if~$J$ is finite).
Let~$T_K^0(V)=K$, $T_K^r(V)$ be the $r$-fold tensor product of~$V$ over~$K$ and set~$T_K(V)=\bigoplus_{r\in\bz_+} T_K^r(V)$. This is a $\bz_+$-graded
associative algebra. In particular, if $A$ is a $\bz_+$-graded associative algebra and~$A[0]$ is commutative semi-simple, we have a canonical 
homomorphism of associative algebras $T_{A[0]}(A[1])\to A$ (cf.~\cite{BGS}).

Let~$\Delta$ be the quiver with~$\Delta_0=J$ and with $\dim e_x V e_y$ arrows $x\from y$
for all~$x,y\in J$. We have an isomorphism of algebras~$K\to \bigoplus_{x\in J} \bc1_x\subset \bc \Delta$. In particular,
we can regard the subspace of $\bc \Delta$ spanned by all arrows as a $K$-bimodule and for any choice of basis 
in~$e_x V e_y$, $x,y\in J$ this subspace is naturally isomorphic to~$V$ as an $K$-bimodule. This isomorphism extends
canonically to an isomorphism of graded associative algebras $T_K(V)\to \bc \Delta$. 
Then, if~$A$ is a quotient of~$T_K(V)$ by an ideal which has the trivial intersection with~$T_K^r(V)$, $r=0,1$, then~$A$ is 
isomorphic to the path algebra $\bc(\Delta,R)$ where~$R$ is the image of~$\ker(T_K(V)\to A)$ in~$\bc\Delta$.

An associative algebra~$A$ is said to be connected if $A=A_1\oplus A_2$ where the~$A_j$ are subalgebras implies that~$A_1=0$ or~$A_2=0$.
Clearly, $\bc(\Delta,R)$ is connected if and only if~$\Delta$ is connected. 

\subsection{}\label{SP100}
Let~$\Psi\subset R^+$ be a fixed extremal set.
\begin{prop}
Let~$F\subset P^+$ be interval closed in the partial order~$\le_\Psi$. Then the algebra $\bt_\Psi^{\lie g}(F)$ is isomorphic, as a $\bz_+$-graded algebra,
to the path algebra of the quiver~$\Delta_\Psi(F)$. In particular, for all~$\lambda,\mu\in F$
$$
|{\Delta_\Psi(F)}(\lambda,\mu)|=\begin{cases}\dim (V(\lambda)^*\tensor T^{d_\Psi(\lambda,\mu)}(\lie g)\tensor V(\mu))^{\lie g},& \lambda\le_\Psi\mu\\
                                        0,&\text{otherwise},
                                       \end{cases}
$$
and if~$F'\subset F$ is interval closed, then~$\Delta_\Psi(F')$ is a convex subquiver of~$\Delta_\Psi(F)$. Furthermore, $\bs_\Psi^{\lie g}(F)$
is isomorphic to the quotient of~$\bc\Delta_\Psi(F)$ by an ideal generated by paths of length~$2$.
\end{prop}
\begin{pf}
By~\cite{CG1}*{Proposition~4.4}, $\bt^{\lie g}_\Psi(F)$ is isomorphic to~$T_{\bt^{\lie g}_\Psi(F)[0]}(\bt^{\lie g}_\Psi(F)[1])$ as a $\bz_+$-graded associative algebra.
Since~$\bt^{\lie g}_\Psi(F)[0]=\bigoplus_{\lambda\in F} \bc 1_\lambda$, it is enough to prove that for all $\lambda,\mu\in F$, the number of
arrows $\lambda\from\mu$ equals $\dim 1_\lambda\bt^{\lie g}_\Psi(F)[1] 1_\mu$. The latter is zero unless $\mu=\lambda+\beta$, $\beta\in \Psi$.
Since by~\propref{PRE100}
$$
1_\lambda\bt^{\lie g}_\Psi(F)[1] 1_{\lambda+\beta}=(V(\lambda)^*\tensor \lie g\tensor V(\lambda+\beta))^{\lie g}\cong \lie g_\beta^\lambda
$$ 
and $\dim\lie g_\beta=1$, it is enough to prove that $\lambda\from\lambda+\beta\in(\Delta_\Psi)_1$ if and only if
$\lie g_\beta^\lambda\not=0$.
Observe first that $\lambda-\varepsilon(\beta)\in P^+$ implies that $\lambda+\beta=\lambda-\varepsilon(\beta)+\varphi(\beta)\in P^+$. 
Since~$\gamma,\gamma+\alpha_i\in R$ implies that
$e_i \lie g_{\gamma}\not=0$, it follows that~$e_i^t \lie g_{\beta}\not=0$ for all~$0\le t\le \varepsilon_i(\beta)$.
Therefore, $\lie g_\beta^\lambda\not=0$ if and only if $\lambda(h_i)\ge \varepsilon_i(\beta)$ for all~$i\in I$. The remaining assertions are straightforward.
\end{pf}

\subsection{}\label{SP30}
For all~$\beta\in\Psi$ and for all $\lambda\from\lambda+\beta\in (\Delta_\Psi)_1$, fix~$0\not=\boldsymbol{a}_{\lambda,\beta}
\in 1_\lambda\bt^{\lie g}_\Psi 1_{\lambda+\beta}=1_\lambda\bs^{\lie g}_\Psi 1_{\lambda+\beta}$.  This choice is unique up to a non-zero scalar.
It follows from~\cite{CG1}*{Proposition~4.4} that the elements $1_\lambda$, $\lambda\in P^+$ and $\boldsymbol{a}_{\lambda,\beta}$, 
$\lambda\from\lambda+\beta\in(\Delta_\Psi)_1$ generate~$\bt_\Psi^{\lie g}$ and~$\bs_\Psi^{\lie g}$. In particular,
for all $\lambda\le_\Psi\mu$ with $d_\Psi(\lambda,\mu)=2$ the set
$$
\{ \boldsymbol{a}_{\lambda,\beta}\boldsymbol{a}_{\lambda+\beta,\beta'}\,:\, \beta,\beta'\in \Psi,\,\mu=\lambda+\beta+\beta',\, 
\lambda\from\lambda+\beta,\lambda+\beta\from\mu\in (\Delta_\Psi)_1\}
$$
is a basis of~$1_\lambda\bt^{\lie g}_\Psi 1_\mu$. By~\propref{PRE100} and~\corref{*}
$$
1_\lambda\bt^{\lie g}_\Psi 1_\mu \cong (T^2(\lie g))_{\mu-\lambda}^{\lambda}=(T^2(\lie n^+_\Psi))_{\mu-\lambda}^{\lambda}.
$$
Let $\Pi_{\lambda}(\beta,\beta')$ be the image of~$\boldsymbol{a}_{\lambda,\beta}\boldsymbol{a}_{\lambda+\beta,\beta'}$ under this isomorphism.
Using~\cite{CG1}*{Lemma~4.2} we obtain the following
\begin{prop}
Let~$\eta\in\Psi+\Psi$, $\lambda,\lambda+\eta\in P^+$ and assume that
$\Delta_\Psi(\lambda,\lambda+\eta)\not=\emptyset$. 
The elements~$\Pi_\lambda(\beta,\beta')$ where~$\beta,\beta'\in\Psi$, $\beta+\beta'=\eta$ and $\lambda\from\lambda+\beta\from\lambda+\eta\in\Delta_\Psi(\lambda,
\lambda+\eta)$ form a basis of~$T^2(\lie n^+_\Psi)^\lambda_\eta$.
In particular, we have a relation 
$$
\sum_{\beta\in\Psi\,:\, \lambda\from\lambda+\beta\from\lambda+\eta\in\Delta_\Psi(\lambda,\lambda+\eta)} x_\beta \boldsymbol{a}_{\lambda,\beta}
\boldsymbol{a}_{\lambda+\eta,\eta-\beta}=0
$$
in~$\bs_\Psi^{\lie g}$ if and only if
\begin{equation*}
\sum_{\beta\in\Psi\,:\, \lambda\from\lambda+\beta\from\lambda+\eta\in\Delta_\Psi(\lambda,\lambda+\eta)} x_\beta \Pi_\lambda(\beta,\eta-\beta)\in\bigwedge^2\lie n^+_\Psi.
\qedhere
\end{equation*}
\end{prop}

\subsection{}\label{SP50}
Thus, to describe the relations, it remains to find a way for describing the elements~$\Pi_{\lambda}(\beta,\beta')$. It turns out that the most
convenient language is provided by $\lie g$-module maps.

Let~$V$ be a finite dimensional $\lie g$-module. Given~$f\in\Hom_{\lie g}(V(\mu),V\tensor V(\lambda))$, note that~$f$
is uniquely determined by~$f(v_\mu)$. 
Using~\propref{PRE100} we obtain an isomorphism of vector spaces
$$
\Hom_{\lie g}(V(\mu),V\tensor V(\lambda))\to V_{\mu-\lambda}^\lambda
$$
given by
$$
f\mapsto v_f:=(1\tensor \xi_{-\lambda})f(v_\mu).
$$
In particular, we have
\begin{equation}\label{SP40.0}
f(v_\mu)=v_f\tensor v_\lambda+ U(\lie n^+)_+ v_f\tensor  U(\lie n^-)_+ v_\lambda.
\end{equation}

Let~$\beta\in\Psi$, $\lambda\in P^+$ and assume that $\lambda\from\lambda+\beta\in(\Delta_\Psi)_1$ and so $\lie g_\beta=\lie g_\beta^\lambda$.
Fix root vectors~$e_\gamma\in\lie g_{\gamma}\setminus\{0\}$, $\gamma\in R^+$. Then by~\eqref{SP40.0}
we have a unique $0\not=p_{\lambda,\beta}\in \Hom_{\lie g}(V(\lambda+\beta),\lie g\tensor V(\lambda))$ satisfying
\begin{equation}\label{SP50.10}
p_{\lambda,\beta}(v_{\lambda+\beta})=e_\beta\tensor v_\lambda+\sum_{\beta<\gamma} e_\gamma\tensor \bou_{\beta,\gamma}(\lambda) v_\lambda,
\end{equation}
where~$\bou_{\beta,\gamma}(\lambda)\in  U(\lie n^-)_{\beta-\gamma}$. Clearly, $p_{\lambda,\beta}(v_{\lambda+\beta})$
spans $(\lie g\tensor V(\lambda))^{\lie n^+}_{\lambda+\beta}$.
Note that the elements $\bou_{\beta,\gamma}(\lambda)$ are uniquely determined modulo~$\Ann_{ U(\lie n^-)} v_\lambda$.

\subsection{}\label{COMP110}
Let~$F(\lie h)$ be the field of fractions of~$S(\lie h)$. Given~$\beta\in\Psi$, let
\begin{align*}
F_\beta(\lie h)&=\{ f g^{-1}\in F(\lie h)\,:\, \lambda\in P^+,\,\lambda\from\lambda+\beta\in(\Delta_\Psi)_1\implies\lambda(g)\not=0\}.
\end{align*}
Clearly, $F_\beta(\lie h)$ is a subring of~$F(\lie h)$.
Given $\lambda\from\lambda+\beta\in (\Delta_\Psi)_1$, note that~$\lambda:S(\lie h)\to \bc$ extends
canonically to a homomorphism $F_\beta(\lie h)\to \bc$ which we also denote by~$\lambda$. 

Furthermore, regard~$ U(\lie b)$ as a right $S(\lie h)$-module via the right multiplication and a left $ U(\lie n^-)$ module via the left multiplication. 
Then $U(\lie b)\tensor_{S(\lie h)} F_\beta(\lie h)$ is a right $S(\lie h)$-module 
and is isomorphic to $ U(\lie n^-)\tensor F_\beta(\lie h)$ as a left~$U(\lie n^-)$-module by the PBW theorem. Thus, $\lambda$
induces a surjective homomorphism of left~$U(\lie n^-)$-modules $\pi_{\lambda,\beta}:U(\lie b)\tensor_{S(\lie h)} F_\beta(\lie h)\to U(\lie n^-)$.

Let~$\lambda\in\lie h^*$. 
The quotient of~$ U(\lie b)$ by the left ideal generated by the kernel of $\lambda:S(\lie h)\to \bc$ is isomorphic to~$ U(\lie n^-)$ as a left $ U(\lie n^-)$-module
and so we have a surjective homomorphism of left $ U(\lie n^-)$-modules $\pi_\lambda:  U(\lie b)\to  U(\lie n^-)$.
Clearly, the restriction of $\pi_\lambda$ to~$ U(\lie n^-)$ is the identity map. 
Furthermore, if~$V$ is a finite dimensional $\lie g$-module, $v\in V_\mu$ and~$x\in  U(\lie b)$, then $x-\pi_\mu(x)\in\Ann_{ U(\lie g)}v$.
\begin{lem} 
Suppose that~$x\in  U(\lie b)$, $y\in  U(\lie b)_{-\eta}$, $\eta\in \bz_+ R^+$. Then
$\pi_\lambda(xy)=\pi_{\lambda-\eta}(x)\pi_\lambda(y)$. Furthermore, if $\lambda\from\lambda+\beta\in (\Delta_\Psi)_1$
then $\pi_{\lambda,\beta}(x\tensor f)=\pi_\lambda(x)\tensor \pi_{\lambda,\beta}(f)$ for all~$f\in F_\beta(\lie h)$.
\end{lem}
\begin{pf}
Note that~$\pi_\lambda(xy)=\pi_\lambda(x)\pi_\lambda(y)$ for all~$x\in  U(\lie n^-)$, $y\in  S(\lie h)$. 
Since~$ U(\lie b)\cong  U(\lie n^-)\tensor  S(\lie h)$ by the PBW theorem, it is enough to show that
$\pi_\lambda(h y)=\pi_{\lambda-\eta}(h)\pi_\lambda(y)$ for all~$h\in\lie h$, $y\in  U(\lie b)_{-\eta}$. We have $\pi_\lambda(hy)=\pi_\lambda(yh)-\eta(h)\pi_\lambda(y)=\pi_\lambda(y)(\lambda-\eta)(h)=\pi_{\lambda-\eta}(h)
\pi_\lambda(y)$. The second assertion is obvious.
\end{pf}
Given~$\beta\in\Psi$, we have a group homomorphism $F_\beta(\lie h)^\times\to (\bc^\times)^{(\Delta_\Psi)_1}$ defined by
$h\mapsto (z_{\lambda,\gamma}(h)\,:\,\gamma\in\Psi,\lambda\from\lambda+\gamma\in(\Delta_\Psi)_1)$, where
$$
z_{\lambda,\gamma}(h)=\begin{cases} \lambda(h),&\gamma=\beta\\
                    1,&\gamma\not=\beta
                   \end{cases}
$$
This yields a natural group homomorphism $\prod_{\beta\in\Psi} F_\beta(\lie h)^\times \to (\bc^\times)^{(\Delta_\Psi)_1}$. We denote its
image by~$G_\Psi$.

\subsection{}\label{SP55}
\begin{defn}
Let~$\beta\in\Psi$. We call a tuple
$$( \bou_{\beta,\gamma}\in U(\lie b)_{\beta-\gamma}\tensor_{S(\lie h)} F_\beta(\lie h)\,:\,\beta\le \gamma, \gamma\in R^+)
$$
an {\em adapted family for~$\beta$} if~$\bou_{\beta,\beta}=1$ and for all~$\lambda\from\lambda+\beta\in(\Delta_\Psi)_1$, the vector
\begin{equation}\label{SP55.10}
\sum_{\gamma\in R^+\,:\,\beta\le\gamma} e_\gamma\tensor \pi_{\lambda,\beta}(\bou_{\beta,\gamma}) v_\lambda
\end{equation}
spans~$(\lie g\tensor V(\lambda))^{\lie n^+}_{\lambda+\beta}$.
\end{defn}
\begin{prop}
Let $\beta\in\Psi$ and suppose that~$(\bou_{\beta,\gamma}\in U(\lie b)\tensor_{S(\lie h)} F_\beta(\lie h)\,:\,\beta\le \gamma, \gamma\in R^+)$
is an adapted family for~$\beta$. Then for all~$\beta'\in\Psi$ and for all~$\lambda\in P^+$ such that~$\lambda\from\lambda+\beta,
\lambda+\beta\from\lambda+\beta+\beta'\in (\Delta_\Psi)_1$ we have, up to a non-zero scalar,
\begin{equation}\label{SP50.20}
\Pi_{\lambda}(\beta',\beta)=e_{\beta}\tensor e_{\beta'}+\sum_{\beta<\gamma\,:\,\gamma,\beta+\beta'-\gamma\in\Psi} e_\gamma\tensor \bou_{\beta,\gamma}(\lambda+\beta') e_{\beta'},
\end{equation}
where
$$
\bou_{\beta,\gamma}(\nu)=\pi_{\nu,\beta}(\bou_{\beta,\gamma})\pmod{\Ann_{U(\lie n^-)}v_\nu},\qquad \nu\from\nu+\beta\in (\Delta_\Psi)_1.
$$
In particular, if~$\beta\in\Psi$ is maximal, $\Pi(\beta',\beta)=e_{\beta}\tensor e_{\beta'}$.
\end{prop}
\begin{pf}
Let~${\nu_1}\le_\Psi{\nu_2}$. Since~$\lie g\cong\lie g^*$,  by~\propref{PRE100}\eqref{PRE100.ii} we have the following canonical isomorphisms of vector spaces
\begin{align*}
1_{\nu_1} \bt^{\lie g}_\Psi 1_{\nu_2}&=(V({\nu_1})^*\tensor T^{d_\Psi({\nu_1},{\nu_2})}(\lie g) \tensor V({\nu_2}))^{\lie g}
\cong (V({\nu_2})^*\tensor T^{d_\Psi({\nu_1},{\nu_2})}(\lie g) \tensor V({\nu_1}))^{\lie g}\\
&\cong \Hom_{\lie g}(V({\nu_2}),T^{d_\Psi({\nu_1},{\nu_2})}(\lie g)\tensor V({\nu_1})).
\end{align*}
Moreover, this isomorphism is compatible with products and compositions, that is, if
$x\in 1_{\nu_1}\bt^{\lie g}_\Psi 1_{\nu_2}$, $y\in 1_{\nu_2}\bt^{\lie g}_\Psi 1_{\nu_3}$, $\nu_1\le_\Psi\nu_2\le_\Psi\nu_3$
and
$$x\mapsto f\in\Hom_{\lie g}(V({\nu_2}),T^{d_\Psi({\nu_1},{\nu_2})}(\lie g)\tensor V({\nu_1})),\quad
y\mapsto g\in\Hom_{\lie g}(V(\nu_3),T^{d_\Psi({\nu_2},\nu_3)}(\lie g)\tensor V({\nu_2})),$$ 
then
$$
x y\mapsto (1\tensor f)\circ g\in\Hom_{\lie g}(V(\nu_3),T^{d_\Psi({\nu_1},\nu)}(\lie g)\tensor V({\nu_1})).
$$
In particular, $\beta\in\Psi$ and
$\lambda\from\lambda+\beta\in(\Delta_\Psi)_1$, we may
assume, without loss of generality, that $p_{\lambda,\beta}\in\Hom_{\lie g}(V(\lambda+\beta),\lie g\tensor V(\lambda))$
is the image of~$\boldsymbol{a}_{\lambda,\beta}$ under the above isomorphism. 
Then~$\Pi_{\lambda}(\beta',\beta)$ is the image of
$$
(1\tensor p_{\lambda,\beta'})\circ p_{\lambda+\beta',\beta}
$$
under the isomorphism 
$$
\Hom_{\lie g}(V(\lambda+\beta+\beta'),T^2(\lie g)\tensor V(\lambda))\to T^2(\lie n^+_\Psi)_{\beta+\beta'}^\lambda.
$$
It is now immediate from~\eqref{SP40.0}, \eqref{SP55.10} and~\propref{PRE100}\eqref{PRE100.0} that
$$
\Pi_{\lambda}(\beta',\beta)=e_{\beta}\tensor e_{\beta'}+\sum_{\gamma\in R^+\,:\,\beta<\gamma} e_\gamma\tensor \bou_{\beta,\gamma}(\lambda+\beta') e_{\beta'}.
$$
Since $\bou_{\beta,\gamma}(\lambda+\beta')e_{\beta'}\in \lie g_{\beta+\beta'-\gamma}$, it follows from~\corref{*} that 
$\bou_{\beta,\gamma}(\lambda+\beta') e_\beta\not=0$ implies that $\gamma,\beta+\beta'-\gamma\in\Psi$.
\end{pf}
The following elementary corollary establishes part~\eqref{mainthm.ii'} of~\thmref{mainthm}.
\begin{cor}
Let~$\beta\in\Psi$, $\lambda\in P^+$. Then $(\lambda\from\lambda+\beta\from\lambda+2\beta)\notin\mathfrak R_\Psi(\lambda,\lambda+2\beta)$.\qed
\end{cor}

\section{First examples}\label{FEX}

The aim of this section is to provide the reader with relatively simple examples of quivers and relations arising from algebras $\bs_\Psi^{\lie g}$, before
we undertake a complete study of all possible relations in these algebras for $\lie g$ of types~$A$ and~$C$. 
We begin with the infinite dimensional example announced in~\cite{CG1} which is independent of type of~$\lie g$.
The same computation allows us to obtain a complete description of relations in~$\bs_\Psi^{\lie g}$ for~$\lie g$ of type~$A_2$. Then we describe the relations in the algebras corresponding
to $\lie g$ of type~$G_2$. The remaining rank~$2$ case is postponed until~\ref{C150}.

Throughout the rest of the paper, given~$\lambda\in P^+$ and~$i\notin I$, we set~$\lambda(h_{i})=+\infty$.

\subsection{}\label{hered}
We begin by excluding the case~$|\Psi|=1$. In this case the algebra~$\bs_\Psi^{\lie g}$ is hereditary and
we have two possibilities. If~$\Psi=\{\theta\}$ then every connected component of
$\Delta_\Psi$ is isomorphic to the quiver~$\mathbb A_\infty^{op}$, where
$$
\mathbb A_\infty = 0\to 1 \to 2 \to \cdots
$$
If~$\Psi=\{\beta\}$ with~$\beta\not=\theta$ then~$\beta\notin P^+$ (in fact, it is easy to check that if the highest short root is
contained in~$\Psi$ then~$|\Psi|>1$) and so the connected components of~$\Delta_\Psi$ are either 
simple one dimensional or of type~$\mathbb A_n$ with the subspace orientation.

\subsection{}\label{type_indep_example}
Suppose that~$\lie g$ is not of type~$A$ or~$C$ (in fact, the computation of the relations works for the type~$A$ as well, but the quiver is more complicated, as we will see below; the corresponding construction for the
type~$C$ will be discussed later). Then there exists a unique~$i_0\in I$ such that~$\theta-\alpha_{i_0}\in R^+$ and it is not hard to see
that~$\Psi=\{\theta,\theta-\alpha_{i_0}\}$ is extremal. 

Recall (cf.~\cite{RinBook}) that a pair~$(\Delta,\tau)$ where~$\Delta=(\Delta_0,\Delta_1)$ is a quiver without multiple arrows and~$\tau:\Delta_0'\to \Delta_0$, $\Delta_0'\subseteq \Delta_0$ is an injective map, is called a translation quiver (and $\tau$ is called the translation map) if  
$(\tau(z))^+=z^-$ for all~$z\in \Delta_0'$. A full embedding of translation quivers $(\Delta,\tau)\to
(\Delta',\tau')$ is
a full embedding of quivers $\Delta\to\Delta'$ which maps the domain of~$\tau$ into the domain of~$\tau'$ and is compatible with the maps~$\tau$, $\tau'$.
If~$(\Delta,\tau)$ is a translation quiver and has no multiple arrows, a relation of the form
$\sum_{y\in x^-} (x\leftarrow y)(y \leftarrow\tau (x))$, $x\in\Delta_0$, is called a mesh relation.

Given a quiver~$\Delta$, a translation quiver~$\bz\Delta$ is defined by 
\begin{gather*}
(\bz\Delta)_0=\bz\times\Delta_0,\quad (\bz\Delta)_1=\{ (n,x)\from (n,y),\, (n+1,y)\from (n,x)\,:\, x\from y\in\Delta_1\},\\
\tau((n,x))=(n-1,x).
\end{gather*}
If~$\Delta$ is a Dynkin quiver, $\bz\Delta$ depends only on~$\bar\Delta$ (cf.~\cite{RinBook}*{\S2.1}).

\begin{prop}
Every connected subalgebra of~$\bs_\Psi^{\lie g}$ is isomorphic to the path algebra of the translation quiver
\begin{equation}\label{SP100.10}
\Gamma=\def\dgeverynode{\scriptscriptstyle}\divide\dgARROWLENGTH by 2\dgVERTPAD=2pt\dgHORIZPAD=2pt
\begin{diagram}
\node[3]{}\node{\vdots}\arrow{s}\node{\vdots}\arrow{s}
\\
\node{}\node{}\node{(0,2)}\arrow{s}\node{(1,2)}\arrow{w}\arrow{s}\node{(2,2)}\arrow{w}\arrow{s}\node{\displaystyle\cdots}\arrow{w}
\\
\node{}\node{(0,1)}\arrow{s}\node{(1,1)}\arrow{w}\arrow{s}\node{(2,1)}\arrow{w}\arrow{s}\node{(3,1)}\arrow{w}\arrow{s}\node{\displaystyle\cdots}\arrow{w}
\\
\node{(0,0)}\node{(1,0)}\arrow{w}\node{(2,0)}\arrow{w}\node{(3,0)}\arrow{w}\node{(4,0)}\arrow{w}\node{\displaystyle\cdots}\arrow{w}
\end{diagram}
\end{equation}
with the translation map~$\tau(m,n)=(m,n+1)$, $m,n\in\bz_+$ and with the mesh relations.
\end{prop}
\begin{pf}
Suppose that~$\lambda\in P^+$ is a sink in~$\Delta_\Psi$. 
Since~$\varphi(\theta)=\varpi_{i_0}$, we must have~$\lambda(h_{i_0})=0$. 
Suppose that $\mu\in P^+$ is a sink in~$\Delta_\Psi[\lambda]$, $\mu\not=\lambda$. Since~$\Delta_\Psi[\lambda]_0\subset (\lambda+\bz\Psi)\cap P^+$,
$\mu=\lambda+m\theta+k\beta$ for some $m,k\in\bz$. Interchanging the role of~$\lambda$ and~$\mu$, if necessary,
we may assume that~$m\ge 0$.
Since~$\theta(h_{i_0})=1$,
$\beta(h_{i_0})=-1$, we have~$m=k>0$. 
On the other hand, for all~$j\not=i_0$ we have $\beta(h_j)=\varphi_j(\beta)$ and so $\mu(h_j)-\lambda(h_j)=
k\varphi_j(\beta)$. Since~$\lambda(h_j)\ge 0$ and $\varphi_{i_0}(\beta)=0$, this implies that $\mu-\varphi(\beta)\in P^+$
which is a contradiction since~$\mu$ is a sink.
Thus, every connected component of~$\Delta_\Psi$ contains a unique sink~$\lambda$ hence
$$
\Delta_\Psi[\lambda]_0\subset (\lambda\le_\Psi)=\{ \lambda+r\theta+s\beta\,:\, 0\le s\le r\}.
$$ 
Note that $\lambda+r\theta+s\beta$, $0\le s\le r$ is connected to~$\lambda$ by a path
$$
\lambda\from\lambda+\theta\from\cdots \from\lambda+r\theta\from \lambda+r\theta+\beta\from\cdots \from\lambda+r\theta+s\beta.
$$
Thus, $\Delta_\Psi[\lambda]_0=(\lambda\le_\Psi)$.
Define a map $\Delta_\Psi[\lambda]_0\to \Gamma_0=\bz_+\times\bz_+$ by $\lambda+r\theta+s\beta\mapsto (r-s,s)$. This map is clearly a bijection.
Furthermore, we have an arrow $\lambda+r\theta+s\beta\from\lambda+(r+1)\theta+s\beta$ and an arrow 
$\lambda+r\theta+s\beta\from\lambda+r\theta+(s+1)\beta$ provided that~$s<r$. 
Since in the quiver~$\Gamma$ we have an arrow $(m,n)\from (m+1,n)$ for all~$m,n\in\bz_+$ and an arrow $(m,n)\from (m-1,n+1)$ for all~$m>0$,
it follows that~$\Delta_\Psi[\lambda]\cong \Gamma$. Finally, if we define~$\tau:\Delta_\Psi[\lambda]_0\to\Delta_\Psi[\lambda]_0$ by
$\tau(\mu)=\mu+\theta+\beta$, we conclude that our isomorphism is in fact an isomorphism of translation quivers.

It remains to compute the relations in our algebra. Since~$\beta<\theta$, by~\propref{SP55} we have $\Pi_{\lambda}(\beta,\theta)=e_{\theta}\tensor e_\beta$. Assuming
that~$[e_{i_0},e_\beta]=e_\theta$ we can easily check that $\bou_{\beta,\beta}=1$ and
$\bou_{\beta,\theta}=-f_{i_0}\tensor h_{i_0}^{-1}\in  U(\lie b)_{-\alpha_{i_0}}\tensor_{S(\lie h)} F_\beta(\lie h)$
form an adapted family for~$\beta$.
Since~$f_{i_0}\notin \Ann_{ U(\lie n^-)} v_\nu$ if~$\nu(h_{i_0})>0$, we conclude that
$\Pi_\lambda(\theta,\beta)=e_{\beta}\tensor e_\theta-(\lambda(h_{i_0})+1)^{-1} e_\theta\tensor e_\beta$. 

Suppose that~$\lambda(h_{i_0})>0$. Then~$t_{\lambda,\theta+\beta}=2$ and, clearly, 
$\lambda(h_{i_0})\Pi_{\lambda}(\beta,\theta)-(\lambda(h_{i_0})+1)\Pi_{\lambda}(\theta,\beta)\in\bigwedge^2 \lie n^+_\Psi$.
Fix the isomorphism 
$\Phi:\bt^{\lie g}_\Psi\to\bc\Delta_\Psi$ by assigning
$$
\boldsymbol{a}_{\lambda,\theta}\mapsto (\lambda\from\lambda+\theta),\qquad \lambda\in P^+
$$
and
$$
\boldsymbol{a}_{\lambda,\beta}\mapsto (-1)^{\lambda(h_{i_0})} (\lambda(h_{i_0}))^{-1} (\lambda\from\lambda+\beta),\qquad \lambda(h_{i_0})>0.
$$ 
Then
it is easy to see that $\mathfrak R_\Psi(\lambda,\lambda+\theta+\beta)$ is spanned by the mesh relation with respect to our translation map. If~$\lambda(h_{i_0})=0$ then
$\Pi_\lambda(\theta,\beta)\in\bigwedge^2 \lie n^+_\Psi$, so the unique path is a zero relation and is again the mesh relation with respect to our translation map.
\end{pf}
Note that we have a full embedding of translation quivers~$\Gamma\hookrightarrow\bz\mathbb A_\infty$ given on the vertices by $(r,s)\mapsto (-r-s,r)$.
The quiver~$\Gamma^{op}$ identifies with the Auslander-Reiten quiver for~$\mathbb A_\infty$ and so a connected subalgebra of $\bs_\Psi^{\lie g}$ can be regarded
as an infinite dimensional analogue of the Auslander algebra of $\bc \mathbb A_\infty$.

\subsection{}\label{A2-1}
In the remainder of the section we will consider $\lie g$ of types~$A_2$ and~$G_2$. Identify~$P$ with~$\bz\times\bz$ and write $(\lambda(h_1),\lambda(h_2))$ for~$\lambda\in P$.

Let~$\lie g$ be of type~$A_2$.
Then~$R^+$ contains two extremal sets with~$|\Psi|>1$, namely $\Psi_i=\{\alpha_i,\theta\}$, $i\in I$. Clearly, it is enough to analyse one of them,
say~$\Psi=\Psi_1$.

Suppose that~$(m,n),(m',n')\in P^+$ are in the same connected component. Then~$(m',n')\in ((m,n)+\bz\Psi)\cap P^+$, that is, $(m',n')=(m+r+2s,n+r-s)$ for some~$r,s\in\bz$. 
This implies that $m'-n'=m-n\pmod 3$.
Since~$\varphi(\alpha_1)=(2,0)$, $\varphi(\theta)=\theta=(1,1)$, the sinks in~$\Delta_\Psi$ are
$(0,m)$, $m\in\bz_+$ and $(1,0)$. 
Let~$0\le r<3$. Then we have
$$
(0,r)\from (1,r+1)\from \cdots \from (2k,2k+r)\to (2(k-1),2k+1+r)\to\cdots \to (0,3k+r),
$$
hence all sinks $(0,3k+r)$ lie in~$\Delta_\Psi[(0,r)]$. Finally, we have $(1,0)\from (2,1)\to (0,2)$ hence~$(1,0)$ belongs to~$\Delta_\Psi[(0,2)]$.
Thus, $\Delta_\Psi$ has three connected components given by
$$
\Delta_\Psi[(0,r)]_0=\{ (m,n)\in\bz_+\times\bz_+\,:\, m-n=r\pmod 3\}
$$
the arrows being $(m,n)\from (m+1,n+1)$ and~$(m,n)\from (m+2,n-1)$, $n>0$. The translation structure is given by~$\tau(m,n)=(m+3,n)$. The computation of
relations performed in~\ref{type_indep_example} implies that all relations are the mesh relations. 

It is easy to see that the quivers~$\Delta_\Psi[(0,r)]$, $r\in\{0,1,2\}$, and hence the corresponding connected subalgebras of~$\bs_\Psi^{\lie g}$, are not isomorphic.
For that, note that~$\Delta_\Psi[(0,r)]$ has a unique sink $\lambda_r$ 
such that~$|\lambda^-|=1$ (indeed, clearly $\lambda_0=(0,0)$, $\lambda_1=(0,1)$ and~$\lambda_2=(1,0)$ have this property).
It follows that any full map of quivers $\Delta_\Psi[(0,r)]\to\Delta_\Psi[(0,s)]$ must send $\lambda_r$ to~$\lambda_s$ and~$\lambda^-_r$ to~$\lambda^-_s$.
On the other hand, $\lambda_r$ belongs to the following full connected subquivers of~$\Delta_\Psi[(0,r)]$, respectively
$$
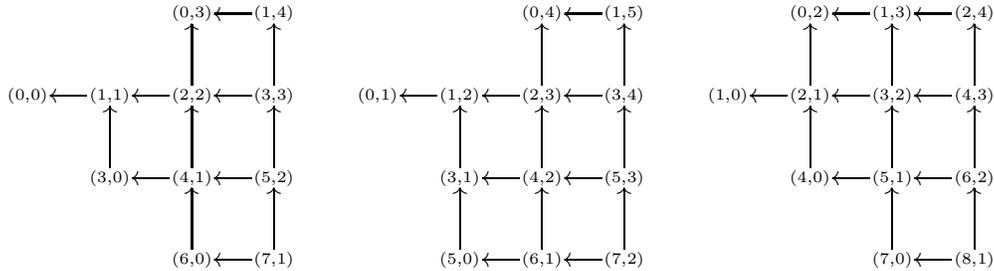
\def\dgeverynode{\scriptscriptstyle}
\divide\dgARROWLENGTH by 2\dgVERTPAD=2pt\dgHORIZPAD=2pt
\begin{diagram}
\node[2]{}\node{(0,3)}\node{(1,4)}\arrow{w}\\
 \node{(0,0)}\node{(1,1)}\arrow{w}\node{(2,2)}\arrow{n}\arrow{w}\node{(3,3)}\arrow{w}\arrow{n}\\
\node{}\node{(3,0)}\arrow{n}\node{(4,1)}\arrow{n}\arrow{w}\node{(5,2)}\arrow{w}\arrow{n}\\
\node[2]{}\node{(6,0)}\arrow{n}\node{(7,1)}\arrow{w}\arrow{n}
\end{diagram}
\qquad
\begin{diagram}
\node[2]{}\node{(0,4)}\node{(1,5)}\arrow{w}
\\
\node{(0,1)}\node{(1,2)}\arrow{w}\node{(2,3)}\arrow{w}\arrow{n}\node{(3,4)}\arrow{w}\arrow{n}\\
\node{}\node{(3,1)}\arrow{n}\node{(4,2)}\arrow{n}\arrow{w}\node{(5,3)}\arrow{n}\arrow{w}\\
\node{}\node{(5,0)}\arrow{n}\node{(6,1)}\arrow{n}\arrow{w}\node{(7,2)}\arrow{w}\arrow{n}
\end{diagram}
\qquad
\begin{diagram}
\node{}\node{(0,2)}\node{(1,3)}\arrow{w}\node{(2,4)}\arrow{w}\\
\node{(1,0)}\node{(2,1)}\arrow{n}\arrow{w}\node{(3,2)}\arrow{w}\arrow{n}\node{(4,3)}\arrow{w}\arrow{n}\\
\node{}\node{(4,0)}\arrow{n}\node{(5,1)}\arrow{w}\arrow{n}\node{(6,2)}\arrow{w}\arrow{n}\\
\node{}\node{}\node{(7,0)}\arrow{n}\node{(8,1)}\arrow{w}\arrow{n}
\end{diagram}
$$
These quivers are obviously non-isomorphic.

\subsection{}\label{G2-10}
Let~$\lie g$ be of type~$G_2$. Let~$\alpha_1$ (respectively, $\alpha_2$) be the long (respectively, the short) simple root. 
Then $R^+=\{\alpha_1,\alpha_2,\alpha_1+\alpha_2,\alpha_1+2\alpha_2,\alpha_1+3\alpha_2,\theta=2\alpha_1+3\alpha_2\}$. It is not hard to show that there
are only two extremal sets of positive roots containing more than one element, namely $\Psi_1=\{\theta-\alpha_1,\theta\}$ and~$\Psi_2=\{\alpha_1,\theta\}$,
which correspond to the two one-dimensional faces of the convex hull of~$R$ having trivial intersection with~$-R^+$.
The set~$\Psi_1$ has already been considered in~\propref{type_indep_example}. We should only note that
since~$\varphi(\theta)=(1,0)$ and~$\varphi(\theta-\alpha_1)=(0,3)$, $(0,r)$, $0\le r<3$ are the only sinks in~$\Delta_\Psi$ and
hence by~\propref{SP100} $\Delta_\Psi$ has three isomorphic connected components.

The situation is rather different if~$\Psi=\Psi_2$. Since~$\varphi(\alpha_1)=(2,0)$ and~$\varphi(\theta)=(1,0)$, it follows that
$(m,n)$ is a sink in~$\Delta_\Psi$ if and only if~$m=0$. Furthermore, since~$\varepsilon(\alpha_1)=(0,3)$, 
we have an arrow $(m,n)\from (m+2,n-3)$ if and only if~$n\ge 3$. Suppose that we have two sinks~$(0,x)$, $(0,y)$ in
the same connected component of~$\Delta_\Psi$. Then we must have $(m+2n,x-3n)=(0,y)$ for some~$m,n\in\bz$, hence
$x=y\pmod 3$.
Furthermore, let~$0\le r\le 2$. Then we have in~$\Delta_\Psi$
$$
(0,r)\from (1,r)\from \cdots \from (2n,r)\to (2(n-1),r+3)\to \cdots\to (2,3(n-1)+r)\to (0,3n+r).
$$ 
Thus, every sink~$(0,3n+r)$, $n\in\bz_+$ lies in~$\Delta_\Psi[(0,r)]$. Therefore, $\Delta_\Psi$ has three isomorphic connected components 
and the quiver~$\Delta_\Psi[(0,r)]$ is
$$
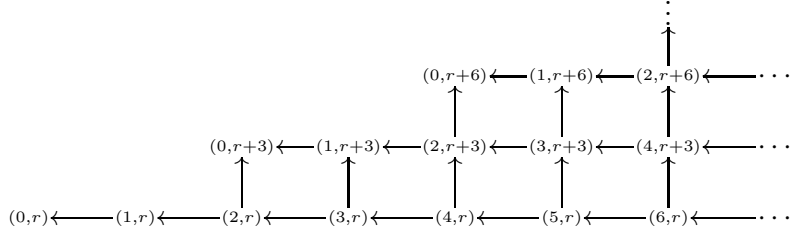
\def\dgeverynode{\scriptscriptstyle}
\divide\dgARROWLENGTH by 2\dgVERTPAD=2pt\dgHORIZPAD=2pt
\begin{diagram}
\node[6]{}\node{\vdots}
\\
\node[4]{}\node{(0,r+6)}\node{(1,r+6)}\arrow{w}\node{(2,r+6)}\arrow{n}\arrow{w}\node{\displaystyle\cdots}\arrow{w} \\
\node{}\node{}\node{(0,r+3)}\node{(1,r+3)}\arrow{w}\node{(2,r+3)}\arrow{n}\arrow{w}\node{(3,r+3)}\arrow{w}\arrow{n}\node{(4,r+3)}\arrow{w}\arrow{n}
\node{\displaystyle\cdots}\arrow{w}\\
\node{(0,r)}\node{(1,r)}\arrow{w}\node{(2,r)}\arrow{n}\arrow{w}\node{(3,r)}\arrow{w}\arrow{n}\node{(4,r)}\arrow{w}\arrow{n}
\node{(5,r)}\arrow{w}\arrow{n}\node{(6,r)}\arrow{w}\arrow{n}\node{\displaystyle\cdots}\arrow{w}
\end{diagram}
$$ 
that is, $\Delta_\Psi[(0,r)]_0=\{ (m,3n+r)\,:\, m,n\in\bz_+\}$ and the arrows are $(m,3n+r)\from (m+1,3n+r)$, $m,n\in\bz_+$,
$(m,3n+r)\from (m+2,3(n-1)+r)$, $n>0$. This is, of course, a translation quiver with~$\tau((m,3k+r))=(m+3,3(k-1)+r)$, $m\in\bz_+$, $k>0$.
In particular, $\Psi$ is our first example of a regular extremal set.
Clearly, there is a full embedding of
$\Delta_\Psi[(0,r)]$ into any of the infinite connected quivers considered in~\ref{A2-1}.

\subsection{}\label{G2-20}
It remains to describe the relations. We write~$x^{(p)}=x^p/p!\in  U(\lie g)$, $x\in\lie g$, $p\in\bz_+$.
Fix root vectors in~$\lie g$ so that $e_{\alpha_1+p\alpha_2}=(\ad e_2)^{(p)} e_1$, $1\le p\le 3$ and $[e_1, e_{\alpha_1+3\alpha_2}]=e_\theta$.
We have only one non-trivial case to consider, namely $\eta=\theta+\alpha_1=(3,3)$. If~$\Delta_\Psi( (m,n),(m+3,n-3))$
is non-empty it always contains two paths. \propref{SP55} immediately implies that
$$
\Pi_\lambda(\alpha_1,\theta)=e_\theta\tensor e_{\alpha_1}.
$$
To find~$\Pi_\lambda(\theta,\alpha_1)$, note that $\{\gamma\in R^+\,:\, \alpha_1\le \gamma\}=\{ \alpha_1,\alpha_1+\alpha_2,\alpha_1+2\alpha_2,\alpha_1+3\alpha_2,\theta\}$.
Since
$\dim U(\lie n^-)_{-(\alpha_1+3\alpha_2)}=4$, the monomials 
$$
f_2^{(a)} f_1^{} f_2^{(3-a)},\qquad 0\le a\le 3
$$
which are of course all possible monomials in the~$f_i$ of weight~$-\alpha_1-3\alpha_2$, 
form a basis of~$ U(\lie n^-)_{-(\alpha_1+3\alpha_2)}$. It is not hard to see that the element
\begin{align*}
U&=f_1^{} f_2^{(3)} (h^{}_2+1)h^{}_2(h^{}_2-1)-f_2^{} f_1^{} f_2^{(2)} (h^{}_2+1)h^{}_2(h^{}_2-2)\\&+f_2^{(2)} f_1^{} f_2^{}(h^{}_2+1)(h^{}_2-1)(h^{}_2-2)-
f_2^{(3)} f_1^{} h^{}_2 (h^{}_2-1)(h^{}_2-2)\in  U(\lie b)
\end{align*}
satisfies
$$
e_1^{} U=6 f_2^{(3)} (h_1^{}+h_2^{}+1)+ U(\lie g)\lie n^+,\qquad e_2^{} U\in U(\lie g)\lie n^+.
$$
Define $\bou_{\alpha_1,\gamma}\in  U(\lie b)_{\alpha_1-\gamma}\tensor_{S(\lie h)} F_{\alpha_1}(\lie h)$, $\alpha_1\le\gamma$ as
\begin{align*}
&\bou_{\alpha_1,\alpha_1+p\alpha_2}=(-1)^p f_2^p\tensor \prod_{t=0}^{p-1} (h_2-t)^{-1},\qquad\qquad \qquad 0\le p\le 3,\\
&\bou_{\alpha_1,\theta}=U\tensor ((h_1+h_2+1)h_2(h_2-1)(h_2-2))^{-1}.
\end{align*}
Then it is easy to see that~$(\bou_{\alpha_1,\gamma}\,:\,\alpha_1\le \gamma)$ is an adapted family for~$\alpha_1$, hence
by~\propref{SP55}
$$
\Pi_\lambda(\theta,\alpha_1)=e_{\alpha_1}\tensor e_\theta+((\lambda_1+\lambda_2+2)\lambda_2(\lambda_2-1)(\lambda_2-2))^{-1} e_\theta\tensor \pi_{\lambda+\theta}(U)e_\theta.
$$
Clearly, $\pi_{\lambda+\theta}(U)=-\lambda_2(\lambda_2-1)(\lambda_2-2) f_2^{(3)} f_1^{}+\Ann_{ U(\lie g)} e_\theta$. Since~$(\lambda+\theta)(h_1)>0$,
we conclude using finite dimensional~$\lie{sl}_2$ theory that $f_2^{(3)} f_1^{}\notin\Ann_{ U(\lie n^-)}v_{\lambda+\theta}$. Thus, we get
$$
\Pi_{\lambda}(\theta,\alpha_1)=e_{\alpha_1}\tensor e_{\theta}-(\lambda_1+\lambda_2+2)^{-1} e_{\theta}\tensor e_{\alpha_1}.
$$
It follows that none of the two paths is a relation and that the relations can be chosen to be the mesh relations.

In particular, we see that although $\Psi_2$ is conjugate to~$\Psi_1$ by the action of the Weyl group of~$\lie g$,
the algebras $\bs_{\Psi_j}^{\lie g}$, $j=1,2$ are not isomorphic.

\section{A recursive family of elements in~\texorpdfstring{$ U(\lie b)$}{\bu(b)}}\label{RECF}
In this section we construct a family of elements of~$ U(\lie b)$ which will play the crucial role in constructing adapted families for $\lie g$
of type~$A$ and~$C$.

\subsection{}\label{A45}
Suppose that~$\lie g$ is of type~$A_\ell$. After~\cite{Lit}, the monomial 
$$
f_1^{a_{1,1}}(f_2^{a_{2,2}} f_1^{a_{2,1}})\cdots (f_\ell^{a_{\ell,\ell}}\cdots f_1^{a_{\ell,1}})\in  U(\lie n^-)
$$
where $a_{j,i+1}\ge a_{j,i}$ for all~$1\le j\le \ell$, $1\le i\le j-1$, is called standard.
Furthermore, let~$\lambda\in P^+$. A standard monomial that satisfies
\begin{equation}\label{A45.10}
\lambda(h_i)\ge a_{j,i}-a_{j,i-1}+\sum_{r=j+1}^\ell (2 a_{r,i}-a_{r,i-1}-a_{r,i+1}), \qquad 1\le j\le \ell,\quad 1\le i\le j
\end{equation}
is called $\lambda$-standard (\cite{Lit}*{Definition~22}). In the above we adopt the convention that~$a_{j,s}=0$ if $s<0$ or~$s>j$.
We have the following
\begin{thm}[\cite{Lit}*{Theorems~17 and~25}]\label{thmlit}
Standard monomials form a basis of~$ U(\lie n^-)$. Moreover, for all~$\lambda\in P^+$, the set 
$$
\{Fv_\lambda\,:\, \text{$F$ is a $\lambda$-standard monomial}\}
$$
is a basis of~$V(\lambda)$.\qed
\end{thm}

Assume now that~$\lie g$ be a simple Lie algebra of rank~$\ell$, $J=\{i,i+1,\dots,j\}$, $1\le i\le j\le \ell$. Suppose that the Lie subalgebra $\lie g_J$ of
$\lie g$ generated by the $e_r, f_r$, $r\in J$ is of type~$A_{j-i+1}$. 
Let $\mu\in P^+$ and let $\eta=\sum_{r\in J} k_r \alpha_r$,
$k_r\in\bz_+$. Set
\begin{multline*}
\cal J(\eta)=\{ \boldsymbol a=(a_{s,r})_{i\le s\le j,\,i\le r\le s}\,:\,  a_{s,r+1}\ge a_{s,r},\,i+1\le s\le j,\,i\le r\le s-1,\\ \sum_{s=r}^j a_{s,r}=k_r,\,i\le r\le j\}.
\end{multline*}
and
$$
\cal J(\eta,\mu)=\{ \boldsymbol a\in J(\eta)\,:\, \mu(h_k)\ge a_{s,k}-a_{s,k-1}+\sum_{r=s+1}^j (2 a_{r,k}-a_{r,k-1}-a_{r,k+1}),\, i\le s,k\le j\},
$$
where we assume that~$a_{s,k}=0$ if $k<i$ or~$k>s$.
Using~\thmref{thmlit}, we immediately obtain
\begin{prop}
The monomials
$$
f_{i}^{a_{i,i}} (f_{i+1}^{a_{i+1,i+1}}f_{i}^{a_{i+1,i}})\cdots (f_{j}^{a_{j,j}}\cdots f_{i}^{a_{j,i}}),\qquad \boldsymbol{a}=(a_{s,k})_{i\le k\le s\le j}\in \cal J(\eta)
$$
form a basis of~$ U(\lie n^-)_{-\eta}$ and the vectors
$$
f_{i}^{a_{i,i}} (f_{i+1}^{a_{i+1,i+1}}f_{i}^{a_{i+1,i}})\cdots (f_{j}^{a_{j,j}}\cdots f_{i}^{a_{j,i}})v_\mu,\qquad \boldsymbol{a}=(a_{s,k})_{i\le k\le s\le j}\in\cal J(\eta,\mu)
$$
form a basis of~$V(\mu)_{\mu-\eta}$. In particular, if~$\mu(h_i)>0$ {\rm(}respectively, $\mu(h_j)>0${\rm)} then $f_{j}\cdots f_i v_\mu\not=0$ {\rm(}respectively,
$f_i\cdots f_j v_\mu\not=0${\rm)}.\qed
\end{prop}
\begin{rem}
The last assertion can of course be established by a simple induction from the elementary theory of finite dimensional $\lie{sl}_2$-modules.
\end{rem}

\subsection{}\label{A48}
Let~$J\subset I$ and assume that~$\lie g_J$ is of type~$A_{|J|}$.
Let~$\Sigma(i,j)$, $i\le j\in J$ be the set of all bijective maps $\sigma:\{i,i+1,\dots,j\}\to \{1,\dots,j-i+1\}$ 
satisfying
$$
\sigma(r+1)<\sigma(r)\,\implies\, \sigma(r+1)=\sigma(r)-1,\qquad i\le r<j.
$$
Given~$\sigma\in\Sigma(i,j)$, let~$\bof_\sigma=f_{\sigma^{-1}(1)}\cdots f_{\sigma^{-1}(j-i-1)}$. Let~$\alpha_{i,j}=\sum_{r=i}^j \alpha_r\in R^+$.
\begin{lem}
The set~$\{\bof_\sigma\,:\, \sigma\in\Sigma(i,j)\}$ is a basis of~$ U(\lie n^-)_{-\alpha_{i,j}}$. 
\end{lem}
\begin{pf}
Clearly, if~$\sigma\in\Sigma(i,j)$ then~$\bof_\sigma$ is a standard monomial, and 
if~$\sigma\not=\sigma'$ then the monomials $\bof_\sigma$, $\bof_{\sigma'}$ are distinct. 
Now, we prove by induction on~$j-i$ that every standard monomial of weight
$-\alpha_{i,j}$ is of the form~$\bof_\sigma$. If~$j=i$ there is nothing to prove.
If~$j>i$, let~$F$ be a standard monomial of weight~$-\alpha_{i,j}$. Removing
$f_j$ from~$F$ we obtain a standard monomial of weight~$-\alpha_{i,j-1}$
which is equal to~$\bof_\tau$, $\tau\in \Sigma(i,j-1)$ by the induction hypothesis. Now, since
$F$ is standard and every $f_r$, $i\le r\le j$ occurs in~$F$ exactly once, it follows that
either $f_j$ occurs in the $(j-i+1)$th position or $f_{j-1}$ occurs immediately after~$f_j$. In the first case,
set~$\sigma(r)=\tau(r)$, $r<j$, $\sigma(j)=j-i+1$.
In the second case, set for all~$r<j$
$$
\sigma(r)=\begin{cases}
           \tau(r),& \tau(r)<\tau(j-1)\\
	   \tau(r)+1,&\tau(r)\ge\tau(j-1)
          \end{cases}
$$
and let~$\sigma(j)=\tau(j-1)$. Then it is easy to see that~$F=\bof_\sigma$ and $\sigma\in\Sigma(i,j)$. 
\end{pf}

\subsection{}\label{A50}
Given~$\eta\in \lie h^*$, the assignment $h\mapsto h-\eta(h)$, $h\in\lie h$, $x\mapsto x$, $x\in U(\lie n^-)$ extends to
an algebra automorphism~$\psi_\eta:U(\lie b)\to U(\lie b)$.
Clearly, $\psi_{\eta}\psi_{\eta'}=\psi_{\eta+\eta'}$ for all~$\eta,\eta'\in \lie h^*$ and
$$
x y=\psi_{\eta}(y) x, \qquad \forall\, x\in  U(\lie g)_\eta,\, y\in  S(\lie h).
$$
Observe also that $\pi_\lambda\circ \psi_\eta=\pi_{\lambda-\eta}$, $\lambda,\eta\in\lie h^*$.

Given~$r\le s\in J$ and $\lambda\in\lie h^*$, set
$$
\cal H_{r,s}:=h_r+\cdots+h_s+s-r\in  S(\lie h).
$$
We use the convention that~$\cal H_{r,s}=0$ if~$r>s$.
Note that $\lambda(\cal H_{r,s})\in\bz_+$ for all~$\lambda\in P^+$ and $\lambda(\cal H_{r,s})=0$, $r\le s$ if and only if~$r=s$ and~$\lambda(h_s)=0$.
Define
$$
{\cal X}^\pm_{i,j,k}:=\sum_{\sigma\in\Sigma(i,j)} \bof_\sigma c_{\sigma}^\pm(k),\qquad {\cal X}^-_{i,j}:=\cal X^-_{i,j,j},\,\cal X^+_{i,j}:=\cal X^+_{i,j,i},
$$
where $c_\sigma^\pm\in  S(\lie h)$ are given by the following formulae
\begin{alignat*}{2}
&c_{\sigma}^-(k)=\prod_{s=i+1}^{j} (-1)^{\delta_{\sigma(s),\sigma(s-1)-1}}
(\cal H_{s,k}+1-\delta_{\sigma(s),\sigma(s-1)-1}),&\qquad &j\le k\in J\\
&c_{\sigma}^+(l)=\prod_{r=i}^{j-1} (-1)^{1+\delta_{\sigma(r+1),\sigma(r)-1}}(\cal H_{l,r}+\delta_{\sigma(r+1),\sigma(r)-1}),&&l\le i\in J.
\end{alignat*}
We let~$\cal X^\pm_{j+1,j,k}=1$, $1\le j\le\ell-1$, $\cal X^\pm_{i,j,k}=0$, $i>j+1$.
\begin{lem}
\begin{subequations}
Let~$i\le j\in J$, $\eta\in P$ and $r\in I$. Then 
\begin{align}
&e_r \psi_{\eta}(\cal X^-_{i,j})=\psi_\eta\Big(\delta_{r,i}\cal X^-_{i+1,j}\cal H_{i,j}
+\cal X^-_{r+1,j}\cal X^-_{i,r-1,j} \eta(h_r)\Big)
+\psi_{\eta+\alpha_r}(\cal X^-_{i,j})e_r,
\label{A50.0a}\\
&e_r \psi_\eta({\cal X}^+_{i,j})=\psi_{\eta}\Big(\delta_{r,j} \cal X^+_{i,j-1}\cal H_{i,j}+\cal X^+_{i,r-1}\cal X^+_{r+1,j,i}\eta(h_r)\Big)
+\psi_{\eta+\alpha_r}(\cal X^+_{i,j})e_r.
\label{A50.0b}\\
\intertext{In particular,}
&e_r \cal X^-_{i,j}=\delta_{r,i}\cal X^-_{i+1,j}\cal H_{i,j}+\psi_{\alpha_r}(\cal X^-_{i,j})e_r,
\label{A50.0c}\\
&e_s {\cal X}^+_{i,j}=\delta_{s,j} \cal X^+_{i,j-1} \cal H_{i,j}+\psi_{\alpha_s}(\cal X^+_{i,j})e_s.
\label{A50.0d}
\end{align}
\end{subequations}
\end{lem}
\begin{pf}
We only establish~\eqref{A50.0a}, the proof of~\eqref{A50.0b} being similar. 
The argument is by induction on~$j-i$.
Note that the induction begins since~$\cal X^-_{i,i}=f_i$ and so
$$
 e_r \psi_{\eta}(\cal X^-_{i,i})=\cal X^-_{i,i}e_r+
\delta_{r,i} h_i =\psi_{\eta+\alpha_r}(\cal X^-_{i,i}) e_r +\delta_{r,i}\psi_\eta(\cal X^-_{i+1,i} (\cal H_{i,i}+\eta(h_i))).
$$
We claim that the~$\cal X^-_{p,q,k}$, $p<q\le k\in J$ satisfy 
\begin{equation}\label{A50.10}
\cal X^-_{p,q,k}=f_p \cal X^-_{p+1,q,k}(\cal H_{p+1,k}+1)-
\cal X^-_{p+1,q,k}f_p\cal H_{p+1,k}.
\end{equation}
Indeed, since~$f_p$ commutes with the~$f_t$, $p+1<t$, $t\in J$, a standard monomial $F$ of weight $-\alpha_{p,q}$ 
equals either $f_p\bof_\tau$ or $\bof_\tau f_p$, $\tau\in \Sigma(p+1,q)$.  Let $\sigma,\sigma'\in\Sigma(p,q)$ be the elements corresponding, respectively,
to $f_p\bof_\tau$ and~$\bof_\tau f_p$. Then
$$
c_{\sigma}^-(k)=(\cal H_{p+1,k}+1)c_\tau^-(k),\qquad c_{\sigma'}^-(k)=-\cal H_{p+1,k} c_\tau^-(k).
$$

Note that $f_i$ commutes with~$\cal X^-_{p,j,k}$, $p>i+1$. If~$r\not=i$, we immediately obtain 
from~\eqref{A50.10}, the induction hypothesis and the properties of~$\psi$ that
\begin{multline*}
e_r \psi_{\eta}(\cal X^-_{i,j})=\psi_{\eta}\Big(f_i \cal X^-_{r+1,j}\cal X^-_{i+1,r-1,j}(\cal H_{i+1,j}+1)-
\cal X^-_{r+1,j}\cal X^-_{i+1,r-1,j}f_i\cal H_{i+1,j}\Big)\eta(h_r)\\
+\delta_{r,i+1} \psi_\eta\Big(f_i\cal X^-_{i+2,j}\cal H_{i+1,j}(\cal H_{i+1,j}+1)-
\cal X^-_{i+2,j}\cal H_{i+1,j}f_i\cal H_{i+1,j}\Big)+\psi_{\eta+\alpha_r}(\cal X^-_{i,j})e_r\\=
\psi_{\eta}(\cal X^-_{r+1,j}\cal X^-_{i,r-1,j})\eta(h_r)+\psi_{\eta+\alpha_r}(\cal X^-_{i,j})e_r.
\end{multline*}
Suppose now that~$r=i$. Then we obtain from~\eqref{A50.10} and the induction hypothesis
\begin{align*}
e_i \psi_{\eta}(\cal X^-_{i,j})&=h_i \psi_\eta(\cal X^-_{i+1,j}(\cal H_{i+1,j}+1))-
\psi_{\eta}(\cal X^-_{i+1,j}\cal H_{i+1,j})+\psi_{\eta+\alpha_i}(\cal X^-_{i,j})e_i\\
&=\psi_{\eta}( \cal X^-_{i+1,j}((h_i+\eta(h_i)+1)(\cal H_{i+1,j}+1)-(h_i+\eta(h_i))\cal H_{i+1,j}))+\psi_{\eta+\alpha_i}(\cal X^-_{i,j})e_i\\
&=\psi_{\eta}(\cal X^-_{i+1,j}(\cal H_{i,j}+\eta(h_i)))+\psi_{\eta+\alpha_i}(\cal X^-_{i,j})e_i.\qedhere
\end{align*}
\end{pf}

\section{Type~\texorpdfstring{$A_\ell$, $\ell>1$}{A\_l}}\label{A}

\subsection{}\label{A10}
We have $R^+=\{\alpha_{i,j}\,:\, 1\le i\le j\le \ell\}$.
In particular, $\theta=\alpha_{1,\ell}$.
In terms of fundamental weights,
$\alpha_{i,j}=\varpi_i+\varpi_j-\varpi_{i-1}-\varpi_{j+1}$, where we set $\varpi_0=\varpi_{\ell+1}=0$. 
Since~$\varepsilon(\alpha_{i,j})=\varpi_{i-1}+\varpi_{j+1}$, we immediately obtain
\begin{lem}
Let~$\alpha_{i,j}\in \Psi$, $\lambda\in P^+$.
Then $\lambda\from\lambda+\alpha_{i,j}\in(\Delta_\Psi)_1$ if and only if~$\lambda(h_{i-1}),\lambda(h_{j+1})>0$. 
\end{lem}

\subsection{}\label{A30}
We now proceed to describe the set of paths of length~$2$ in~$\Delta_\Psi$. Suppose that
$\alpha_{i,m},\alpha_{j,k}\in\Psi$, $i\le j$. If~$m+1<j$ we have $\alpha_{i,m}+\alpha_{j,k}=\alpha_{i,k}-\alpha_{m+1,j-1}$ which is a contradiction by~\corref{*},
while if~$j=m+1$, $\alpha_{i,m}+\alpha_{j,k}=\alpha_{i,k}\in R^+$ which is again a contradiction. Thus, we must have~$j\le m$. If~$j=i$ or~$m=k$, there is only
one way of writing~$\alpha_{i,m}+\alpha_{j,k}$ as a sum of roots.
Otherwise we may assume without loss of generality that~$i<j\le k<m$ and so we have 
$$
\alpha_{i,m}+\alpha_{j,k}=\alpha_{i,k}+\alpha_{j,m}
$$
It is easy to check that the sets
\begin{alignat*}{2}
&\{ \alpha_{i,j},\alpha_{i,k}\},&\qquad &i\le j<k\\
&\{ \alpha_{i,k},\alpha_{j,k}\},&\qquad &i<j\le k\\
&\{ \alpha_{i,m},\alpha_{j,k},\alpha_{i,k},\alpha_{j,m}\},&\qquad &i<j\le k<m
\end{alignat*}
are extremal and so all cases listed above actually occur.

Now we can list all paths of length~$2$ in~$\Delta_\Psi$. First, let $\eta=\alpha_{i,j}+\alpha_{i,k}$, $i\le j<k$. 
Suppose that~$\lambda+\eta\in P^+$. Then~$\lambda(h_{i-1})>1$ and either $\lambda(h_{j+1}),\lambda(h_{k+1})>0$
or $\lambda(h_{k+1})>0$, $j=k+1$ and~$\lambda(h_k)=\lambda(h_{j+1})=0$. 
Using~\lemref{A10} we see that
${\Delta_\Psi}(\lambda,\lambda+\eta)$ is non-empty only if~$\lambda(h_{i-1})>1$, $\lambda(h_{k+1})>0$ and we have
\begin{equation}\label{A30.2p.l}
{\Delta_\Psi}(\lambda,\lambda+\eta)=\begin{cases}
\{ \lambda\from\lambda+\alpha_{i,j}\from\lambda+\eta,\lambda\from\lambda+\alpha_{i,k}\from\lambda+\eta\},& \lambda(h_{j+1})>0\\
\{ \lambda\from\lambda+\alpha_{i,k}\from\lambda+\eta\},& k=j+1,\,\lambda(h_{j+1})=0.
\end{cases}
\end{equation}
Similarly, if~$\eta=\alpha_{i,k}+\alpha_{j,k}$, $i<j\le k$, then~${\Delta_\Psi}(\lambda,\lambda+\eta)$ is non-empty only 
if~$\lambda(h_{k+1})>1$, $\lambda(h_{i-1})>0$ and
\begin{equation}\label{A30.2p.r}
{\Delta_\Psi}(\lambda,\lambda+\eta)=\begin{cases}
\{ \lambda\from\lambda+\alpha_{i,k}\from\lambda+\eta,\lambda\from\lambda+\alpha_{j,k}\from\lambda+\eta\},& \lambda(h_{j-1})>0\\
\{ \lambda\from\lambda+\alpha_{i,k}\from\lambda+\eta\},& i=j-1,\,\lambda(h_i)=0.
\end{cases}
\end{equation}
Finally, let~$\eta=\alpha_{i,m}+\alpha_{j,k}=\alpha_{i,k}+\alpha_{j,m}$, $i<j\le k<m$. If~$\lambda+\eta\in P^+$, we must have
$\lambda(h_{i-1}),\lambda(h_{m+1})>0$. Using~\lemref{A10} again we see that~${\Delta_\Psi}(\lambda,\lambda+\eta)$, if non-empty,
has one of the following forms.
\begin{alignat}{2}\label{A30.4p}
&\{\lambda\from\lambda+\alpha_{r,s}\from\lambda+\eta\,:\, (r,s)\in\{ (i,m),(i,k),(j,m),(j,k)\} \},\quad\lambda(h_{j-1}),\lambda(h_{k+1})>0&\\
&\{ \lambda\from\lambda+\alpha_{i,m}\from\lambda+\eta,\lambda\from\lambda+\alpha_{i,k}\from\lambda\},\qquad i=j-1,\lambda(h_{j-1})=0,\lambda(h_{k+1})>0&\\
&\{ \lambda\from\lambda+\alpha_{i,m}\from\lambda+\eta,\lambda\from\lambda+\alpha_{j,m}\from\lambda\}, \qquad k=m-1,\lambda(h_{j-1})>0,\lambda(h_{k+1})=0&\\
&\{\lambda\from\lambda+\alpha_{i,m}\from\lambda+\eta\},\phantom{\lambda\from\lambda}\qquad\qquad i=j-1,k=m-1,\lambda(h_{j-1})=\lambda(h_{k+1})=0.&
\end{alignat}
In particular, we have the following
\begin{lem}
An extremal set $\Psi\subset R^+$ is regular if and only if $\alpha_{i,j},\alpha_{i,k}\in \Psi$, $j<k$,
{\em(}respectively, $\alpha_{i,k},\alpha_{j,k}\in\Psi$, $i<j${\em)}, implies that~$k>j+1$ {\em(}respectively, $j>i+1${\em)}.
\end{lem}

\subsection{}\label{A32}\label{TQF10}
Fix~$r,s>0$. Given $\boldsymbol{m}=(m_1,\dots,m_r)\in(\bz_+\cup\{+\infty\})^r$,
$\boldsymbol{n}=(n_1,\dots,n_s)\in
(\bz_+\cup\{+\infty\})^{s}$ and~$a\in\bz$, $-|\boldsymbol{n}|\le a\le |\boldsymbol{m}|$, we define a quiver~$\Gamma_a(\boldsymbol{m},\boldsymbol{n})$ as follows. We set
\begin{multline*}
\Gamma_a(\boldsymbol{m},\boldsymbol{n})_0=\{ (\boldsymbol{x},\boldsymbol{y})=((x_1,\dots,x_r),(y_1,\dots,y_s))\in\bz_+^r\times\bz_+^s\,:\, x_i\le m_i,\, 1\le i\le r,\\ y_j\le n_j,
\,1\le j\le s,\,|\boldsymbol{x}|=|\boldsymbol{y}|+a\}.
\end{multline*}
In other words, $\Gamma(\boldsymbol{m},\boldsymbol{n})$ is just the set of lattice points in the $(r+s)$-dimensional rectangular parallelepiped $[0,m_1]\times\cdots\times [0,m_r]
\times [0,n_1]\times\cdots\times [0,n_s]$
which lie on the hyperplane $z_1+\cdots+z_r-z_{r+1}-\cdots-z_{r+s}=a$.
The arrows are 
$$
(\boldsymbol{x},\boldsymbol{y})\from (\boldsymbol{x}+\boldsymbol{e}_i^{(r)},\boldsymbol{y}+\boldsymbol{e}_j^{(s)}),\qquad x_i<m_i,\,y_j<m_j,\,
1\le i\le r,\, 1\le j\le s,
$$
Note that the map $(\boldsymbol{x},\boldsymbol{y})\mapsto (\boldsymbol{m}-\boldsymbol{x},\boldsymbol{n}-\boldsymbol{y})$ yields an isomorphism 
of quivers~$\Gamma_a(\boldsymbol{m},\boldsymbol{n})\cong \Gamma_{|\boldsymbol{m}|-|\boldsymbol{n}|-a}(\boldsymbol{m},\boldsymbol{n})^{op}$.

For example, $\Gamma_0((n,n),(n))\cong \Gamma_n((n,n),(n))^{op}$ is isomorphic to the following quiver
\begin{equation}\label{mesh-quiv}
\def\dgeverynode{\scriptscriptstyle}\divide\dgARROWLENGTH by 2\dgHORIZPAD=3pt\dgVERTPAD=3pt
\begin{diagram}
\node{(n,0)}\arrow{se}\node{}\node{(n-1,1)}\arrow{se}\arrow{sw}\node{\cdots}\node{(1,n-1)}\arrow{se}\arrow{sw}\node{}\node{(0,n)}\arrow{sw}\\
\node{}\node{\cdots}\arrow{se}\node{\cdots}\node{\cdots}\arrow{sw}\arrow{se}\node{\cdots}\node{\cdots}\arrow{sw}\\
\node{}\node{}\node{(1,0)}\arrow{se}\node{}\node{(0,1)}\arrow{sw}\node{}\\
\node{}\node{}\node{}\node{(0,0)}\node{} \node{}
\end{diagram}
\end{equation}
This is a translation quiver, with~$\tau((x,y))=(x+1,y+1)$, $0\le x+y\le n-2$.
It is easy to see that there is a full embedding of translation quivers of the above quiver into~$\bz\Gamma_{2n+1}$ where~$\Gamma_{2n+1}$ is
any quiver of type~$\mathbb A_{2n+1}$. On the vertices, that embedding is given by $(x,y)\mapsto (-y,n+y-x)$, where we assume that the
vertices of~$\Gamma_{2n+1}$ are numbered from~$0$ to~$2n$. There is also a full embedding of the above quiver into the
Auslander-Reiten quiver of the hereditary algebra of type~$\mathbb A_{2n+1}$ where the~$n$th node is the unique source.

Clearly $\Gamma_{2}((1^3),(1))$ is the quiver of type~$\mathbb D_4$ in which the triple node is the unique source. Two more small examples 
(respectively, $\Gamma_3((1^4),(1^2))$ and~$\Gamma_3((1^4),(1^3))$) are shown below
$$
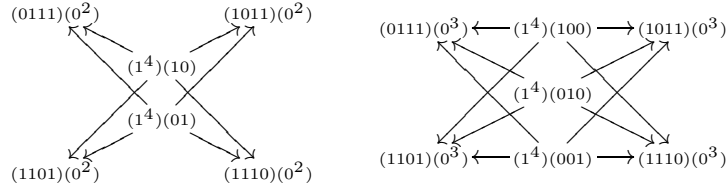
\def\dgeverynode{\scriptscriptstyle}\divide\dgARROWLENGTH by 2\dgHORIZPAD=3pt\dgVERTPAD=3pt
\begin{diagram}
\node{(0111)(0^2)}\node{}\node{(1011)(0^2)}\\
\node{}\node{(1^4)(10)}\arrow{nw}\arrow{ne}\arrow{ssw}\arrow{sse}\\
\node{}\node{(1^4)(01)}\arrow{sw}\arrow{se}\arrow{nnw}\arrow{nne}\\
\node{(1101)(0^2)}\node{}\node{(1110)(0^2)}
\end{diagram}\qquad
\begin{diagram}
\node{(0111)(0^3)}\node{(1^4)(100)}\arrow{w}\arrow{e}\arrow{ssw}\arrow{sse}\node{(1011)(0^3)}\\
\node{}\node{(1^4)(010)}\arrow{ne}\arrow{nw}\arrow{se}\arrow{sw}\\
\node{(1101)(0^3)}\node{(1^4)(001)}\arrow{e}\arrow{w}\arrow{nne}\arrow{nnw}\node{(1110)(0^3)}
\end{diagram}
$$
\begin{lem}
The quiver~$\Gamma_a(\boldsymbol{m},\boldsymbol{n})$ is connected.
\end{lem}
\begin{pf}
Clearly, every vertex in~$\Gamma_a(\boldsymbol{m},\boldsymbol{n})$ is connected to a sink and a vertex $(\boldsymbol{x},\boldsymbol{y})\in\Gamma_a(\boldsymbol{m},\boldsymbol{n})$ is a sink if and only if either~$\boldsymbol{x}=\boldsymbol{0}\in\bz_+^r$ 
or~$\boldsymbol{y}=\boldsymbol{0}\in\bz_+^s$. In particular, if~$a=0$ then~$\Gamma_0(\boldsymbol{m},\boldsymbol{n})$ has a unique sink and
hence is connected. If~$a>0$ (respectively, $a<0$) the sinks in~$\Gamma_a(\boldsymbol{m},\boldsymbol{n})$ are the vertices
$(\boldsymbol{x},\boldsymbol{0})$ (respectively, $(\boldsymbol{0},\boldsymbol{y})$) with~$|\boldsymbol{x}|=a$ (respectively, $|\boldsymbol{y}|=-a$). Suppose
that~$a>0$, the other case being similar. If~$a=|\boldsymbol{m}|$ then we have a unique sink which is also a source. Otherwise, let~$S=\{ \boldsymbol{x}\in
\bz_+^r\,:\, x_i\le m_i,|\boldsymbol{x}|=a\}$ and let~$\prec$ be the lexicographic order on~$S$. Let
$(\boldsymbol{x},\boldsymbol{0})$, $\boldsymbol{x}\in S$ be a sink and suppose that $\boldsymbol{x}$ is not the minimal element of~$S$. Let~$1\le j\le r$
be maximal such that~$x_j<m_j$. If there is~$1\le i<j$ minimal such that~$x_i>0$, we have $(\boldsymbol{x},\boldsymbol{0})\from (\boldsymbol{x}+\boldsymbol{e}_j^{(r)},\boldsymbol{e}_1^{(s)})\to (\boldsymbol{x}-\boldsymbol{e}_i^{(r)}+
\boldsymbol{e}_j^{(r)},\boldsymbol{0})$ and $\boldsymbol{x}-\boldsymbol{e}_i^{(r)}+
\boldsymbol{e}_j^{(r)}\prec\boldsymbol{x}$. Suppose that~$x_i=0$ for all~$i<j$. Since~$\boldsymbol{x}$ is not minimal,
there exists~$\boldsymbol{x}'\in S$ such that $\boldsymbol{x}'\prec \boldsymbol{x}$, that is  $x_i'=0$, $1\le i< j$ and $x_j'<x_j$. Since~$|\boldsymbol{x}'|=a=|\boldsymbol{x}|$, we must have~$x_k'>x_k$
for some~$j<k\le r$, which is a contradiction by the choice of~$j$.
Thus, the connected component of $(\boldsymbol{x},\boldsymbol{0})$ contains 
a sink~$(\boldsymbol{x}',\boldsymbol{0})$ with~$\boldsymbol{x}'\prec\boldsymbol{x}$. The assertion is now immediate.
\end{pf}

\subsection{}\label{A33}
Fix~$1\le i_1<\cdots<i_r<j_1<\cdots<j_s\le \ell\in I$
and consider $\Psi=\{ \alpha_{i_p,j_q}\,:\,1\le p\le r,\, 1\le q\le s\}$. It is easy to see that~$\Psi$ is extremal. Assume further that~$i_{p+1}\not=i_p+1$, $j_{q+1}\not=j_q+1$ for all~$1\le p<r$, $1\le q<s$, 
and so by~\lemref{A30} ~$\Psi$ is
regular. 
\begin{prop}
Let~$\lambda\in P^+$.
Then the quiver~$\Delta_\Psi[\lambda]$
is isomorphic to~$\Gamma_a(\boldsymbol{m},\boldsymbol{n})$ where~$\boldsymbol{m}=(\lambda(h_{i_p-1})+\lambda(h_{i_p}))_{1\le p\le r}$, $\boldsymbol{n}=(
\lambda(h_{j_q})+\lambda(h_{j_q+1}))_{1\le q\le s}$ and $a=\sum_{p=1}^r\lambda(h_{i_p})-\sum_{q=1}^s\lambda(h_{j_q})$.
\end{prop}
\begin{pf}
Let~$J=\{i_p,i_p-1\,:\,1\le p\le r\}\cup\{j_q, j_q+1\,:\,1\le q\le s\}$.
Suppose that~$\mu\in\Delta_\Psi[\lambda]_0$. Since~$\Delta_\Psi[\lambda]_0\subset (\lambda+\bz\Psi)\cap P^+$,
we have~$\mu(h_j)=\lambda(h_j)$, $j\notin J$, and
\begin{alignat*}{3}
&\mu(h_{i_p})=\lambda(h_{i_p})+\sum_{q=1}^s x_{p,q},&\qquad&\mu(h_{i_p-1})=\lambda(h_{i_p-1})-\sum_{q=1}^s x_{p,q},&\qquad &1\le p\le r
\\
&\mu(h_{j_q})=\lambda(h_{j_q})+\sum_{p=1}^r x_{p,q},&&\mu(h_{j_q+1})=\lambda(h_{j_q+1})-\sum_{p=1}^r x_{p,q},&& 1\le q\le s
\end{alignat*}
where~$x_{p,q}\in\bz$, $1\le p\le r$, $1\le q\le s$. It follows that~$\Delta_\Psi[\lambda]_0$ is contained in the set~$S(\lambda)$ of~$\mu\in P^+$
satisfying the following conditions
\begin{align*}
&\mu(h_{i_p-1})+\mu(h_{i_p})=\lambda(h_{i_p-1})+\lambda(h_{i_p}),\quad \mu(h_{j_q})+\mu(h_{j_q+1})=\lambda(h_{j_q})+\lambda(h_{j_q+1})\\
&\sum_{p=1}^r \mu(h_{i_p})-\sum_{p=1}^s\mu(h_{j_q})=\sum_{p=1}^r \lambda(h_{i_p})-\sum_{p=1}^s\lambda(h_{j_q}).
\end{align*}
Clearly, if~$\mu\in S(\lambda)$ then~$\mu^-\subset S(\lambda)$ and so~$S(\lambda)$ defines a convex subquiver of~$\Delta_\Psi$ containing
$\Delta_\Psi[\lambda]$ as a full connected subquiver.
Define a map $S(\lambda)\to\Gamma_a(\boldsymbol{m},\boldsymbol{n})_0$ by
$$
\mu\mapsto ((\mu(h_{i_1}),\dots,\mu(h_{i_r})),(\mu(h_{j_1}),\dots,\mu(h_{j_s}))).
$$
This map is clearly a bijection and it is easy to see that it induces an isomorphism of quivers. Since by~\lemref{A32} the quiver $\Gamma_a(\boldsymbol{m},\boldsymbol{n})$ is connected, the assertion follows.
\end{pf}

\subsection{}\label{A90}
For~$1\le i<j\le\ell$, we fix root vectors~$e_{i,j}\in \lie g_{\alpha_{i,j}}$ such that
\begin{equation}\label{A90.0}
[e_r,e_{p,q}]=\delta_{r,p-1} e_{r,q}-\delta_{r,q+1} e_{p,r},\qquad  [f_r,e_{p,q}]=\delta_{r,p} e_{r+1,q}-\delta_{r,q} e_{p,r-1}.
\end{equation}
For example, the standard basis of the matrix realisation of~$\lie{sl}_{\ell+1}$ has these properties.

Fix~$\alpha_{i,j}\in \Psi$. 
Clearly $\{\gamma\in R^+\,:\, \alpha_{i,j}\le \gamma\}=\{\alpha_{p,q}\,:\, 1\le p\le i,j\le q\le\ell\}$. If~$\lambda,\lambda+\alpha_{i,j}\in P^+$, we have
$\lambda(h_{i-1}),\lambda(h_{j+1})>0$ and so
$$
\lambda(\cal H_{t,i-1})\ge \lambda(h_{i-1})>0,\quad \lambda(\cal H_{j+1,t})\ge \lambda(h_{j+1})>0,\qquad
1\le t\le i-1,\, j+1\le t\le \ell.
$$ 
Therefore
$$
\cal H_{r,i-1},\cal H_{j+1,s}\in F_{\alpha_{i,j}}(\lie h)^\times,\qquad 1\le r\le i-1,\,j+1\le s\le \ell.
$$
For all~$1\le p\le i$, $j\le q\le \ell$, define $\bou_{\alpha_{i,j},\alpha_{p,q}}\in  U(\lie b)_{\alpha_{i,j}-\alpha_{p,q}}\tensor_{S(\lie h)} F_{\alpha_{i,j}}(\lie h)$ by
\begin{equation}\label{A90.10}
\bou_{\alpha_{i,j},\alpha_{p,q}}=(-1)^{i-p} \cal X^-_{p,i-1} \cal X^+_{j+1,q} \tensor \prod_{t=p}^{i-1} \cal H_{t,i-1}^{-1} \prod_{t=j+1}^{q} \cal H_{j+1,t}^{-1}.
\end{equation}
\begin{lem}
Let $\alpha_{i,j}\in\Psi$, $1\le i<j\le \ell$. Then $(\bou_{\alpha_{i,j},\alpha_{p,q}}\,:\, 1\le p\le i,\,j\le q\le \ell)$ is
an adapted family for~$\alpha_{i,j}$.
\end{lem}
\begin{pf}
We have $\pi_{\lambda,\alpha_{i,j}}(\bou_{\alpha_{i,j},\alpha_{p,q}})=\cal X^-_{p,i-1}(\lambda) \cal X^+_{j+1,q} (\lambda) B_{p,q}(i,j,\lambda)$,
where
\begin{equation}\label{A90.10b}
B_{p,q}(i,j,\lambda):=(-1)^{i-p} \prod_{t=p}^{i-1} (\lambda(\cal H_{t,i-1}))^{-1} \prod_{t=j+1}^{q} \lambda(\cal H_{j+1,t})^{-1}.
\end{equation}
Write~$B_{p,q}=B_{p,q}(i,j,\lambda)$ to shorten the notation.
It is immediate that the~$B_{p,q}$ satisfy 
\begin{equation}\label{A90.20}
  B_{p+1,q}=-\lambda(\cal H_{p,i-1}) B_{p,q},\qquad 
  B_{p,q-1}=\lambda(\cal H_{j+1,q}) B_{p,q}
\end{equation}
Let
$$
v=\sum_{\alpha_{i,j}\le\gamma} e_\gamma\tensor \pi_{\lambda,\alpha_{i,j}}(\bou_{\alpha_{i,j},\gamma}) v_\lambda.
$$
Since~$v\in (\lie g\tensor V(\lambda))_{\lambda+\beta}$ and~$\bou_{\alpha_{i,j},\alpha_{i,j}}=1$,
it remains to prove that~$e_r v=0$
for all~$1\le r\le \ell$. 
It follows immediately from~\lemref{A50} and~\eqref{A90.0} that
$e_r v=0$, $i\le r\le j$. 
Note that~$\cal X^-_{p,i-1}(\lambda)$ and~$\cal X^+_{j+1,q}(\lambda)$ commute for all~$1\le p\le i$, $j+1\le q\le\ell$.
Let $1\le r\le i-1$. By~\lemref{A50}
$$
e_{r} \cal X^-_{p,i-1}(\lambda)\cal X^+_{j+1,q}(\lambda)v_\lambda =
\delta_{r,p} \lambda(\cal H_{r,i-1})\cal X^-_{r+1,i-1}(\lambda)\cal X^+_{j+1,q}(\lambda)v_\lambda,
$$
hence, using~\eqref{A90.0} and~\eqref{A90.20}, we obtain
$$
e_{r} v=\sum_{q=j+1}^{\ell} (B_{r+1,q}+B_{r,q}\lambda(\cal H_{r,i-1}))
e_{r,q}\tensor \cal X^-_{r+1,i-1}(\lambda)\cal X^+_{j+1,q}(\lambda)v_\lambda=0.
$$
Finally, suppose that $j+1\le r\le \ell$. 
By~\lemref{A50}
$$
e_{r} \cal X^-_{p,i-1}(\lambda)\cal X^+_{j+1,q}(\lambda)v_\lambda =
\delta_{r,q} \lambda(\cal H_{j+1,r})\cal X^-_{p,i-1}(\lambda)\cal X^+_{j+1,r-1}(\lambda)v_\lambda
$$
and so by~\eqref{A90.0} and~\eqref{A90.20} 
\begin{equation*}
e_{r} v=
\sum_{p=1}^{i} (-B_{p,r-1}+B_{p,r}\lambda(\cal H_{j+1,r}))e_{p,r}\tensor \cal X^-_{p,i-1}(\lambda)\cal X^+_{j+1,r-1}(\lambda)v_\lambda=0.\qedhere
\end{equation*}
\end{pf}

\subsection{}\label{A95}
To describe the relations without ambiguity, we need to fix an isomorphism $\bt_\Psi^{\lie g}\to\bc\Delta_\Psi$ which amounts to
fixing an element $\boldsymbol{z}\in (\bc^\times)^{(\Delta_\Psi)_1}$. 
Given~$\alpha_{i,j}\in\Psi$, let
\begin{equation}\label{A95.0}
\cal Z_{\alpha_{i,j},\Psi}=\prod_{1\le t<i:\alpha_{t,j}\in\Psi} 
\cal H_{t,i-1}
\prod_{j<t\le \ell\,:\, \alpha_{i,t}\in\Psi} 
\cal H_{j+1,t}.
\end{equation}
Since~$\lambda+\alpha_{i,j}\in P^+$ implies that~$\lambda(h_{i-1}),\lambda(h_{j+1})>0$, $\cal Z_{\alpha_{i,j},\Psi}\in F_{\alpha_{i,j}}(\lie h)^\times$.
Let~$\boldsymbol{z}$ be the image of $(\cal Z_{\beta,\Psi}^{-1})_{\beta\in\Psi}\in \prod_{\beta\in\Psi} F_\beta(\lie h)^\times$ in~$G_\Psi$.
In other words,
we fix the isomorphism~$\Phi$ by requiring
$$
\boldsymbol{a}_{\lambda,\alpha_{i,j}}\mapsto  \lambda(\cal Z_{\alpha_{i,j},\Psi})^{-1}(\lambda\from\lambda+\alpha_{i,j}),\qquad \lambda,\lambda+\alpha_{i,j}\in P^+.
$$

Now we are ready to compute all relations in the algebra~$\bs_\Psi^{\lie g}$. For readers convenience, we describe different cases
in separate propositions.  

\begin{prop} Let~$\Psi$ be an extremal set, $|\Psi|>1$.
\begin{enumerate}[{\rm(i)}]
\item\label{A95.i}
Let $\alpha_{i,j},\alpha_{i,k}\in\Psi$, $1\le i\le j<k\le \ell$. Then for all~$\lambda\in P^+$ such that 
$t_{\lambda,\alpha_{i,j}+\alpha_{i,k}}=2$,
$\mathfrak R_\Psi(\lambda,\lambda+\alpha_{i,j}+\alpha_{i,k})$ is spanned by the commutativity relation.
If~$t_{\lambda,\alpha_{i,j}+\alpha_{i,k}}=1$, $\dim\mathfrak R_\Psi(\lambda,\lambda+\alpha_{i,j}+\alpha_{i,k})=1$.

\item\label{A95.iii} Let $\alpha_{i,k},\alpha_{j,k}\in \Psi$,  $1\le i<j\le k\le \ell$. Then for all~$\lambda\in P^+$ such that
$t_{\lambda,\alpha_{i,k}+\alpha_{j,k}}=2$,
$\mathfrak R_\Psi(\lambda,\lambda+\alpha_{i,k}+\alpha_{j,k})$ is spanned by the commutativity relation.
If~$t_{\lambda,\alpha_{i,k}+\alpha_{j,k}}=1$, $\dim \mathfrak R_\Psi(\lambda,\lambda+\alpha_{i,k}+\alpha_{j,k})=1$. 
\end{enumerate}
In particular, if~$\eta\in\Psi+\Psi$ satisfies $m_\eta=2$, then $\cal N_\eta=\emptyset$.
\end{prop}
\begin{pf}
We present a detailed argument here since the computations of this kind will be used repeatedly in the rest of this paper and in the future we 
will omit most of the details.

Retain the notations of the proof of~\lemref{A90}.
To prove~\eqref{A95.i}, note that $\alpha_{i,j}<\alpha_{i,k}$. It follows from~\propref{SP55} and~\lemref{A90} that
\begin{align*}
&\Pi_\lambda(\alpha_{i,j},\alpha_{i,k})=e_{i,k}\tensor e_{i,j}
\\
&\Pi_\lambda(\alpha_{i,k},\alpha_{i,j})=e_{i,j}\tensor e_{i,k}+e_{i,k}\tensor \bou_{\alpha_{i,j},\alpha_{i,k}}(\mu) e_{i,k},\qquad \mu=\lambda+\alpha_{i,k},
\end{align*}
and
$$
\bou_{\alpha_{i,j},\alpha_{i,k}}(\mu) =\pi_{\mu,\alpha_{i,j}}(\bou_{\alpha_{i,j},\alpha_{i,k}})= B_{i,k}(i,j,\mu)\cal X^+_{j+1,k}(\mu).
$$
Note that~$\mu(h_{k})=\lambda(h_{k})+1>0$. Let~$\sigma\in\Sigma(j+1,k)$ and suppose that~$\bof_\sigma\notin\Ann_{ U(\lie g)} e_{i,k}$. Then~\eqref{A90.0} we
must have $\sigma(k)=k-j$ and so~$\bof_\sigma e_{i,k}=-\bof_{\sigma'} e_{i,k-1}$, where~$\sigma'\in\Sigma(j+1,k-1)$ is the restriction of~$\sigma$.
Following this way we conclude that $\sigma(r)=r-j$, $j+1\le r\le k$, that is~$\bof_\sigma=f_{j+1}\cdots f_k$. Since~$\mu(h_{k})>0$, using~\propref{A45} and \lemref{A50} we
conclude that
\begin{equation}\label{A95.5}
\cal X^+_{j+1,k}(\mu)=(-1)^{k-j-1} \prod_{t=j+1}^{k-1} \mu(\cal H_{j+1,t}) f_{j+1}\cdots f_k+\Ann_{ U(\lie g)} v_\mu\cap \Ann_{ U(\lie g)} e_{i,k}
\end{equation}
and $f_{j+1}\cdots f_k\notin\Ann_{U(\lie n^-)} v_\mu$.
Therefore, 
$$
\bou_{\alpha_{i,j},\alpha_{i,k}}(\mu) e_{i,k}=-(\mu(\cal H_{j+1,k}))^{-1} e_{i,j}
$$
hence
\begin{equation}\label{A95.10}
\Pi_\lambda(\alpha_{i,k},\alpha_{i,j})=e_{i,j}\tensor e_{i,k}-(\lambda(\cal H_{j+1,k})+1)^{-1} e_{i,k}\tensor e_{i,j}.
\end{equation}
It is easy to see that the intersection of
$\bc \Pi_{\lambda}(\alpha_{i,j},\alpha_{i,k})+\bc \Pi_{\lambda}(\alpha_{i,k},\alpha_{i,j})$ with~$\bigwedge^2 \lie n^+_\Psi$ 
is spanned by
$$
\lambda(\cal H_{j+1,k})\Pi_{\lambda}(\alpha_{i,j},\alpha_{i,k})-(\lambda(\cal H_{j+1,k})+1) \Pi_{\lambda}(\alpha_{i,k},\alpha_{i,j}).
$$
Thus, $\mathfrak R_\Psi(\lambda,\lambda+\eta)$, $\eta=\alpha_{i,j}+\alpha_{i,k}$ is spanned by
\begin{multline*}
\lambda(\cal H_{j+1,k}) (\lambda(\cal Z_{\alpha_{i,j},\Psi})(\lambda+\alpha_{i,j})(\cal Z_{\alpha_{i,k},\Psi}))^{-1}
(\lambda\from\lambda+\alpha_{i,j}\from\lambda+\eta)\\-(\lambda(\cal H_{j+1,k})+1)(\lambda(\cal Z_{\alpha_{i,k},\Psi})(\lambda+\alpha_{i,k})(\cal Z_{\alpha_{i,j},
\Psi}))^{-1}(\lambda\from\lambda+\alpha_{i,k}\from\lambda+\eta)
\end{multline*}
Since $\alpha_{t,j}+\alpha_{i,k}=\alpha_{t,k}+\alpha_{i,j}$, $1\le t<i$, \corref{*} implies 
that for all~$1\le t<i$, $\alpha_{t,j}\in\Psi$ if and only if~$\alpha_{t,k}\in\Psi$. Then
\begin{multline*}
\lambda(\cal Z_{\alpha_{i,k},\Psi}\cal Z_{\alpha_{i,j},\Psi}^{-1})(\lambda+\alpha_{i,k})(\cal Z_{\alpha_{i,j},\Psi})
(\lambda+\alpha_{i,j})(\cal Z_{\alpha_{i,k},\Psi}^{-1})
\\=
\frac{(\lambda(\cal H_{j+1,k})+1)\displaystyle\prod_{1\le t<i\,:\, \alpha_{t,k}\in\Psi}\lambda(\cal H_{t,i-1})
\prod_{1\le t<i:\alpha_{t,j}\in\Psi}
(\lambda(\cal H_{t,i-1})-1)}{\lambda(\cal H_{j+1,k}) \displaystyle\prod_{1\le t<i\,:\, \alpha_{t,j}\in\Psi}\lambda(\cal H_{t,i-1})
\prod_{1\le t<i:\alpha_{t,k}\in\Psi}
(\lambda(\cal H_{t,i-1})-1)}=\frac{\lambda(\cal H_{j+1,k})+1}{\lambda(\cal H_{j+1,k})},
\end{multline*}
and so $\mathfrak R_\Psi(\lambda,\lambda+\eta)$ is spanned by the commutativity relation.
If~$t_{\lambda,\eta}=1$, by~\eqref{A30.2p.l} we have~$\lambda(\cal H_{j+1,k})=\lambda(h_{j+1})=0$ and it follows 
from~\eqref{A95.10} that~$\Pi_{\lambda}(\alpha_{i,k},\alpha_{i,j})
\in\bigwedge^2 \lie n^+_\Psi$. Thus, the unique path~$(\lambda\from\lambda+\alpha_{i,j+1}\from\lambda+\eta)$ in~${\Delta_\Psi}(\lambda,\lambda+\eta)$ is a relation.

To prove~\eqref{A95.iii}, note that~$\alpha_{j,k}<\alpha_{i,k}$ and so
\begin{align*}
&\Pi_{\lambda}(\alpha_{j,k},\alpha_{i,k})=e_{i,k}\tensor e_{j,k},
\\
&\Pi_{\lambda}(\alpha_{i,k},\alpha_{j,k})=e_{j,k}\tensor e_{i,k}+e_{i,k}\tensor \bou_{\alpha_{j,k},\alpha_{i,k}}(\nu)e_{i,k}
\qquad \nu=\lambda+\alpha_{i,k},
\end{align*}
where
$$
\bou_{\alpha_{j,k},\alpha_{i,k}}(\nu)=\pi_{\nu,\alpha_{j,k}}(\bou_{\alpha_{j,k},\alpha_{i,k}})=B_{i,k}(j,k,\nu)\cal X^-_{i,j-1}(\nu).
$$
An argument similar to the above shows that
\begin{equation}\label{A95.15}
\cal X^-_{i,j-1}(\nu)=(-1)^{j-i-1} \prod_{t=i+1}^{j-1} \nu(\cal H_{t,j-1}) f_{j-1}\cdots f_i +\Ann_{ U(\lie g)} v_\nu\cap \Ann_{ U(\lie g)} e_{i,k},
\end{equation}
hence
\begin{equation}\label{A95.20}
 \Pi_{\lambda}(\alpha_{i,k},\alpha_{j,k})=e_{j,k}\tensor e_{i,k}-(\lambda(\cal H_{i,j-1})+1)^{-1} e_{i,k}\tensor e_{j,k}.
\end{equation}
To finish the computation, we observe that $\alpha_{i,k}+\alpha_{j,t}=\alpha_{j,k}+\alpha_{i,t}$, $k<t\le\ell$, hence by~\corref{*}
$\alpha_{i,t}\in\Psi$ if and only if~$\alpha_{j,t}\in\Psi$ for all~$k<t\le \ell$. 
This implies that
$$
\lambda(\cal Z_{\alpha_{i,k},\Psi}\cal Z_{\alpha_{j,k},\Psi}^{-1})(\lambda+\alpha_{i,k})(\cal Z_{\alpha_{j,k},\Psi})(\lambda+\alpha_{j,k})
(\cal Z_{\alpha_{i,k},\Psi}^{-1})
=(\lambda(\cal H_{i,j-1})+1)(\lambda(\cal H_{i,j-1}))^{-1}.
$$
Finally, if~$t_{\lambda,\alpha_{i,k}+\alpha_{j,k}}=1$, \eqref{A30.2p.r} implies 
that~$\lambda(\cal H_{i,j-1})=0$ hence~$\Pi_{\lambda}(\alpha_{i,k},\alpha_{j,k})\in\bigwedge^2\lie n^+_\Psi$ and so
the corresponding path is a relation.
\end{pf}

\begin{ex}
Fix~$i<j<k\in I$ with~$k\not=i+1$ and either~$i\not=1$ or~$k\not=\ell$. Let~$\lambda=m(\varpi_{i-1}+\varpi_{j+1}+\varpi_{k+1})$. Then by~\propref{TQF10}
and by the above,
$\bs_\Psi^{\lie g}(\lambda\le_\Psi)$ has global dimension~$2$ and is isomorphic to the path algebra of the translation quiver~\eqref{mesh-quiv} with the mesh relations. In particular,
it is isomorphic to a subalgebra of the Auslander algebra of the path algebra of the quiver of type~$\mathbb A_{2m+1}$, where the node preserved by the diagram automorphism is
the unique sink.
\end{ex}

\subsection{}\label{A110}
The following proposition completes the proof of~\thmref{mainthm} for~$\lie g$ of type~$A$.
\begin{prop}
Let~$\Psi\subset R^+$, $|\Psi|\ge4$ be extremal. Suppose that 
$$
\{\alpha_{i,k},\alpha_{j,k},\alpha_{i,m},\alpha_{j,m}\}\subset\Psi, \qquad i<j\le k<m
$$ 
and let~$\eta=\alpha_{i,k}+\alpha_{j,m}=\alpha_{i,m}+\alpha_{j,k}$. Let~$x_\lambda=\lambda(\cal H_{i,j-1})$, $y_\lambda=\lambda(\cal H_{k+1,m})$.
\begin{enumerate}[{\rm(i)}]
 \item\label{A110.i} Suppose that
$t_{\lambda,\eta}=4$, $x_\lambda\not=y_\lambda$. Then $\mathfrak R_\Psi(\lambda,\lambda+\eta)$ is spanned by
\begin{align*}
&(x_\lambda+1)(y_\lambda+2)(\lambda\from\lambda+\alpha_{i,k}\from\lambda+\eta)-(x_\lambda+2)(y_\lambda+1)(\lambda\from\lambda+\alpha_{j,m}\from\lambda+\eta)\\
&\qquad-
(x_\lambda-y_\lambda) (\lambda\from\lambda+\alpha_{i,m}\from\lambda+\eta)
\\
\intertext{and}
&x_\lambda(y_\lambda+1)(\lambda\from\lambda+\alpha_{i,k}\from\lambda+\eta)-(x_\lambda+1)y_\lambda(\lambda\from\lambda+\alpha_{j,m}\from\lambda+\eta)\\
&\qquad-(x_\lambda-y_\lambda)(\lambda\from\lambda+\alpha_{j,k}\from\lambda+\eta)
\end{align*}
\item\label{A110.i''}
Suppose that~$t_{\lambda,\eta}=4$ and~$x_\lambda=y_\lambda$. Then $\mathfrak R_\Psi(\lambda,\lambda+\eta)$ is spanned by
\begin{align*}
&(\lambda\from\lambda+\alpha_{i,k}\from\lambda+\eta)-(\lambda\from\lambda+\alpha_{j,m}\from\lambda+\eta)\\
\intertext{and}
&2(\lambda\from\lambda+\alpha_{i,k}\from\lambda+\eta)-x_\lambda(\lambda\from\lambda+\alpha_{i,m}\from\lambda+\eta)\\
&\qquad-(x_\lambda+2)(\lambda\from\lambda+\alpha_{j,k}\from\lambda+\eta).
\end{align*}

\item\label{A110.ii}  Suppose that $t_{\lambda,\eta}=2$. Then $x_\lambda\not=y_\lambda$ and either~$i=j-1$, $x_\lambda=0$
and the relation is
$$
y_\lambda(\lambda\from\lambda+\alpha_{i,m}\from\lambda+\eta)+(y_\lambda+2)(\lambda\from\lambda+\alpha_{i,k}\from\lambda+\eta).
$$
or~$m=k+1$, $y_\lambda=0$ and 
the relation is
$$
x_\lambda(\lambda\from\lambda+\alpha_{i,m}\from\lambda+\eta)+(x_\lambda+2)(\lambda\from\lambda+\alpha_{j,m}\from\lambda+\eta).
$$

\item\label{A110.iv} Suppose that~$t_{\lambda,\eta}=1$. Then~$i=j-1$, $m=k+1$, $x_\lambda=y_\lambda=0$ and
and~$\mathfrak R_\Psi(\lambda,\lambda+\eta)=0$. 
\end{enumerate}
Thus, $\cal N_{\eta}$ is contained in the set $P^+\cap \{\xi\in\lie h^*\,:\, \xi(\cal H_{i,j-1}-\cal H_{k+1,m})=0\}$ and coincides with this set if~$\Psi$
is regular.
\end{prop}

\begin{pf}
We have $\alpha_{j,k}<
\alpha_{i,k},\alpha_{j,m}<\alpha_{i,m}$ while $\alpha_{i,k}$,
$\alpha_{j,m}$ are not comparable in the standard partial order.

To prove~\eqref{A110.i} we compute using Lemma~\ref{A90}, \propref{SP55} and~\eqref{A95.5}, \eqref{A95.15}
\begin{subequations}
\begin{align}
&\Pi_\lambda(\alpha_{j,k},\alpha_{i,m})=e_{i,m}\tensor e_{j,k},\label{A110.p1}\\
&\Pi_\lambda(\alpha_{i,m},\alpha_{j,k})=
e_{j,k}\tensor e_{i,m}-(x_\lambda+1)^{-1} e_{i,k}\tensor e_{j,m}-(y_\lambda+1)^{-1} e_{j,m}\tensor e_{i,k}\nonumber\\&
\phantom{\Pi_\lambda(\alpha_{i,m},\alpha_{j,k})=}+(x_\lambda+1) (y_\lambda+1)^{-1} e_{i,m}\tensor e_{j,k}\label{A110.p2}\\
&\Pi_\lambda(\alpha_{j,m},\alpha_{i,k})= e_{i,k}\tensor e_{j,m}-(y_\lambda+1)^{-1} e_{i,m}\tensor e_{j,k}
\label{A110.p3}\\
&\Pi_\lambda(\alpha_{i,k},\alpha_{j,m})= e_{j,m}\tensor e_{i,k}-(x_\lambda+1)^{-1} e_{i,m}\tensor e_{j,k}.\label{A110.p4}
\end{align}
\end{subequations}
In particular, we see that none of the paths in~${\Delta_\Psi}(\lambda,\lambda+\eta)$ is a relation.
Furthermore, we have
\begin{alignat*}{2}
&\alpha_{t,k}\in\Psi\iff \alpha_{t,m}\in\Psi,&\qquad &1\le t< j,\\
&\alpha_{i,t}\in\Psi\iff \alpha_{j,t}\in\Psi,&\qquad &k<t\le \ell.
\end{alignat*}
Indeed this follows from~\corref{*} by observing that $\alpha_{j,m}+\alpha_{t,k}=\alpha_{t,m}+\alpha_{j,k}$, $\alpha_{i,t}+\alpha_{j,k}=\alpha_{i,k}+\alpha_{j,t}$.
Then if we set~$z=\lambda(\cal Z_{\alpha_{j,k},\Psi})(\lambda+\alpha_{j,k})(\cal Z_{\alpha_{i,m},\Psi})$,
\begin{align*}
&\lambda(\cal Z_{\alpha_{i,m},\Psi})(\lambda+\alpha_{i,m})(\cal Z_{\alpha_{j,k},\Psi})=(x_\lambda+1)(y_\lambda+1)%
x_\lambda^{-1} y_\lambda^{-1} z
\\
&\lambda(\cal Z_{\alpha_{j,m},\Psi})(\lambda+\alpha_{j,m})(\cal Z_{\alpha_{i,k},\Psi})=(y_\lambda+1)y_\lambda^{-1} z\\
&\lambda(\cal Z_{\alpha_{i,k},\Psi})(\lambda+\alpha_{i,k})(\cal Z_{\alpha_{j,m},\Psi})=(x_\lambda+1) x_\lambda^{-1} z.
\end{align*}
The relations in~\eqref{A110.i} and in~\eqref{A110.i''} are now straightforward. 

To prove~\eqref{A110.ii} observe that in these cases we have, respectively,
\begin{align*}
&\lambda(\cal Z_{\alpha_{i,m},\Psi}\cal Z_{\alpha_{i+1,m},\Psi}^{-1})(\lambda+\alpha_{i,m})(\cal Z_{\alpha_{i+1,k},\Psi})(\lambda+\alpha_{i+1,m})(\cal Z_{\alpha_{i,k},\Psi}^{-1})=(y_\lambda+1)y_\lambda^{-1},\\
&\lambda(\cal Z_{\alpha_{i,k+1},\Psi}\cal Z_{\alpha_{j,k+1},\Psi}^{-1})(\lambda+\alpha_{i,k+1})(\cal Z_{\alpha_{j,k},\Psi})(\lambda+\alpha_{j,k+1})
(\cal Z_{\alpha_{i,k},\Psi}^{-1})
=(x_\lambda+1)x_\lambda^{-1}.
\end{align*}
The relations now follow easily from the above and~\eqref{A110.p2},\eqref{A110.p4} (respectively, \eqref{A110.p2} and~\eqref{A110.p3}).
Finally, in the last case
$\Pi_{\lambda}(\alpha_{i,k+1},\alpha_{i+1,k})\notin\bigwedge^2 \lie n^+_\Psi$, hence the unique path $\lambda\from\lambda+\alpha_{i,k+1}
\from\lambda+\eta$ is not a relation.  
\end{pf}

\subsection{}
Retain the notations and the assumptions of~\ref{A33}. 
Then by~\propref{A33}, a connected subalgebra of $\bt^{\lie g}_\Psi$ corresponding to a connected component of~$\Delta_\Psi$
is isomorphic to the path algebra of the quiver~$\Gamma_a(\boldsymbol{m},\boldsymbol{n})$
for some~$\boldsymbol{m}\in(\bz_+\cup\{+\infty\})^r$, $\boldsymbol{n}\in(\bz_+\cup\{+\infty\})^s$ and~$-|\boldsymbol{n}|\le a\le |\boldsymbol{m}|$.
However, this isomorphism looses some information which is necessary for describing relations in~$\bs_\Psi^{\lie g}$, since the latter depend on
$\mu(\cal H_{i_p,i_{p'}-1})$, $\mu(\cal H_{j_q+1,j_{q'}})$. Given~$\lambda\in P^+$, set
$${z}(\lambda)_p^-=\lambda(\cal H_{i_p+1,i_{p+1}-2})+2,\qquad 
{z}_q^+(\lambda)=\lambda(\cal H_{j_q+2,j_{q+1}-1})+2.
$$
These parameters are obviously constant on connected components of~$\Delta_\Psi$ and can take arbitrary integer values $\ge2$.
Let~$(\boldsymbol{x},\boldsymbol{y})=((x_1,\dots,x_r),(y_1,\dots,y_s))$ be the image of~$\mu\in\Delta_\Psi[\lambda]_0$
under the isomorphism of quivers constructed in~\propref{A33}. Then we have
$$
\mu(\cal H_{i_p,i_{p+1}-1})=x_p+m_{p+1}-x_{p+1}+{z}(\lambda)_p^-,\qquad \mu(\cal H_{j_q+1,j_{q+1}})=n_q-y_q+y_{q+1}+{z}(\lambda)_q^+
$$
and
\begin{align*}
&M_{p,p'}(\boldsymbol{x}):=\mu(\cal H_{i_p,i_{p'}-1})=x_p-x_{p'}+\sum_{k=p}^{p'-1} (m_{k+1}+{z}(\lambda)_{k}^-)+p'-p-1,\qquad 1\le p<p'\le r,
\\
&N_{q,q'}(\boldsymbol{y}):=\mu(\cal H_{j_q+1,j_{q'}})=y_{q'}-y_{q}+\sum_{k=q}^{q'-1} (n_{k}+{z}(\lambda)_{k}^+)+q'-q-1,\qquad 1\le q<q'\le s.
\end{align*}
Thus, the isomorphism of the connected subalgebra of~$\bt^{\lie g}_\Psi$ corresponding to~$\Delta_\Psi[\lambda]$, $\lambda\in P^+$
onto $\bc\Gamma_a(\boldsymbol{m},\boldsymbol{n})$ provided by~\propref{A33} induces the following relations on~$\Gamma_a(\boldsymbol{m},\boldsymbol{n})$. 
First, for all~$1\le p\le r$, $1\le q<q'\le s$
and for all~$(\boldsymbol{x},\boldsymbol{y})$ such that~$(\boldsymbol{x}+2\boldsymbol{e}^{(r)}_p,\boldsymbol{y}+
\boldsymbol{e}^{(s)}_q+\boldsymbol{e}^{(s)}_{q'})\in\Gamma_a(\boldsymbol{m},\boldsymbol{n})_0$,
we have a commutativity relation
\begin{multline*}
((\boldsymbol{x},\boldsymbol{y})\from (\boldsymbol{x}+\boldsymbol{e}^{(r)}_p,\boldsymbol{y}+\boldsymbol{e}^{(s)}_q)\from (\boldsymbol{x}+2\boldsymbol{e}^{(r)}_p,\boldsymbol{y}+
\boldsymbol{e}^{(s)}_q+\boldsymbol{e}^{(s)}_{q'}))\\-((\boldsymbol{x},\boldsymbol{y})\from (\boldsymbol{x}+\boldsymbol{e}^{(r)}_p,\boldsymbol{y}+\boldsymbol{e}^{(s)}_{q'})\from (\boldsymbol{x}+2\boldsymbol{e}^{(r)}_p,\boldsymbol{y}+
\boldsymbol{e}^{(s)}_q+\boldsymbol{e}^{(s)}_{q'})).
\end{multline*}
Similarly, for all~$1\le p<p'\le r$ and for all~$1\le q\le s$ such that $(\boldsymbol{x}+\boldsymbol{e}^{(r)}_p+\boldsymbol{e}^{(r)}_{p'},\boldsymbol{y}+
2\boldsymbol{e}^{(s)}_q))\in\Gamma_a(\boldsymbol{m},\boldsymbol{n})_0$, we have the commutativity relation
\begin{multline*}
((\boldsymbol{x},\boldsymbol{y})\from (\boldsymbol{x}+\boldsymbol{e}^{(r)}_p,\boldsymbol{y}+\boldsymbol{e}^{(s)}_q)\from (\boldsymbol{x}+\boldsymbol{e}^{(r)}_p+\boldsymbol{e}^{(r)}_{p'},\boldsymbol{y}+
2\boldsymbol{e}^{(s)}_q))\\-((\boldsymbol{x},\boldsymbol{y})\from (\boldsymbol{x}+\boldsymbol{e}^{(r)}_{p'},\boldsymbol{y}+\boldsymbol{e}^{(s)}_{q})\from (\boldsymbol{x}+\boldsymbol{e}^{(r)}_p+\boldsymbol{e}^{(r)}_{p'},\boldsymbol{y}+
2\boldsymbol{e}^{(s)}_q)).
\end{multline*}
Finally, for all~$1\le p<p'\le r$, $1\le q<q'\le s$, let~$\boldsymbol{x}'=\boldsymbol{x}+\boldsymbol{e}^{(r)}_p+\boldsymbol{e}^{(r)}_{p'}$, $\boldsymbol{y}'=\boldsymbol{y}+\boldsymbol{e}^{(s)}_q+\boldsymbol{e}^{(s)}_{q'}$. Assume that~$(\boldsymbol{x}',\boldsymbol{y}')\in\Gamma_a(\boldsymbol{m},\boldsymbol{n})_0$.
If~$M_{p,p'}(\boldsymbol{x})\not=N_{q,q'}(\boldsymbol{y})$
we have
\begin{align*}
&(M_{p,p'}(\boldsymbol x)+1)(N_{q,q'}(\boldsymbol y)+2)((\boldsymbol x,\boldsymbol y)\from (\boldsymbol x+\boldsymbol{e}^{(r)}_p,\boldsymbol{y}+\boldsymbol{e}^{(s)}_q)\from (\boldsymbol{x}',\boldsymbol{y}'))\\&\qquad-(M_{p,p'}(\boldsymbol x)+2)(N_{q,q'}(\boldsymbol y)+1)((\boldsymbol x,\boldsymbol y)\from(\boldsymbol x+\boldsymbol{e}^{(r)}_{p'},\boldsymbol y+\boldsymbol{e}^{(s)}_{q'})\from(\boldsymbol x',\boldsymbol y'))\\
&\qquad-
(M_{p,p'}(\boldsymbol x)-N_{q,q'}(\boldsymbol y)) ((\boldsymbol x,\boldsymbol y)\from (\boldsymbol x+\boldsymbol{e}^{(r)}_p,\boldsymbol{y}+\boldsymbol{e}^{(s)}_{q'})\from(\boldsymbol x',\boldsymbol y'))
\\
\intertext{and}
&M_{p,p'}(\boldsymbol x)(N_{q,q'}(\boldsymbol y)+1)((\boldsymbol x,\boldsymbol y)\from(\boldsymbol x+\boldsymbol{e}^{(r)}_p,\boldsymbol{y}+\boldsymbol{e}^{(s)}_q)
\from(\boldsymbol x',\boldsymbol y'))\\&\qquad-(M_{p,p'}(\boldsymbol x)+1)N_{q,q'}(\boldsymbol y)((\boldsymbol x,\boldsymbol y)\from(\boldsymbol x+\boldsymbol{e}^{(r)}_{p'},\boldsymbol y+\boldsymbol{e}^{(s)}_{q'})\from(\boldsymbol x',\boldsymbol y'))\\
&\qquad-(M_{p,p'}(\boldsymbol x)-N_{q,q'}(\boldsymbol y))((\boldsymbol x,\boldsymbol y)\from(\boldsymbol x+\boldsymbol{e}^{(r)}_{p'},\boldsymbol{y}+\boldsymbol{e}^{(s)}_q)\from(\boldsymbol x',\boldsymbol y')).
\end{align*}
Finally, if~$M_{p,p'}(\boldsymbol{x})=N_{q,q'}(\boldsymbol{y})$, we have
\begin{align*}
&((\boldsymbol x,\boldsymbol y)\from (\boldsymbol x+\boldsymbol{e}^{(r)}_p,\boldsymbol{y}+\boldsymbol{e}^{(s)}_q)\from (\boldsymbol{x}',\boldsymbol{y}'))
-((\boldsymbol x,\boldsymbol y)\from(\boldsymbol x+\boldsymbol{e}^{(r)}_{p'},\boldsymbol y+\boldsymbol{e}^{(s)}_{q'})\from(\boldsymbol x',\boldsymbol y'))\\
\intertext{and}
&2((\boldsymbol x,\boldsymbol y)\from (\boldsymbol x+\boldsymbol{e}^{(r)}_p,\boldsymbol{y}+\boldsymbol{e}^{(s)}_q)\from (\boldsymbol{x}',\boldsymbol{y}'))\\
&\qquad-M_{p,p'}(\boldsymbol{x})((\boldsymbol x,\boldsymbol y)\from (\boldsymbol x+\boldsymbol{e}^{(r)}_p,\boldsymbol{y}+\boldsymbol{e}^{(s)}_{q'})\from (\boldsymbol{x}',\boldsymbol{y}'))\\&\qquad
+(M_{p,p'}(\boldsymbol{x})+2)((\boldsymbol x,\boldsymbol y)\from(\boldsymbol x+\boldsymbol{e}^{(r)}_{p'},\boldsymbol{y}+\boldsymbol{e}^{(s)}_q)\from(\boldsymbol x',\boldsymbol y')).
\end{align*}
Note that the coefficients in these relations, and in particular their genericity, depend on a family 
of $(r+s)$ positive integer parameters $z_p^-(\lambda)$, $z_q^+(\lambda)$, $1\le p\le r$, $1\le q\le s$ which are independent of~$\boldsymbol{m}$,
$\boldsymbol{n}$. The resulting algebra is Koszul, has global dimension at most~$rs$ and is finite dimensional if and only if~$i_1>1$ and~$j_s<\ell$.

\section{Type~\texorpdfstring{$C_\ell$, $\ell\ge 2$}{C\_l}}\label{C}

\subsection{}\label{C0}
Let $\beta_{i,j}=\beta_{j,i}=\alpha_{i,\ell-1}+\alpha_{j,\ell-1}+\alpha_\ell$,
$1\le i\le j<\ell$ and~$\beta_{\ell,\ell}=\alpha_\ell$. In particular, $\beta_{1,1}=\theta$.
The roots~$\alpha_{i,j}$ and~$\beta_{i,j}$, $i<j$ are short 
while the roots~$\beta_{i,i}$, $i\in I$ are long and 
$$
R^+=\{ \alpha_{i,j}\,:\, 1\le i\le j<\ell\}\cup\{\beta_{i,j}\,: i\le j\in I\}.
$$
In terms of fundamental weights, $\beta_{i,j}=\varpi_i+\varpi_j-\varpi_{i-1}-\varpi_{j-1}$,
where we set as before~$\varpi_0=0$. 

Let~$\Psi$ be an extremal set of positive roots.
Our first observation is that~$\alpha_{i,j}\notin\Psi$ for all~$1\le i\le j<\ell$ since $2\alpha_{i,j}=\beta_{i,i}-\beta_{j+1,j+1}$ and 
hence if~$\alpha_{i,j}\in\Psi$ we get a contradiction by~\corref{*}.
Furthermore, since~$2\beta_{i,j}=\beta_{i,i}+\beta_{j,j}$ we conclude by~\corref{*} that~$\beta_{i,j}\in\Psi$ if and only if~$\beta_{i,i},\beta_{j,j}\in\Psi$.
From this observation, it is immediate that all extremal sets in~$R^+$ are of the form $\Psi(i_1,\dots,i_k):=\{ \beta_{i_r,i_s}\,:\, 1\le r\le s\le k\}$,
$1\le i_1<\cdots<i_k\le \ell$, $1\le k\le \ell$ (see also~\cite{CRD}).

Since~$\varepsilon(\beta_{i,j})=\varpi_{i-1}+\varpi_{j-1}$, we immediately obtain the following
\begin{lem}
Let~$\beta_{i,j}\in \Psi$, $i<j\in I$. Then for all~$\lambda\in P^+$,
$\lambda\from\lambda+\beta_{i,j}\in(\Delta_\Psi)_1$ if and only if $\lambda(h_{i-1}),\lambda(h_{j-1})>0$.
Furthermore, $\lambda\from\lambda+\beta_{i,i}\in(\Delta_\Psi)_1$ if and only if~$\lambda(h_{i-1})>1$.\qed
\end{lem}

\def\cref#1{(C\ref{#1})}
\subsection{}\label{C10}
We now proceed to describe all paths of length~$2$ in~$\Delta_\Psi$. Let~$\eta\in\Psi+\Psi$.
It follows from~\ref{C0} that apart from the trivial case~$\eta=2\beta_{i,i}$, $\beta_{i,i}\in\Psi$, we have four cases to consider.

\begin{enumerate}[(C1)]
\item \label{C10.c2} Assume that~$i<j\in I$. Let~$\eta=\beta_{i,i}+\beta_{i,j}$. Then~${\Delta_\Psi}(\lambda,\lambda+\eta)=\emptyset$
unless~$\lambda(h_{i-1})>2$ and 
$$
{\Delta_\Psi}(\lambda,\lambda+\eta)=
\begin{cases}
\{ (\lambda\from\lambda+\beta_{i,i}\from\lambda+\eta), (\lambda\from\lambda+\beta_{i,j}\from\lambda+\eta)\},& \lambda(h_{j-1})>0\\
\{ (\lambda\from\lambda+\beta_{i,i}\from\lambda+\eta)\}, &i=j-1,\,\lambda(h_{j-1})=0.
\end{cases}
$$
Similarly, if~$\eta=\beta_{i,j}+\beta_{j,j}$, ${\Delta_\Psi}(\lambda,\lambda+\eta)=\emptyset$ unless~$\lambda(h_{j-1})>1$ and~$\lambda(h_{i-1})>0$. Then
$$
{\Delta_\Psi}(\lambda,\lambda+\eta)=
\begin{cases}
\{ (\lambda\from\lambda+\beta_{i,j}\from\lambda+\eta), (\lambda\from\lambda+\beta_{j,j}\from\lambda+\eta)\},& \lambda(h_{j-1})>2\\
\{ (\lambda\from\lambda+\beta_{i,j}\from\lambda+\eta)\}, &i=j-1,\,\lambda(h_{j-1})=2.
\end{cases}
$$

\item\label{C10.c3} Let~$\eta=\beta_{i,i}+\beta_{i,j}=2\beta_{i,j}$, $i\in j\in I$. Then~${\Delta_\Psi}(\lambda,\lambda+\eta)=\emptyset$ unless~$\lambda(h_{i-1})>1$. 
If~$\lambda(h_{j-1})>1$ then ${\Delta_\Psi}(\lambda,\lambda+\eta)$ contains all three possible paths. Otherwise, 
${\Delta_\Psi}(\lambda,\lambda+\eta)$ is empty unless~$i=j-1$. If~$i=j-1$ we have
$$
{\Delta_\Psi}(\lambda,\lambda+\eta)=\begin{cases}
                                            \{ (\lambda\from\lambda+\beta_{i,i}\from\lambda+\eta),(\lambda\from\lambda+\beta_{i,j}\from\lambda+\eta)\},&
\lambda(h_{j-1})=1,\\
\{ (\lambda\from\lambda+\beta_{i,i}\from\lambda+\eta)\},&\lambda(h_{j-1})=0.
                                           \end{cases}
$$

\item\label{C10.c4} Assume that~$i<j<k$. First, let~$\eta=\beta_{i,i}+\beta_{j,k}=\beta_{i,j}+\beta_{i,k}$. Then~${\Delta_\Psi}(\lambda,\lambda+\eta)=\emptyset$ unless
$\lambda(h_{i-1})>1$ and~$\lambda(h_{j-1})+\lambda(h_{k-1})>0$.
If~$\lambda(h_{j-1}),\lambda(h_{k-1})>0$ then we have all four possible paths. Otherwise,
$$
{\Delta_\Psi}(\lambda,\lambda+\eta)=\begin{cases}
\{ (\lambda\from\lambda+\beta_{i,i}\from\lambda+\eta),(\lambda\from\lambda+\beta_{i,k}\from\lambda+\eta)\},& i=j-1,\,\lambda(h_{j-1})=0,\\
\{ (\lambda\from\lambda+\beta_{i,j}\from\lambda+\eta)\},& j=k-1,\lambda(h_{k-1})=0.
\end{cases}
$$
Next, let~$\eta=\beta_{j,j}+\beta_{i,k}=\beta_{i,j}+\beta_{j,k}$. Then~${\Delta_\Psi}(\lambda,\lambda+\eta)=\emptyset$ unless
$\lambda(h_{i-1})>0$ and~$\lambda(h_{j-1})+\lambda(h_{k-1})>1$. If~$\lambda(h_{j-1})>1$ and~$\lambda(h_{k-1})>0$ then
we have all possible paths. Otherwise,
$$
{\Delta_\Psi}(\lambda,\lambda+\eta)=\begin{cases}
\{ (\lambda\from\lambda+\beta_{i,k}\from\lambda+\eta)\},& i=j-1,\lambda(h_{j-1})=1,
\\
\{ (\lambda\from\lambda+\beta_{j,j}\from\lambda+\eta),(\lambda\from\lambda+\beta_{i,j}\from\lambda+\eta)\},& j=k-1,\lambda(h_{k-1})=0.
\end{cases}
$$
Finally, if~$\eta=\beta_{k,k}+\beta_{i,j}=\beta_{i,k}+\beta_{j,k}$, then ${\Delta_\Psi}(\lambda,\lambda+\eta)=\emptyset$ unless
$\lambda(h_{i-1})>0$ and~$\lambda(h_{j-1})+\lambda(h_{k-1})>1$. If~$\lambda(h_{k-1})>1$ and~$\lambda(h_{j-1})>0$ then we have all possible paths.
Otherwise, 
$$
{\Delta_\Psi}(\lambda,\lambda+\eta)=\begin{cases}
\{ (\lambda\from\lambda+\beta_{i,k}\from\lambda+\eta)\},& i=j-1,\lambda(h_{j-1})=0,
\\
\{ (\lambda\from\lambda+\beta_{i,j}\from\lambda+\eta),(\lambda\from\lambda+\beta_{j,k}\from\lambda+\eta)\},& j=k-1,\lambda(h_{k-1})=1.
\end{cases}
$$

\item\label{C10.c5} Finally, let~$i<j<k<l\in I$, $\eta=\beta_{i,j}+\beta_{k,l}=\beta_{i,k}+\beta_{j,l}=\beta_{i,l}+\beta_{j,k}$.
Then~${\Delta_\Psi}(\lambda,\lambda+\eta)=\emptyset$ unless
$$
\lambda(h_{i-1}),\lambda(h_{j-1})+\lambda(h_{k-1}),\lambda(h_{k-1})+\lambda(h_{l-1})>0.
$$
If~$\lambda(h_{r-1})>0$, $r\in\{i,j,k,l\}$ we have all possible paths. Furthermore, if~$\lambda(h_{k-1})=0$ we must have~$j=k-1$ and 
$$
{\Delta_\Psi}(\lambda,\lambda+\eta)=\{ (\lambda\from\lambda+\beta_{i,j}\from\eta),(\lambda\from\lambda+\beta_{j,l}\from\lambda+\eta)\}.
$$
Finally, if~$\lambda(h_{k-1})>0$ we have
$$
{\Delta_\Psi}(\lambda,\lambda+\eta)=\left\{\setlength\arraycolsep{-2.5pt}
\begin{array}{ll}
\{ (\lambda\from\lambda+\beta_{i,k}\from\lambda+\eta),&(\lambda\from\lambda+\beta_{i,l}\from\lambda+\eta)\},\\&\quad i=j-1,\lambda(h_{j-1})=0,\lambda(h_{l-1})>0\\
\{ (\lambda\from\lambda+\beta_{i,k}\from\lambda+\eta),&(\lambda\from\lambda+\beta_{j,k}\from\lambda+\eta)\},\\&\quad k=l-1,\lambda(h_{j-1})>0,\lambda(h_{l-1})=0\\
\{ (\lambda\from\lambda+\beta_{i,k}\from\lambda+\eta)\},& \quad k=l-1,i=j-1,\lambda(h_{j-1})=\lambda(h_{l-1})=0.
\end{array}\right.
$$
\end{enumerate}
The cases~\cref{C10.c2}--\cref{C10.c3} (respectively, \cref{C10.c4}, \cref{C10.c5}) occur if $\Psi\supset \Psi(i,j)$
(respectively, $\Psi\supset\Psi(i,j,k)$, $\Psi\supset\Psi(i,j,k,l)$).
In particular, we obtain the following
\begin{lem}
The set~$\Psi(i_1,\dots,i_k)$ is regular if and only if~$i_{r+1}\not=i_r+1$ for all~$1\le r<k$. 
\end{lem}

\subsection{}\label{C2}
Retain the notations of~\ref{M95}. A straightforward
induction on~$r$ shows that 
\begin{equation}\label{C2.0}
\#\Xi_0(\boldsymbol{m})_0=\left\lceil \frac12 (m_1+1)\cdots(m_r+1)\right\rceil,\qquad \#\Xi_1(\boldsymbol{m})_0=
\left\lfloor \frac12(m_1+1)\cdots(m_r+1)\right\rfloor.
\end{equation}
It is immediate that $\Xi_a(\boldsymbol{m})\cong \Xi_a(\boldsymbol{m}')$ if $\boldsymbol{m}'$ is a permutation of~$\boldsymbol{m}$
or is obtained from~$\boldsymbol{m}$ by adding or removing zeroes. Clearly, $\Xi_a((m))\cong \Xi_{a'}((m'))$, $a,a'\in\{0,1\}$ if and only if $\lfloor(m-a)/2\rfloor=
\lfloor (m'-a')/2\rfloor$. Note also that~$\Xi_{0}((1,1))\cong \Xi_0((2))\cong \Xi_0((3))\cong \Xi_{1}((3))\cong \Xi_1((4))$.
\begin{prop}
Let~$\boldsymbol{m}=(m_1,\dots,m_r)\in\bz_+^r$, $m_1\ge \cdots \ge m_r>0$, $r>1$.
The quivers 
$\Xi_0(\boldsymbol{m})$ and $\Xi_1(\boldsymbol{m})$, $\boldsymbol{m}\not=(1,1)$ are connected and pairwise non-isomorphic.
Furthermore, $\Xi_0(\boldsymbol{m})\cong \Xi_1(\boldsymbol{m})^{op}$ if and only if~$|\boldsymbol{m}|$ is odd.
\end{prop}
\begin{pf}
Observe that every vertex in~$\Xi_a(\boldsymbol{m})$ is connected to a sink. 
Clearly, $\boldsymbol{0}=(0,\dots,0)$ is the unique sink in $\Xi_0(\boldsymbol{m})$, hence~$\Xi_0(\boldsymbol{m})$ is connected. On the other hand,
the $\boldsymbol{e}_j:=\boldsymbol{e}_j^{(r)}$ are the only sinks in~$\Xi_1(\boldsymbol{m})$. 
If~$r=2$ and~$\boldsymbol{m}\not=(1,1)$, then $m_1>1$ and so we have $\boldsymbol{e}_1\from 2\boldsymbol{e}_1+\boldsymbol{e}_2\to
\boldsymbol{e}_2$. If~$r>2$ then for all $1\le i<j<k\le r$, we have
$$
\divide\dgARROWLENGTH by 2
\begin{diagram}
\node{\boldsymbol{e}_i}\node{\boldsymbol{e}_i+\boldsymbol{e}_j+\boldsymbol{e}_k}\arrow{e}\arrow{w}\arrow{s}\node{\boldsymbol{e}_j}\\
\node{}\node{\boldsymbol{e}_k}
\end{diagram} 
$$
Thus, all sinks in~$\Xi_1(\boldsymbol{m})$, $\boldsymbol{m}\not=(1,1)$ lie in the same connected component hence
$\Xi_1(\boldsymbol{m})$ is connected and $\Xi_1(\boldsymbol{m})\not\cong \Xi_0(\boldsymbol{m}')$ for all~$\boldsymbol{m}'=(m_1',\dots,m_k')$,
$m_1'\ge \cdots\ge m'_k>0$, $k>1$.

Given~$\boldsymbol{m}=(m_1,\dots,m_r)$, $m_1\ge \cdots\ge m_r$, let $n_p(\boldsymbol{m})=\#\{j\,:\, m_j=p\}$, $p>0$ and 
$\ell(\boldsymbol{m})=\sum_{p>0} n_p(\boldsymbol{m})$.
Suppose that~$\Xi_0(\boldsymbol m)$ is isomorphic to~$\Xi_0(\boldsymbol{m}')$, $\boldsymbol{m}'=(m_1',\dots,m_k')$, $m_1'\ge \cdots\ge m_k'$. 
Suppose first that~$k=\ell(\boldsymbol{m}')>\ell(\boldsymbol{m})=r$.
Then 
$\#\boldsymbol{0}^-=r(r+1)/2-n_1(\boldsymbol{m})$. Since~$\boldsymbol{0}$ (respectively, $\boldsymbol{0}'$) is the unique 
sink in~$\Xi_0(\boldsymbol{m})$ (respectively, in~$\Xi_0(\boldsymbol{m}')$),
we must have~$n_1(\boldsymbol{m}')=k(k+1)/2-r(r+1)/2+n_1(\boldsymbol{m})$. Since~$n_1(\boldsymbol{m}')\le k$, this implies that
$n_1(\boldsymbol{m})\le r(r+1)/2-k(k-1)\le 0$ with the equality if and only if~$k=r+1$ which in turn implies that~$n_1(\boldsymbol{m}')=k$. 
Since~$r>1$, by~\eqref{C2.0} we obtain
$\#\Xi(\boldsymbol{m})_0\ge 3^r/2>2^r=\#\Xi(\boldsymbol{m}')_0$ which is a contradiction. Thus,
$k=r$ and $n_1(\boldsymbol{m})=n_1(\boldsymbol{m}')$.

Furthermore, note that $\boldsymbol{x}\in\Xi(\boldsymbol{m})_0$ is a source if and only if $|\boldsymbol{x}| \ge |\boldsymbol{m}|-1$. It follows that
if $|\boldsymbol{m}|\in2\bz_+$ then~$\boldsymbol{m}$ is the unique source 
$\Xi_0(\boldsymbol{m})$. Otherwise, the $\boldsymbol{m}-\boldsymbol{e}_j$, $1\le j\le r$ are sources. Therefore, $|\boldsymbol{m}|=|\boldsymbol{m}'|\pmod 2$.
Since the length of any path in~$\Xi_0(\boldsymbol{m})$ from a source to the unique sink is $\lfloor |\boldsymbol{m}|/2\rfloor$, it follows that~$|\boldsymbol{m}|=
|\boldsymbol{m}'|$. Clearly, if~$n_1(\boldsymbol{m})\ge r-1$ or~$m_j\le 2$ for all~$1\le j\le r$, this implies that  $\boldsymbol{m}=\boldsymbol{m}'$, so we may assume that~$n_1(\boldsymbol{m})<r-1$ and~$m_j>2$ for some~$1\le j\le r$.

Note that if~$\boldsymbol{x}\in\Xi(\boldsymbol{m})_0$ satisfies
$\#\boldsymbol{x}^+ =1$ then either $\boldsymbol{x}=2x\boldsymbol{e}_j$ for some~$1\le j\le r$, $0\le x\le m_j/2$, or
$\boldsymbol{x}=\boldsymbol{e}_j+\boldsymbol{e}_k$, $1\le j<k\le r$. 
Given~$1\le j\le r$ with~$m_j>2$ consider a path
$$
\boldsymbol{0}\from 2 \boldsymbol{e}_j\from \cdots \from 2k \boldsymbol{e}_j,\qquad k=\lfloor m_j/2\rfloor
$$  
in~$\Xi_0(\boldsymbol{m})$. Then its image in~$\Xi_0(\boldsymbol{m}')$ under our isomorphism of quivers  must be
$$
\boldsymbol{0}\from 2\boldsymbol{e}_{j'}\from\cdots\from 2k \boldsymbol{e}_{j'},
$$
for some~$1\le j'\le r$ with~$\lfloor m_{j'}/2\rfloor=\lfloor m_j/2\rfloor$.
Furthermore, it is easy to check that 
\begin{equation}\label{C2.10}
\boldsymbol{x}\in\Xi_0(\boldsymbol{m})_0,\, \#\boldsymbol{x}^+\le 3\implies\boldsymbol{x}=x_i \boldsymbol{e}_i+
x_j\boldsymbol{e}_j,\qquad 1\le i<j\le r.
\end{equation}

Suppose first that~$r=2$ and~$n_1(\boldsymbol{m})=0$. Since~$m_1+m_2=m_1'+m_2'$, we may assume, without loss of generality, that~$m_1>m_1'$. 
By the above, we must have~$\lfloor m_1/2\rfloor=\lfloor m_1'/2\rfloor$
hence~$m_1=2a+1$, $m_1'=2a$, $a\ge 1$ and so~$m_2=2b$, $m_2'=2b+1$, $b\ge 1$. Since~$\#\Xi_0((m_1,m_2))_0=
\#\Xi_0((m_1',m_2'))_0$, we conclude that~$a=b$, which is a contradiction since~$m_1'\ge m_2'$. 

Suppose now that~$r>2$. For~$1\le i<j\le s$ fixed let $\Xi^{i,j}_0(\boldsymbol{m})$ be the full subquiver of~$\Xi_0(\boldsymbol{m})$ defined 
by $\{ x_i \boldsymbol{e}_i+x_j\boldsymbol{e}_j\,:\, x_i\le m_i,\,x_j\le m_j,\,x_i+x_j\in2\bz_+\}$.
Clearly for all~$\boldsymbol{x}\in \Xi^{i,j}_0(\boldsymbol{m})_0$ the set of direct successors of~$\boldsymbol{x}$ in~$\Xi_0(\boldsymbol{m})$ 
is contained in $\Xi^{i,j}_0(\boldsymbol{m})_0$, hence~$\Xi_0^{i,j}(\boldsymbol{m})$
is a convex connected subquiver of~$\Xi_0(\boldsymbol{m})$. It is clearly isomorphic to~$\Xi_0((m_i,m_j))$.
It follows from~\eqref{C2.10} that the isomorphism of quivers $\Xi_0(\boldsymbol{m})\to\Xi_0(\boldsymbol{m}')$ induces an isomorphism of quivers
$\Xi_0^{i,j}(\boldsymbol{m})\to \Xi_0^{i',j'}(\boldsymbol{m}')$ for 
some~$1\le i'<j'\le r$ which by the~$r=2$ case implies that~$m_i=m_{i'}'$, $m_j=m_{j'}'$. Therefore, $\boldsymbol{m}=\boldsymbol{m}'$.

Suppose that~$\Xi_1(\boldsymbol{m})\cong \Xi_1(\boldsymbol{m}')$. Since $\Xi_1(\boldsymbol{m})$ contains $\ell(\boldsymbol{m})$ sinks,
it follows that~$\ell(\boldsymbol{m})=\ell(\boldsymbol{m}')=r$. Furthermore, we have 
$$
|\boldsymbol{e}_i^-|=\begin{cases}\binom{r}{2}-n_1(\boldsymbol{m}),& m_i=1\\
                      \binom{r}{2}+r-1-n_1(\boldsymbol{m}),& m_i=2.
			\\
			\binom{r}{2}+r-n_1(\boldsymbol{m}),&m_i>2.
                     \end{cases}
$$
It follows that~$n_p(\boldsymbol{m})=n_p(\boldsymbol{m}')$, $p=1,2$. Since~$\Xi_1(\boldsymbol{m})$ contains a unique source if~$|\boldsymbol{m}|$ is odd
and $r$ sources otherwise, it follows that~$|\boldsymbol{m}|=|\boldsymbol{m}'| \pmod 2$. Since the length of a path from a source to a sink is
$(|\boldsymbol{m}|-1)/2$ if~$|\boldsymbol{m}|$ is odd and~$|\boldsymbol{m}|/2-1$ if~$|\boldsymbol{m}|$ is even, it follows that~$|\boldsymbol{m}|=|\boldsymbol{m}'|$.
Furthermore, note that $\boldsymbol{x}\in\Xi_1(\boldsymbol{m})$, $\#\boldsymbol{x}^+\le 3$ implies that~$\boldsymbol{x}\in\Xi^{i,j}_1(\boldsymbol{m})$
for some~$1\le i<j\le r$ or~$\boldsymbol{x}=\boldsymbol{e}_i+\boldsymbol{e_j}+\boldsymbol{e}_k$, $1\le i<j<k\le r$. On the other hand, a vertex of the second
type is connected to three sinks in~$\Xi_1(\boldsymbol{m})$ by arrows, while a vertex of the first type can be connected to at most two sinks. Thus, we conclude
as before that the image of~$\Xi^{i,j}_1(\boldsymbol{m})$ under the isomorphism $\Xi_1(\boldsymbol{m})\to\Xi_1(\boldsymbol{m}')$ is contained 
in~$\Xi_1^{i',j'}(\boldsymbol{m}')$ for some~$1\le i'<j'\le r$.
The rest of the argument is similar to that in the ``even'' case and is omitted.

To prove the last assertion, note that if~$|\boldsymbol{m}|$ is odd, then at least one of the~$m_r$ is odd, hence $\#\Xi_0(\boldsymbol{m})_0=\#\Xi_1(\boldsymbol{m})_0$
and the map $\Xi_0(\boldsymbol{m})_0\to \Xi_1(\boldsymbol{m})_0$, $\boldsymbol{x}\mapsto \boldsymbol{m}-\boldsymbol{x}$, is a bijection. This map induces
the desired isomorphism of quivers. Conversely, if~$|\boldsymbol{m}|$ is even, then $\Xi_1(\boldsymbol{m})$ contains~$\ell(\boldsymbol{m})>1$ sources.
Since~$\Xi_0(\boldsymbol{m})$ has a unique sink, $\Xi_0(\boldsymbol{m})$ and $\Xi_1(\boldsymbol{m})^{op}$ cannot be isomorphic.
\end{pf}

\subsection{}\label{C3}\label{C5}
We can now describe all connected components of~$\Delta_\Psi$ for~$\Psi$ regular.
\begin{prop}
Let~$\Psi=\Psi(i_1,\dots,i_k)$, $1\le i_1<\cdots<i_k\le \ell$ and suppose that~$\Psi$ is regular.
Let~$\lambda\in P^+$ and assume that~$|\lambda^+\cup\lambda^-|>0$. Then~$\Delta_\Psi[\lambda]$ 
is isomorphic to the quiver $\Xi_a(\boldsymbol{m})$ where 
$\boldsymbol{m}=(m_1,\dots,m_k)\in(\bz_+\cup\{+\infty\})^k$, $m_r=\lambda(h_{i_r-1})+\lambda(h_{i_r})$, $1\le r\le k$
and $a=\lambda(h_{i_1})+\cdots+\lambda(h_{i_k})\pmod2$.
\end{prop}

\begin{pf}
Let~$J=\{i_r\,:\, 1\le r\le n\}\cup \{ i_r-1\,:\, 1\le r\le n\}$.
Suppose that~$\mu\in\Delta_\Psi[\lambda]_0$. Since~$\Delta_\Psi[\lambda]_0\subset (\lambda+\bz\Psi)\cap P^+$, we have
\begin{align*}
\mu(h_{i_r})&=\lambda(h_{i_r})+\sum_{s=1}^{r-1} x_{s,r}+2 x_{r,r}+\sum_{s=r+1}^{k} x_{r,s},\\
\mu(h_{i_r-1})&=\lambda(h_{i_r-1})- \sum_{s=1}^{r-1} x_{s,r}-2 x_{r,r}-\sum_{s=r+1}^{k} x_{r,s},
\end{align*}
where~$x_{p,q}\in\bz$, $1\le p\le q\le k$. It follows that
\begin{subequations}
\begin{gather}
\mu(h_{i_r-1})+\mu(h_{i_r})=\lambda(h_{i_r-1})+\lambda(h_{i_r}),\qquad 1\le r\le k,\qquad \mu(h_j)=\lambda(h_j),\,j\notin J\label{C5.10a}\\
\sum_{r=1}^k \mu(h_{i_r})=\sum_{r=1}^k \lambda(h_{i_r})\pmod 2.\label{C5.10b}
\end{gather}
\end{subequations}
Let~$S(\lambda)$ be the set of~$\mu\in P^+$ satisfying these conditions. 
Then~$\Delta_\Psi[\lambda]_0\subset S(\lambda)$ and for all $\mu\in S(\lambda)$, $\mu^-\subset S(\lambda)$.
Thus, $S(\lambda)$ defines
a convex subquiver~$\Gamma$ of~$\Delta_\Psi$ with~$\Gamma_0=S(\lambda)$ and~$\Delta_\Psi[\lambda]$ is a full connected subquiver of~$\Gamma$.

Let~$m_r=\lambda(h_{i_r-1})+\lambda(h_{i_r})$, $a=\lambda(h_{i_1})+\cdots+\lambda(h_{i_k})\pmod 2$. Then
we have a bijective map
\begin{align*}
\Gamma_0&\longrightarrow\Xi_a(\boldsymbol{m})_0 \\
\mu&\longmapsto (\mu(h_{i_1}),\dots,\mu(h_{i_k})).
\end{align*}
It is easy to see that this induces an isomorphism of quivers~$\Gamma\to\Xi_a(\boldsymbol{m})$. To complete the argument, observe that
the assumption that $|\lambda^+\cup\lambda^-|>0$ implies that we cannot have~$a=1$ and $m_r=\delta_{r,p}+\delta_{r,q}$ for some $1\le p<q\le k$.
Then~$\Xi_a(\boldsymbol{m})$ is connected by~\propref{C3}. Therefore, $\Gamma$ is connected hence $\Delta_\Psi[\lambda]=\Gamma$.
\end{pf}

\subsection{}\label{C50}
Fix root vectors~$e_{\beta_{i,j}}\in\lie g_{\beta_{i,j}}\setminus\{0\}$, $1\le i\le j\le\ell$ 
so that
\begin{equation}
\label{C50.10}
\begin{split}
&[e_i,e_{\beta_{j,k}}]=\delta_{i,j-1} e_{\beta_{i,k}}+\delta_{i,k-1}(1+\delta_{i,j}) e_{\beta_{j,i}},\qquad j<k
\\
&[e_i,e_{\beta_{j,j}}]=\delta_{i,j-1} e_{\beta_{i,j}}
\end{split}
\end{equation}
and
\begin{equation}
\label{C50.20}
\begin{split}
&[f_i,e_{\beta_{j,k}}]=\delta_{i,j}(1+\delta_{j+1,k}) e_{\beta_{j+1,k}}+
\delta_{i,k}e_{\beta_{j,k+1}},\qquad j<k
\\
&[f_i,e_{\beta_{j,j}}]=\delta_{i,j} e_{\beta_{j,j+1}}.
\end{split}
\end{equation}
For example, we can use the standard presentation of~$\lie g$ as the matrix Lie algebra~$\lie{sp}_{2\ell}$.
The subalgebra~$\lie g_J$ of~$\lie g$ with~$J=I\setminus\{\ell\}$ is of course a simple Lie algebra
of type~$A_{\ell-1}$. Note  that~$[e_\ell,e_{\beta_{i,j}}]=0=
[f_\ell,e_{\beta_{i,j}}]$, $1\le i\le j\le\ell$. Due to this observation, we 
can perform our computations in~$ U(\lie g_{J})$.

\subsection{}\label{C60}
Retain the notations of~\ref{A50}.
Fix~$1\le i\le j<\ell$.
Given any pair~$1\le r\le s$ such that
$r\le i+1$, $s\le j+1$ set
$$
\cal U_{r,s,i,j}=e_{s-1}\cdots e_r \cal X^-_{r,i}\cal X^-_{r,j}.
$$
In particular, $\cal U_{r,r,i,j}=\cal X^-_{r,i}\cal X^-_{r,j}\in  U(\lie b)$. 
Clearly, $\cal U_{r,s,i,j}\in  U(\lie g)_{-\alpha_{r,i}-\alpha_{s,j}}$. We set~$\cal U_{r,s,i,j}=0$ if~$r>i+1$.
\begin{lem}
The elements~$\cal U_{r,s,i,j}$ satisfy
\begin{subequations}
\begin{align}
&e_s \,\cal U_{r,s,i,j}=\cal U_{r,s+1,i,j},\label{C60.10b}\\
&e_r\,\cal U_{r,s,i,j}=(1+\delta_{r+1,s})\cal U_{r+1,s,i,j}(\cal H_{r,i}-\delta_{i,j})\cal H_{r,j}+ U(\lie g)\lie n^+,\label{C60.10c}\qquad r<s,\,\\
&e_k \,\cal U_{r,s,i,j}\in  U(\lie g)\lie n^+,\qquad k\not=r,s.\label{C60.10a}
\end{align}
\end{subequations}
\end{lem}
\begin{pf}
The first identity is obvious. 
To prove~\eqref{C60.10c} and~\eqref{C60.10a}, we need to show first that
\begin{subequations}
\begin{align}
&e_k \cal U_{r,r,i,j}\in  U(\lie g)e_k,\qquad k\not=r\label{C60.20a}\\
&e_r^2 \cal U_{r,r,i,j}=2 \cal U_{r+1,r+1,i,j} (\cal H_{r,i}-\delta_{i,j})\cal H_{r,j}+ U(\lie g)e_r
\label{C60.20b}\\
&e_{r-1}e_r \cal U_{r,r,i,j}\in  U(\lie g)\lie n^+.\label{C60.20c}
\end{align}
\end{subequations}
Using~\lemref{A50} we immediately deduce~\eqref{C60.20a} and the following identity 
\begin{equation}\label{C60.25}
\cal U_{r,r+1,i,j}=e_r \cal U_{r,r,i,j}=\cal X^-_{r+1,i}\cal X^-_{r,j}(\cal H_{r,i}-1-\delta_{i,j})+
\psi_{\alpha_r}(\cal X^-_{r,i})\cal X^-_{r+1,j}\cal H_{r,j}+ U(\lie g)e_r.
\end{equation}
Then~\eqref{C60.20b} and~\eqref{C60.20c} are easy to obtain using~\eqref{A50.0a}.
To prove~\eqref{C60.10c} for~$s>r$, note that the case~$s=r+1$ is immediate from~\eqref{C60.20b}.
Assume that~$s>r+1$.
Clearly $e_r$ commutes with the~$e_t$, $r+1<t\le s-1$. Since~$(\ad e_a)^2 e_b=0$
for all $1\le a\not=b<\ell$ with~$|a-b|=1$, we have in~$ U(\lie g)$
\begin{equation}\label{C60.30}
e_a^2 e_b^{}-2 e_a^{} e_{b}^{} e_a^{}+e_b^{}e_a^2=0.
\end{equation}
Therefore,
$$
e_r\cal U_{r,s,i,j}=e_{s-1}\cdots e_{r+2} e_r e_{r+1} e_r\cal U_{r,r,i,j}=
e_{s-1}\cdots e_{r+1} \cal U_{r+1,r+1,i,j}(\cal H_{r,i}-\delta_{i,j})\cal H_{r,j}+ U(\lie g)\lie n^+,
$$
where we used~\eqref{C60.20a} and~\eqref{C60.20b}. 
To prove~\eqref{C60.10a} for~$s>r$, note that for
$k<r-1$ or~$k>s$ this is an immediate consequence of~\eqref{C60.20a}. Thus, if~$s=r+1$ there is nothing to do.
Assume that~$s>r+1$.
If~$k=r-1$, the assertion
follows from~\eqref{C60.20c}. If~$k=s-1$, it follows from~\eqref{C60.30} that
$$
e_{s-1} \cal U_{r,s,i,j}=e_{s-1}^2 e_{s-2}\cdots e_r\cal U_{r,r,i,j}
=-e_{s-2}\cdots e_r e_{s-1}^2 \cal U_{r,r,i,j}+2 e_{s-1}\cdots e_r e_{s-1} \cal U_{r,r,i,j}
$$
which is contained in~$ U(\lie g)\lie n^+$ by~\eqref{C60.20a}. If~$r<k<s-1$ we can write, using~\eqref{C60.30}
\begin{multline*}
e_k\cal U_{r,s,i,j}=e_{s-1}\cdots e_k e_{k+1} e_ke_{k-1}\cdots e_r \cal U_{r,r,i,j}=
\frac12 e_{s-1}\cdots e_{k+1} e_k^2 e_{k-1}\cdots e_r\cal U_{r,r,i,j}\\+
\frac12 e_{s-1}\cdots e_k^2 e_{k+1} e_{k-1}\cdots e_r\cal U_{r,r,i,j}.
\end{multline*}
The second term is in~$ U(\lie g)\lie n^+$ by~\eqref{C60.20a} since~$e_{k+1}$ commutes with the~$e_t$,
$t\le k-1$. Applying~\eqref{C60.30} again, we obtain
$$
e_k\cal U_{r,s,i,j}=\Big(e_{s-1}\cdots e_r e_k-\frac12 e_{s-1}\cdots e_{k+1} e_{k-1}\cdots e_r e_k^2\Big) \cal U_{r,r,i,j}+ U(\lie g)\lie n^+\in  U(\lie g)\lie n^+,
$$
where we used~\eqref{C60.20a}.
\end{pf}
\subsection{}\label{C115}
For our purposes, we need to find the projection $\bar{\cal U}_{r,s,i,j}$ of~$\cal U_{r,s,i,j}$ onto~$U(\lie b)$.
\begin{lem}
Let~$1\le i\le j<\ell$ and suppose that~$r<s$, $s\le j$, $r\le i$. Then
\begin{equation}\label{C115.10}
\begin{split}
\bar{\cal U}_{r,s,i,j}=\cal X^-_{s,i}\cal X^-_{r,j}\prod_{t=r}^{s-1} (\cal H_{t,i}-\delta_{r,t}-\delta_{i,j})+\psi_{\alpha_{r,s-1}}(\cal X^-_{r,i})
\cal X^-_{s,j}\prod_{t=r}^{s-1} \cal H_{t,j}\\
-\cal X^-_{s,i}\sum_{t=r+1}^{s-1} \psi_{\alpha_{r,t-1}}(\cal X^-_{r,t-1,i})\cal X^-_{t,j}\prod_{p=r}^{t-1}\cal H_{p,j}\prod_{p=t+1}^{s-1} (\cal H_{p,i}-\delta_{i,j}).
\end{split}
\end{equation}
In particular,
\begin{equation}\label{C115.20}
\bar{\cal U}_{i+1,s,i,j}=\cal X^-_{s,j}\prod_{t=i+1}^{s-1} \cal H_{t,j},\qquad i+1\le s\le j+1
\end{equation}
and $\cal U_{r,s,i,j}(\mu):=\pi_\mu(\bar{\cal U}_{r,s,i,j})$ is given by the following formulae
\begin{align*}
\cal U_{r,r,i,j}&(\mu)=\cal X^-_{r,i}(\mu-\varpi_j)\cal X^-_{r,j}(\mu)
\\
\cal U_{r,s,i,j}&(\mu)=\cal X^-_{s,i}(\mu-\varpi_j)\cal X^-_{r,j}(\mu)\prod_{t=r}^{s-1} (\mu(\cal H_{t,i})-\delta_{r,t}-\delta_{i,j})+
\cal X^-_{r,i}(\mu-\varpi_j)
\cal X^-_{s,j}(\mu)\prod_{t=r}^{s-1} \mu(\cal H_{t,j})\\&
-\cal X^-_{s,i}(\mu-\varpi_j)\sum_{t=r+1}^{s-1} \cal X^-_{r,t-1,i}(\mu-\varpi_j)\cal X^-_{t,j}(\mu)\prod_{p=r}^{t-1}\mu(\cal H_{p,j})\prod_{p=t+1}^{s-1} (\mu(\cal H_{p,i})-\delta_{i,j})
,\quad r<s.
\end{align*}
\end{lem}
\begin{pf}
The elements~$\bar{\cal U}_{r,s,i,j}$ are uniquely determined by the conditions that~$\bar{\cal U}_{r,s,i,j}=\cal U_{r,s,i,j}+U(\lie g)\lie n^+$ and
$\bar{\cal U}_{r,s,i,j}\in U(\lie b)$.
The argument is by induction on~$s-r$, the induction base being~\eqref{C60.25}.
To prove the inductive step, note that by~\lemref{A50} and the induction hypothesis we have
\begin{multline}\label{C115.30}
\cal U_{r,s+1,i,j}=e_s\cal U_{r,s,i,j}=\cal X^-_{s+1,i}\cal X^-_{r,j}\prod_{k=r}^{s} (\cal H_{k,i}-\delta_{r,k}-\delta_{i,j})+
e_{s} \psi_{\alpha_{r,s-1}}(\cal X^-_{r,i})
\cal X^-_{s,j}\prod_{k=r}^{s-1} \cal H_{k,j}\\
-e_s \cal X^-_{s,i}\sum_{t=r+1}^{s-1} \psi_{\alpha_{r,t-1}}(\cal X^-_{r,t-1,i})\cal X^-_{t,j}\prod_{p=r}^{t-1}\cal H_{p,j}\prod_{p=t+1}^{s-1} (\cal H_{p,i}-\delta_{i,j})
+ U(\lie g)\lie n^+.
\end{multline}
Applying~\lemref{A50} to the second term, we obtain
\begin{multline*}
e_s \psi_{\alpha_{r,s-1}}(\cal X^-_{r,i})\cal X^-_{s,j}\prod_{t=r}^{s-1}\cal H_{t,j}=
\psi_{\alpha_{r,s-1}}(\cal X^-_{s+1,i}\cal X^-_{r,s-1,i})(\alpha_{r,s-1})(h_s)\\+
\psi_{\alpha_{r,s}}(\cal X^-_{r,i})e_s \cal X^-_{s,j}\prod_{t=r}^{s-1}\cal H_{t,j}\\
=-\cal X^-_{s+1,i}\psi_{\alpha_{r,s-1}}(\cal X^-_{r,s-1,i})\cal X^-_{s,j}\prod_{t=r}^{s-1}\cal H_{t,j}+\psi_{\alpha_{r,s}}(\cal X^-_{r,i})\cal X^-_{s+1,j}\prod_{t=r}^s\cal H_{t,j}+ U(\lie g)\lie n^+,
\end{multline*}
where we noted that~$\psi_{\alpha_{r,s-1}}(\cal X^-_{s+1,i})=\cal X^-_{s+1,i}$. Finally, the last term in~\eqref{C115.30} can be written as follows
\begin{multline*}
-\cal X^-_{s+1,i}\cal H_{s,i} \sum_{t=r+1}^{s-1} \psi_{\alpha_{r,t-1}}(\cal X^-_{r,t-1,i})\cal X^-_{t,j}\prod_{p=r}^{t-1}\cal H_{p,j}\prod_{p=t+1}^{s-1} (\cal H_{p,i}-\delta_{i,j})\\
-\psi_{\alpha_s}(\cal X^-_{s,i}) \sum_{t=r+1}^{s-1}e_s \psi_{\alpha_{r,t-1}}(\cal X^-_{r,t-1,i})\cal X^-_{t,j}\prod_{p=r}^{t-1}\cal H_{p,j}\prod_{p=t+1}^{s-1}
(\cal H_{p,i}-\delta_{i,j}).
\end{multline*}
Since~$\cal X^-_{r,t-1,i}=\sum_{\tau\in\Sigma(r,t-1)} \bof_\tau c_{\tau}^-(i)$, it follows that~$e_s\psi_{\eta}(\cal X^-_{r,t-1,i})
=\psi_{\eta+\alpha_s}(\cal X^-_{r,t-1,i})e_s$ and so we get
$$
\psi_{\alpha_s}(\cal X^-_{s,i}) e_s\psi_{\alpha_{r,t-1}}(\cal X^-_{r,t-1,i})\cal X^-_{t,j}=
\psi_{\alpha_{r,t-1}+\alpha_s}(\cal X^-_{r,t-1,i})\psi_{\alpha_s}(\cal X^-_{t,j})e_s\in  U(\lie g)\lie n^+.
$$
Thus the last term in~\eqref{C115.30} equals
$$
-\cal X^-_{s+1,i}\sum_{t=r+1}^{s-1} \psi_{\alpha_{r,t-1}}(\cal X^-_{r,t-1,i})\cal X^-_{t,j}\prod_{p=r}^{t-1}\cal H_{p,j}\prod_{p=t+1}^s
(\cal H_{p,i}-\delta_{i,j})+ U(\lie g)\lie n^+.
$$
The inductive step is now straightforward.
\end{pf}
\begin{cor}
For all~$1\le i\le j<\ell$, $1\le r\le i+1$, $r\le s\le j+1$, $e_\ell \bar{\cal U}_{r,s,i,j}\in U(\lie g)\lie n^+$.
\end{cor}

\subsection{}\label{C70}
We can now construct adapted families for all~$\beta\in\Psi$.
Suppose that~$\beta_{i,j}\in\Psi$, $i\le j\in I$. If~$\lambda\from\lambda+\beta_{i,j}\in (\Delta_\Psi)_1$, we have 
$\lambda(\cal H_{t,i-1})\ge \lambda(h_{i-1})>0$, $1\le t\le i-1$ and $\lambda(\cal H_{r,j-1})\ge \lambda(h_{j-1})>0$, $1\le r\le j-1$. 
Furthermore, if~$i=j$, $\lambda(h_{i-1})>1$, hence $\lambda(\cal H_{t,i-1}-1)>0$. Therefore,
$$
\cal H_{t,i-1}-\delta_{i,j},\cal H_{r,j-1}\in F_{\beta_{i,j}}(\lie h)^\times, \qquad 1\le t\le i-1,\,1\le r\le j-1.
$$

Clearly, $\{ \gamma\in R^+\,:\, \beta_{i,j}\le \gamma\}=\{ \beta_{r,s}\,:\, 1\le r\le i,\, r\le s\le j\}$. Define the elements
$\bou_{\beta_{i,j},\beta_{r,s}}\in U(\lie b)_{\beta_{i,j}-\beta_{r,s}}\tensor_{S(\lie h)} F_{\beta_{i,j}}(\lie h)$ by
\begin{equation}\label{C70.10}
\bou_{\beta_{i,j},\beta_{r,s}}=(-1)^{i+j+r+s}\, \frac{1+\delta_{r,s}}{1+\delta_{i,j}}\,\bar{\cal U}_{r,s,i-1,j-1}\tensor \prod_{t=r}^{i-1} (\cal H_{t,i-1}-\delta_{i,j})^{-1}
\prod_{t=r}^{j-1}\cal H_{t,j-1}^{-1}.
\end{equation}
\begin{lem}
Let~$\beta_{i,j}\in\Psi$. Then then $\{\bou_{\beta_{i,j},\beta_{r,s}}\,:\,1\le r\le i,\, r\le s\le j\}$ is an adapted
family for~$\beta_{i,j}$.
\end{lem}
\begin{pf}
We have
$$
\pi_{\lambda,\beta_{i,j}}(\bou_{\beta_{i,j},\beta_{r,s}})=C_{r,s}(i,j,\lambda)\cal U_{r,s,i-1,j-1}(\lambda),
$$
where
\begin{equation}
C_{r,s}(i,j,\lambda)=(-1)^{i+j+r+s}\, \frac{1+\delta_{r,s}}{1+\delta_{i,j}}\, \prod_{t=r}^{i-1} (\lambda(\cal H_{t,i-1})-\delta_{i,j})^{-1}
\prod_{t=r}^{j-1}(\lambda(\cal H_{t,j-1}))^{-1}\label{C70.10b}
\end{equation}
To shorten the notation, we denote $C_{r,s}(i,j,\lambda)$ (respectively, $\cal U_{r,s,i-1,j-1}(\lambda)$) by
$C_{r,s}$ (respectively, $\cal U_{r,s}$).
Observe that for all~$r\le s<j$ and for all~$k<i\le s$ we have
\begin{equation}\label{C70.20}
C_{r,s}=-\frac{1+\delta_{r,s}}{1+\delta_{r,s+1}}\, C_{r,s+1},\quad C_{k+1,s}=-\frac{1+\delta_{k+1,s}}{1+\delta_{k,s}}
(\lambda(\cal H_{k,i-1})-\delta_{i,j})\lambda(\cal H_{k,j-1}).
\end{equation}
Set
$$
u=\sum_{\gamma\in R^+\,:\,\beta_{i,j}\le\gamma} e_\gamma\tensor \pi_{\lambda,\beta_{i,j}}(\bou_{\beta_{i,j},\gamma})\in (\lie g\tensor V(\lambda+\beta_{i,j}))_{\lambda+\beta_{i,j}}.
$$
Since by~\eqref{C115.20}, $\bou_{\beta_{i,j},\beta_{i,j}}=1$, 
it remains to show
that~$e_k u=0$ for all~$1\le k
\le \ell$. For~$k\ge j$ this is immediate from~\eqref{C60.10a} and~\eqref{C50.10}. Suppose that~$i\le k<j$. Using~\eqref{C50.10}
and~\eqref{C60.10b} we obtain
$$
e_k u=\sum_{r=1}^i (C_{r,k+1}(1+\delta_{r,k})+C_{r,k}) e_{\beta_{r,k}}\tensor \cal U_{r,k+1} v_\lambda,
$$
which equals zero by~\eqref{C70.20}. Furthermore, if~$k< i$, it follows from~\eqref{C50.10}, \eqref{C60.10b} and~\eqref{C60.10c} that
\begin{multline*}
e_k u=
\sum_{s=k+1}^j \Big(C_{k+1,s}+C_{k,s}(1+\delta_{k+1,s})(\lambda(\cal H_{k,i-1})-\delta_{i,j})
\lambda(\cal H_{k,j-1})\Big)e_{\beta_{k,s}}\tensor \cal U_{k+1,s}v_\lambda\\
+\sum_{r=1}^k (C_{r,k+1}(1+\delta_{r,k})+C_{r,k}) e_{\beta_{r,k}}\tensor \cal U_{r,k+1}
v_\lambda.
\end{multline*}
Using~\eqref{C70.20} it is easy to see that~$e_k u=0$. 
\end{pf}

\subsection{}\label{C120}
Now we have all necessary ingredients to describe the relations. We set
$$
\cal Z_{\beta_{i,j},\Psi}=\prod_{t<i\,:\, \beta_{t,i}\in\Psi} 
(\cal H_{t,i-1}-\delta_{i,j})\prod_{t<j\,:\,\beta_{t,j}\in\Psi} \cal H_{t,j-1}\in F_{\beta_{i,j}}(\lie h)^\times
$$
and fix the isomorphism $\bt_{\Psi}^{\lie g}\to \bc\Delta_\Psi$ corresponding to the image of
$(\cal Z_{\beta,\Psi})_{\beta\in\Psi}\in \prod_{\beta\in\Psi} F_\beta(\lie h)^\times$ in~$G_\Psi$ (cf.~\ref{A95}).
\begin{prop}
Let~$\eta=\beta_{i,j}+\beta_{i,i}$, $i\not=j$. If~$t_{\lambda,\eta}=2$
then $\mathfrak R_\Psi(\lambda,\lambda+\eta)$ is spanned by the commutativity relation.
If~$t_{\lambda,\eta}=1$ then $\dim \mathfrak R_\Psi(\lambda,\lambda+\eta)=1$.
In particular, $\cal N_\eta=\emptyset$.
\end{prop}
\begin{pf}
Let~$i<j$. We have two different cases to consider.

$1^\circ.$ 
Let~$\eta=\beta_{i,i}+\beta_{i,j}\in\Psi+\Psi$, $i<j$. Suppose that $t_{\lambda,\eta}=2$, hence
$\lambda(h_{i-1})\ge 3$ and~$\lambda(h_{j-1})\ge 1$ by~\ref{C10}\cref{C10.c2}.
Then $\beta_{i,j}<\beta_{i,i}$ and it follows from~\propref{SP55} and~\lemref{C70} that
$$
\Pi_\lambda(\beta_{i,j},\beta_{i,i})=e_{i,i}\tensor e_{i,j},\quad
\Pi_\lambda(\beta_{i,i},\beta_{i,j})=e_{\beta_{i,j}}\tensor e_{\beta_{i,i}}+
e_{\beta_{i,i}}\tensor \bou_{\beta_{i,j},\beta_{i,i}}(\mu),
e_{\beta_{i,i}},\quad\mu=\lambda+\beta_{i,i}.
$$
Using~\eqref{C50.20} and an argument similar to that in the proof of~\propref{A95}, we conclude that
\begin{equation}\label{C120.0}
\cal U_{i,i,i-1,j-1}
=(-1)^{j-i-1}\prod_{t=i+1}^{j-1} \mu(\cal H_{t,j-1}) f_{j-1}\cdots f_i+\Ann_{ U(\lie g)} v_\mu\cap \Ann_{ U(\lie g)} e_{\beta_{i,i}},
\end{equation}
hence
\begin{equation}\label{C120.10}
\Pi_\lambda(\beta_{i,i},\beta_{i,j})=e_{\beta_{i,j}}\tensor e_{\beta_{i,i}}-2(\lambda(\cal H_{i,j-1})+2)^{-1} e_{\beta_{i,i}}\tensor e_{\beta_{i,j}}.
\end{equation}
To complete the computation of relations in this case, it remains to observe that
$$
\lambda(\cal Z_{\beta_{i,i},\Psi}\cal Z_{\beta_{i,j},\Psi}^{-1})(\lambda+\beta_{i,i})(\cal Z_{\beta_{i,j},\Psi})(\lambda+\beta_{i,j})(\cal Z_{\beta_{i,i},\Psi}^{-1})
=(\lambda(\cal H_{i,j-1})+2)\lambda(\cal H_{i,j-1})^{-1},
$$
and it is now easy to see that~$\mathfrak R_\Psi(\lambda,\lambda+\eta)$ is spanned by the commutativity relation.
Furthermore, if $t_{\lambda,\eta}=1$ then it follows from~\ref{C10}\cref{C10.c2} that~$\lambda(\cal H_{i,j-1})=0$ and the path
in question $(\lambda\from\lambda+\beta_{i,i}\from\lambda+\eta)$ is a relation by~\eqref{C120.10}.

$2^\circ.$ Let~$\eta=\beta_{i,j}+\beta_{j,j}$, $i<j$ and suppose that~$t_{\lambda,\eta}=2$. Then 
$\beta_{j,j}<\beta_{i,j}$ and by~\propref{SP55} and~\lemref{C70}
\begin{align*}
&\Pi_\lambda(\beta_{j,j},\beta_{i,j})=e_{\beta_{i,j}}\tensor e_{\beta_{j,j}},\\
&\Pi_\lambda(\beta_{i,j},\beta_{j,j})= e_{\beta_{j,j}}\tensor e_{\beta_{i,j}}+
e_{\beta_{i,j}}\tensor \bou_{\beta_{j,j},\beta_{i,j}}(\nu) e_{\beta_{i,j}},
\end{align*}
where~$\nu=\lambda+\beta_{i,j}$. Note that~$\cal X^-_{k,j-1}(\nu)e_{\beta_{i,j}}=0$,
$i<k\le j-1$.
It follows from~\lemref{C115} that
\begin{multline*}
\cal U_{i,j,j-1,j-1}
=\cal X^-_{i,j-1}(\nu) \prod_{k=i}^{j-1}(\nu(\cal H_{k,j})-\delta_{i,k}-1)
\\+\cal X^-_{i,j-1}(\nu-\varpi_{j-1})\prod_{k=i}^{j-1}\nu(\cal H_{k,j})+
\Ann_{ U(\lie g)}v_\nu\cap \Ann_{ U(\lie g)}e_{\beta_{i,j}}.
\end{multline*}
Furthermore, it is easy to see from~\eqref{C50.20} that for any~$\xi\in P$,
$$
\cal X^-_{i,j-1}(\xi)=
(-1)^{j-i-1} \prod_{k=i+1}^{j-1} \xi(\cal H_{k,j-1})f_{j-1}\cdots f_i\pmod{\Ann_{ U(\lie g)}e_{\beta_{i,j}}}.
$$
Since~$\nu(h_i)=\lambda(h_i)+1>0$, $f_{j-1}\cdots f_i\notin\Ann_{ U(\lie g)}v_\nu$ by~\corref{A45}. Therefore,
$$
\Pi_\lambda(\beta_{i,j},\beta_{j,j})=e_{\beta_{j,j}}\tensor e_{\beta_{i,j}}-
2 \lambda(\cal H_{i,j-1})^{-1} e_{\beta_{i,j}}\tensor e_{\beta_{j,j}}.
$$
Since
$$
\lambda(\cal Z_{\beta_{i,j},\Psi}\cal Z_{\beta_{j,j},\Psi}^{-1})(\lambda+\beta_{i,j})(\cal Z_{\beta_{j,j},\Psi})(\lambda+\beta_{j,j})
(\cal Z_{\beta_{i,j},\Psi}^{-1})=
(\lambda(\cal H_{i,j-1})(\lambda(\cal H_{i,j-1})-2)^{-1}
$$
it is now easy to see that we again obtain the commutativity relation.

Finally, if $t_{\lambda,\eta}=1$ then by~\ref{C10}\cref{C10.c2}, $\lambda(\cal H_{i,j-1})=2$ and so~$\Pi_\lambda(\beta_{i,j},\beta_{j,j})\in
\bigwedge^2 \lie n^+_\Psi$. Therefore, the corresponding path is a relation.
\end{pf}

\subsection{}\label{C140}
The next case $\eta=2\beta_{i,j}$, $i<j$ is more interesting since this is the only case in this paper
where $m_\eta>1$ for~$\eta\in 2\Psi$.
\begin{prop}
Suppose that~$\Psi\supset \Psi(i,j)$, $i<j\in I$ and
let~$\eta=2\beta_{i,j}=\beta_{i,i}+\beta_{j,j}\in\Psi+\Psi$. Assume that~$\lambda\in P^+$ is such that
that~${\Delta_\Psi}(\lambda,\lambda+\eta)\not=\emptyset$.
\begin{enumerate}[{\rm(i)}]
\item\label{C140.i} If~$t_{\lambda,\eta}=3$,
the unique relation is
\begin{equation}\label{C140.0}\begin{split}
&\lambda(\cal H_{i,j-1})^2(\lambda\from\lambda+\beta_{i,i}\from\lambda+\eta)-(\lambda(\cal H_{i,j-1})+2)^2 (\lambda\from\lambda+\beta_{j,j}\from\lambda+\eta)\\
&+(\lambda(\cal H_{i,j-1})+1)
(\lambda\from\lambda+\beta_{i,j}\from\lambda+\eta).
\end{split}
\end{equation}
\item\label{C140.ii} If~$t_{\lambda,\eta}=2$, the unique relation is 
$$
(\lambda\from\lambda+\beta_{i,i}\from\lambda+\eta)+2 (\lambda\from\lambda+\beta_{i,j}\from\lambda+\eta).
$$
\item\label{C140.iii} If~$t_{\lambda,\eta}=1$ then $\lambda(\cal H_{i,j-1})=0$ and $\mathfrak R_\Psi(\lambda,\lambda+\eta)=0$.
\end{enumerate}
Thus, $\mathfrak R_\Psi(\lambda,\lambda+\eta)$ has dimension $\lfloor t_{\lambda,\eta}/2\rfloor$ and is
generic if and only if~$t_{\lambda,\eta}>1$.
\end{prop}
\begin{pf}
Suppose first that~$t_{\lambda,\eta}=3$.
We have $\beta_{j,j}<\beta_{i,j}<\beta_{i,i}$, hence
\begin{align*}
&\Pi_\lambda(\beta_{j,j},\beta_{i,i})=e_{\beta_{i,i}}\tensor e_{\beta_{j,j}},
\\
&\Pi_\lambda(\beta_{i,i},\beta_{j,j})=e_{\beta_{j,j}}\tensor e_{\beta_{i,i}}+
e_{\beta_{i,i}}\tensor 
\bou_{\beta_{j,j},\beta_{i,i}}(\mu) e_{\beta_{i,i}}+e_{\beta_{i,j}}\tensor\bou_{\beta_{j,j},\beta_{i,j}}(\mu)e_{\beta_{i,i}},\\
&\Pi_\lambda(\beta_{i,j},\beta_{i,j})=e_{\beta_{i,j}}\tensor e_{\beta_{i,j}}+e_{\beta_{i,i}}\tensor\bou_{\beta_{i,j},\beta_{i,i}}(\nu) e_{\beta_{i,j}},
\end{align*}
where~$\mu=\lambda+\beta_{i,i}$, $\nu=\lambda+\beta_{i,j}$.
Using~\lemref{C115} and~\eqref{C120.0} we can write
$$
\bou_{\beta_{j,j},\beta_{i,i}}(\mu)
=((\mu(\cal H_{k,j-1})-1)\mu(\cal H_{k,j-1}))^{-1}
(f_{j-1}\cdots f_i)^2+\Ann_{U(\lie g)} e_{\beta_{i,i}}\cap \Ann_{U(\lie g)}v_\mu.
$$
We claim that~$(f_{j-1}\cdots f_i)^2\notin \Ann_{ U(\lie n^-)} v_\mu$. Indeed, since
$\mu(h_i)=\lambda(h_i)+2\ge 2$, it follows from~\corref{A45} that $F=f_{j-2}\cdots f_i f_{j-1}\cdots f_i
\notin\Ann_{ U(\lie n^-)} v_\mu$.
Now, if~$\zeta$ is the weight of $F v_\mu$, we have $\zeta(h_{j-1})=\lambda(h_{j-1})\ge 2$. Furthermore, $e_{j-1}^2 Fv_\mu=0$. It now follows from the elementary $\lie{sl}_2$ theory that~$f_{j-1} F v_\mu\not=0$.
Thus,
$$\bou_{\beta_{j,j},\beta_{i,i}}(\mu) e_{\beta_{i,i}}
=2 ((\lambda(\cal H_{i,j-1})+1)(\lambda(\cal H_{i,j-1})+2))^{-1} e_{\beta_{j,j}}.
$$
A computation similar to that of~$\Pi_{\lambda}(\beta_{i,j},\beta_{j,j})$ in~\ref{C120} yields
$$
\bou_{\beta_{j,j},\beta_{i,j}}(\mu) e_{\beta_{i,i}}=-(\lambda(\cal H_{i,j-1})+2)^{-1}e_{\beta_{i,j}}.
$$
Thus,
\begin{multline}\label{C140.10}
\Pi_\lambda(\beta_{i,i},\beta_{j,j})=
e_{\beta_{j,j}}\tensor e_{\beta_{i,i}}-(\lambda(\cal H_{i,j-1})+2)^{-1} e_{\beta_{i,j}}
\tensor e_{\beta_{i,j}}\\+2((\lambda(\cal H_{i,j-1})+1)(\lambda(\cal H_{i,j-1})+2))^{-1}e_{\beta_{i,i}}\tensor e_{\beta_{j,j}}.
\end{multline}
The computation of~$\Pi_\lambda(\beta_{i,j},\beta_{i,j})$ is similar to that of
$\Pi_\lambda(\beta_{i,i},\beta_{i,j})$ in~\ref{C120} and yields
$$
\Pi_\lambda(\beta_{i,j},\beta_{i,j})=e_{\beta_{i,j}}\tensor e_{\beta_{i,j}}-
4 (\lambda(\cal H_{i,j-1}))^{-1} e_{\beta_{i,i}}\tensor e_{\beta_{j,j}}.
$$
Note that none of these paths is a relation. 
To complete the computation of relations in this case, it remains to note that
\begin{align*}
&\lambda(\cal Z_{\beta_{i,i},\Psi})(\lambda+\beta_{i,i})(\cal Z_{\beta_{j,j},\Psi})=(\lambda(\cal H_{i,j-1})+2)(\lambda(\cal H_{i,j-1})+1) z,\\
&\lambda(\cal Z_{\beta_{i,j},\Psi})(\lambda+\beta_{i,j})(\cal Z_{\beta_{i,j},\Psi})=(\lambda(\cal H_{i,j-1}))^2 z,\\
&\lambda(\cal Z_{\beta_{j,j},\Psi})(\lambda+\beta_{j,j})(\cal Z_{\beta_{i,i},\Psi})=\lambda(\cal H_{i,j-1})(\lambda(\cal H_{i,j-1})-1)z,
\end{align*}
where~$z\in\bc^\times$. The relation~\eqref{C140.0} is now straightforward. Since all coefficients in it are positive integers, 
$\mathfrak R_\Psi(\lambda,\lambda+\eta)$ is generic.

If~$t_{\lambda,\eta}=2$ then by~\ref{C10}\cref{C10.c3},
$i=j-1$ and $\lambda(\cal H_{i,j-1})=1$ and we immediately obtain the relation using the above formulae.
Finally, if~$t_{\lambda,\eta}=1$ then $\lambda(\cal H_{i,j-1})=0$ and it is easy to see from~\eqref{C140.10} that the
corresponding path is not a relation.
\end{pf}

\subsection{}\label{C150}
We now present an infinite dimensional example which in particular includes the remaining rank~$2$ case. 
Let~$\Psi=\Psi(1,2)$, $\ell\ge 2$. Since~$\beta_{1,1}=2\varpi_1=\theta$, $\beta_{1,2}=\varpi_2$ and~$\beta_{2,2}=-2\varpi_1+2\varpi_2$ it is clear 
that $\lambda,\mu\in P^+$ are in the same connected component of~$\Delta_\Psi$ only if~$\lambda(h_i)=\mu(h_i)$, $2<i\le \ell$. Therefore, it 
is enough to describe the connected components of~$\Delta_\Psi$ in the case~$\ell=2$. Identify~$P$ with~$\bz\times \bz$ and write $(\lambda(h_1),\lambda(h_2))$
for~$\lambda\in P$. Since~$\varphi(\theta)=(2,0)$, $\varphi(\beta_{1,2})=(1,1)$ and~$\varphi(\beta_{2,2})=(0,2)$, we conclude that
the only sinks in~$\Delta_\Psi$ are $(0,0)$, $(0,1)$ and~$(1,0)$. Furthermore, if~$(m,n)$ and~$(m',n')$ are in the same
connected component it is immediate that~$m=m'\pmod 2$. Since we have $(0,0)\from (2,0)\from (2,1)\to (0,1)$, we conclude that
there are two connected components, $\Delta_\Psi[(r,0)]$, $r=0,1$ (if~$\ell>2$, each of these components has infinite multiplicity).
We have $\Delta_\Psi[(r,0)]_0=\{ (m,n)\,:\, m,n\in\bz_+, m=r\pmod 2\}$ and the arrows are $(m,n)\from (m+2,n)$, $m,n\in\bz_+$,
$(m,n)\from (m,n+1)$, $m>0$, $n\in\bz_+$ and $(m,n)\from (m-2,n+2)$, $m\ge 2$, $n\in\bz_+$. Thus, the quivers~$\Delta_\Psi[(0,r)]$, $r=0,1$
are, respectively,
$$
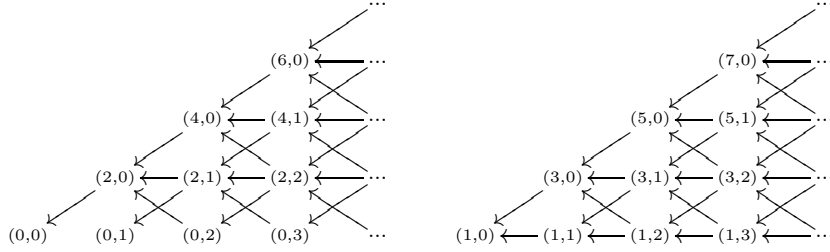
\def\dgeverynode{\scriptscriptstyle}
\divide\dgARROWLENGTH by 2\dgVERTPAD=4pt\dgHORIZPAD=4pt
\begin{diagram}
\node{}\node{}\node{}\node{}\node{\scriptstyle\dots}\arrow{sw}
\\
\node{}\node{}\node{}\node{(6,0)}\arrow{sw}\node{\scriptstyle\dots}\arrow{w}\arrow{sw}
\\
\node{}\node{}\node{(4,0)}\arrow{sw}\node{(4,1)}\arrow{w}\arrow{sw}\node{\scriptstyle\dots}\arrow{nw}\arrow{sw}\arrow{w}
\\
\node{}\node{(2,0)}\arrow{sw}\node{(2,1)}\arrow{w}\arrow{sw}\node{(2,2)}\arrow{nw}\arrow{w}\arrow{sw}
\node{\scriptstyle\dots}\arrow{nw}\arrow{sw}\arrow{w}\\
\node{(0,0)}\node{(0,1)}\node{(0,2)}\arrow{nw}\node{(0,3)}\arrow{nw}\node{\scriptstyle\dots}\arrow{nw}\\
\end{diagram}
\qquad
\begin{diagram}
\node{}\node{}\node{}\node{}\node{\scriptstyle\dots}\arrow{sw}
\\
\node{}\node{}\node{}\node{(7,0)}\arrow{sw}\node{\scriptstyle\dots}\arrow{w}\arrow{sw}
\\
\node{}\node{}\node{(5,0)}\arrow{sw}\node{(5,1)}\arrow{w}\arrow{sw}\node{\scriptstyle\dots}\arrow{nw}\arrow{sw}\arrow{w}
\\
\node{}\node{(3,0)}\arrow{sw}\node{(3,1)}\arrow{w}\arrow{sw}\node{(3,2)}\arrow{nw}\arrow{w}\arrow{sw}
\node{\scriptstyle\dots}\arrow{nw}\arrow{sw}\arrow{w}\\
\node{(1,0)}\node{(1,1)}\arrow{w}\node{(1,2)}\arrow{nw}\arrow{w}\node{(1,3)}\arrow{nw}\arrow{w}\node{\scriptstyle\dots}\arrow{nw}\arrow{w}\\
\end{diagram}
$$
Both are translation quivers with~$\tau((m,n))=(m,n-2)$, $m>0$, $n\ge 2$. The relations are: the commutativity relations in
$$
\def\dgeverynode{\scriptscriptstyle}
\divide\dgARROWLENGTH by 2
\dgVERTPAD=4pt\dgHORIZPAD=4pt
\begin{diagram}
\node{}\node{(m+2,n)}\arrow{sw}\node{(m+2,n+1)}\arrow{w}\arrow{sw}\\
\node{(m,n)}\node{(m,n+1)}\arrow{w} 
\end{diagram}\qquad
\begin{diagram}
\node{(m+2,n)}\node{(m+2,n+1)}\arrow{w}\\\node{}\node{(m,n+2)}\arrow{nw}\node{(m,n+3)}\arrow{w}\arrow{nw}
\end{diagram}
$$
for all~$m>0$, $n\in\bz_+$, the zero relations $(2,n)\from (2,n+1)\from (0,n+3)$, $n\ge 0$
and
\begin{multline*}
m^2((m,n)\from (m+2,n)\from (m,n+2))-(m+2)^2((m,n)\from (m-2,n+2)\from (m,n+2))\\+(m+1)((m,n)\from(m,n+1)\from(m,n+2)),\qquad m>1,
\end{multline*}
and, finally, $((1,n)\from (3,n)\from (1,n+2))+2((1,n)\from(1,n+1)\from(1,n+2))$.
Thus, if~$\ell=2$ and~$|\Psi|>1$, the algebra $\bs_\Psi^{\lie g}$ is the direct sum of two non-isomorphic connected Koszul subalgebras of left global dimension~$3$.

\subsection{}\label{C160}
Next, we consider~$\eta\in\Psi+\Psi$ with~$m_\eta=4$.
We have three possibilities here, and the computations turn out to be rather different.
\begin{prop}
Let~$i<j<k\in I$ and let~$x_\lambda=\lambda(\cal H_{i,j-1})$, $y_\lambda=\lambda(\cal H_{j,k-1})$.
Assume that $\Psi\supset\Psi(i,j,k)$ and let~$\eta\in \Psi(i,j,k)+\Psi(i,j,k)$, $m_\eta=4$.
\begin{enumerate}[{\rm(i)}]
 \item\label{C160.i} Let~$\eta=\beta_{i,i}+\beta_{j,k}=\beta_{i,j}+\beta_{i,k}$. If~$t_{\lambda,\eta}=4$ then~$\mathfrak R_\Psi(\lambda,\lambda+\eta)$ is spanned by
\begin{subequations}
\begin{multline}\label{C160.i.4a}
2 (1 + y_\lambda) (\lambda\from\lambda+\beta_{i, i}\from\lambda+\eta) - (3 + x_\lambda) (x_\lambda+y_\lambda+3) (\lambda\from\lambda+\beta_{i, j}\from\lambda+\eta)\\+ 
(2 + x_\lambda) (x_\lambda+y_\lambda+4) (\lambda\from\lambda+\beta_{i,k}\from\lambda+\eta)
\end{multline}
and
\begin{multline}\label{C160.i.4b}
x_\lambda (x_\lambda+y_\lambda+2) (\lambda\from\lambda+\beta_{i, i}\from\lambda+\eta)+ (2x_\lambda+y_\lambda+4) (\lambda\from\lambda+\beta_{i,j}\from\lambda+\eta)\\
- (2 + x_\lambda) (x_\lambda+y_\lambda+4) (\lambda\from\lambda+\beta_{j, k}\from\lambda+\eta)
\end{multline}
\end{subequations}
If~$t_{\lambda,\eta}=2$ then~$\mathfrak R_\Psi(\lambda,\lambda+\eta)$ is spanned by
\begin{equation}\label{C160.i.2}
(x_\lambda+y_\lambda+1)(\lambda\from\lambda+\beta_{i, i}\from\lambda+\eta) + (x_\lambda+y_\lambda+4)(\lambda\from\lambda+\beta_{i,k}\from\lambda+\eta).
\end{equation}
Finally, if~$t_{\lambda,\eta}=1$ then the unique path is a relation. 

\item\label{C160.ii} Let~$\eta=\beta_{j,j}+\beta_{i,k}=\beta_{i,j}+\beta_{j,k}$. 
If~$t_{\lambda,\eta}=4$, $\mathfrak R_\Psi(\lambda,\lambda+\eta)$ is spanned by
\begin{subequations}
\begin{multline}\label{C160.ii.4a}
(x_\lambda-1) (2 + y_\lambda) (\lambda\from\lambda+\beta_{i, k}\from\lambda+\eta) - (1 + x_\lambda) y_\lambda 
(\lambda\from\lambda+\beta_{j, j}\from\lambda+\eta)\\- (x_\lambda - y_\lambda-1) (\lambda\from\lambda+\beta_{j, k}\from\lambda+\eta)
\end{multline}
and
\begin{multline}\label{C160.ii.4b}
(1 + x_\lambda) y_\lambda (\lambda\from\lambda+\beta_{i, j}\from\lambda+\eta) + 
 2 (2 + x_\lambda + y_\lambda) (\lambda\from\lambda+\beta_{i, k}\from\lambda+\eta) \\- (2 + x_\lambda) (1 + y_\lambda) (\lambda\from\lambda+\beta_{j, k}\from\lambda+\eta)
\end{multline}
\end{subequations}
and~$\mathfrak R_\Psi(\lambda,\lambda+\eta)$ is generic unless $x_\lambda=y_\lambda+1$.
If~$t_{\lambda,\eta}=2$ then $\mathfrak R_\Psi(\lambda,\lambda+\eta)$ is spanned by
\begin{equation}\label{C160.ii.2}
(x_\lambda-1)(\lambda\from\lambda+\beta_{i,j}\from\lambda+\eta) + (x_\lambda+2)(\lambda\from\lambda+\beta_{j,j}\from\lambda+\eta),
\end{equation}
If~$t_{\lambda,\eta}=1$  then~$\lambda(\cal H_{i,j-1})=1$ and~$\mathfrak R_\Psi(\lambda,\lambda+\eta)=0$.

\item\label{C160.iii} Let~$\eta=\beta_{i,j}+\beta_{k,k}=\beta_{i,k}+\beta_{j,k}$. If~$t_{\lambda,\eta}=4$, $\mathfrak R_\Psi(\lambda,\lambda+\eta)$ is spanned by
\begin{subequations}
\begin{multline}\label{C160.iii.4a}
(1 + y_\lambda) (3+x_\lambda+y_\lambda) (\lambda\from\lambda+\beta_{j, k}\from\lambda+\eta) - (2 + y_\lambda) (2+x_\lambda+y_\lambda) (\lambda\from\lambda+\beta_{i, k}\from\lambda+\eta) \\
+ 2 (1 + x_\lambda) (\lambda\from\lambda+\beta_{i, j}\from\lambda+\eta)
\end{multline}
and
\begin{multline}\label{C160.iii.4b}
2 (1 + x_\lambda) (\lambda\from\lambda+\beta_{k, k}\from\lambda+\eta) + (y_\lambda-1)(x_\lambda+y_\lambda+1)(\lambda\from\lambda+\beta_{j, k}\from\lambda+\eta)
\\- y_\lambda (x_\lambda+y_\lambda) (\lambda\from\lambda+\beta_{i, k}\from\lambda+\eta).
\end{multline}
\end{subequations}
If~$t_{\lambda,\eta}=2$, then the unique relation is
\begin{equation}\label{C160.iii.2}
(x_\lambda+4) (\lambda\from\lambda+\beta_{i,j}\from\lambda+\eta)+(x_\lambda+1)(\lambda\from\lambda+\beta_{j,k}\from\lambda+\eta).
\end{equation}
If~$t_{\lambda,\eta}=1$ then~$\dim\mathfrak R_\Psi(\lambda,\lambda+\eta)=1$.
\end{enumerate}
In particular, in cases~\eqref{C160.i} and~\eqref{C160.iii}, $\cal N_\eta=\emptyset$ while in the case~\eqref{C160.ii}, 
$\{\lambda\in \cal N_\eta\,:\, t_{\lambda,\eta}>1\}=P^+\cap \{\xi\in\lie h^*\,:\, \xi(\cal H_{i,j-1}-\cal H_{j,k-1}-1)=0\}$.
\end{prop}
\begin{pf}
Suppose first that~$\eta=\beta_{i,j}+\beta_{i,k}$. We have~$\beta_{j,k}<\beta_{i,k}<\beta_{i,j}<\beta_{i,i}$.
Suppose first that
$t_{\lambda,\eta}=4$. 
The first two paths are straightforward
\begin{align}
&\Pi_\lambda(\beta_{j,k},\beta_{i,i})=e_{\beta_{i,i}}\tensor e_{\beta_{j,k}},
\label{C160.0a}\\
&\Pi_\lambda(\beta_{i,k},\beta_{i,j})=
e_{\beta_{i,j}}\tensor e_{\beta_{i,k}}-
2(x_\lambda+1)^{-1} e_{\beta_{i,i}}\tensor e_{\beta_{j,k}}.\label{C160.0b}
\end{align}
The next path is more involved. We have
$$
\Pi_\lambda(\beta_{i,j},\beta_{i,k})=
e_{\beta_{i,k}}\tensor e_{\beta_{i,j}}+
e_{\beta_{i,i}}\tensor 
\bou_{\beta_{i,k},\beta_{i,i}}(\mu) e_{\beta_{i,j}}+
e_{\beta_{i,j}}\tensor 
\bou_{\beta_{i,k},\beta_{i,j}}(\mu) e_{\beta_{i,j}}
$$
where~$\mu=\lambda+\beta_{i,j}$.
We claim that
\begin{multline}\label{C160.10}
\cal U_{i,i,i-1,k-1}(\mu)=\cal X^-_{i,k-1}(\mu)=
(-1)^{k-i-1} \Big(\prod_{t=i+1}^{k-1} \mu(\cal H_{t,k-1}) f_{k-1}\cdots f_i\\-
\prod_{t=i+1}^{k-1} (\mu(\cal H_{t,k-1})+\delta_{t,j}) f_{j-1}\cdots f_i f_{k-1}\cdots f_j\Big)
+\Ann_{ U(\lie g)}e_{\beta_{i,j}}\cap \Ann_{ U(\lie g)} v_\mu.
\end{multline}
Indeed, suppose that~$\bof_\sigma\notin \Ann_{ U(\lie g)} e_{\beta_{i,j}}$, $\sigma\in\Sigma(i,k-1)$. Then by~\eqref{C50.20} we must have
$\sigma(i)=k-i$ or~$\sigma(j)=k-i$. If~$\sigma(i)=k-i$ then~$\bof_\sigma=f_{k-1}\cdots f_i$ by the definition of~$\Sigma(i,k-1)$. Otherwise, 
we must have $\sigma(r)=k-r$, $r\le j\le k-1$ and so~$\bof_\sigma=\bof_{\sigma'} f_{k-1}\cdots f_j$, where
$\sigma'\in\Sigma(i,j-1)$. Then $\bof_\sigma e_{\beta_{i,j}}=\bof_{\sigma'} e_{\beta_{i,k}}$ hence
$\bof_{\sigma'}=f_{j-1}\cdots f_i$.
Since~$\mu(h_{r})=\lambda(h_r)+1>0$, $r=i,j$, it follows from~\corref{A45} that the vectors 
$
f_{k-1}\cdots f_i v_\mu$, $f_{j-1}\cdots f_i f_{k-1}\cdots f_j v_\mu
$
are non-zero and linearly independent. Therefore,
\begin{multline*}
\bou_{\beta_{i,k},{i,i}}(\mu)
e_{\beta_{i,j}}=-2\prod_{t=i}^{k-1}\mu(\cal H_{t,k-1})^{-1} \Big(2\prod_{t=i+1}^{k-1} \mu(\cal H_{t,k-1})-\prod_{t=i+1}^{k-1} (\mu(\cal H_{t,k-1})+\delta_{t,j})\Big)
e_{\beta_{j,k}}\\=
-2 (\mu(\cal H_{j,k-1})-1) (\mu(\cal H_{i,k-1})\mu(\cal H_{j,k-1}))^{-1}
e_{\beta_{j,k}}.
\end{multline*}
An already familiar computation yields~$\bou_{\beta_{i,k},\beta_{i,j}}(\mu) e_{\beta_{i,j}}=-\mu(\cal H_{j,k-1})^{-1} e_{\beta_{i,k}}$ and we obtain
\begin{equation}
\Pi_\lambda(\beta_{i,j},\beta_{i,k})=
 e_{\beta_{i,k}}
\tensor e_{\beta_{i,j}}-
(y_\lambda+1)^{-1}\, e_{\beta_{i,j}}\tensor e_{\beta_{i,k}}-
\frac{2 y_\lambda}{(y_\lambda+1)(x_\lambda+y_\lambda+2)}\, e_{\beta_{i,i}}\tensor e_{\beta_{j,k}}.\label{C160.0c}
\end{equation}
Finally, 
\begin{equation}\label{C160.0d}
\begin{split}
\Pi_{\lambda}(\beta_{i,i},\beta_{j,k})
&=
  e_{\beta_{j,k}}\tensor e_{\beta_{i,i}}-
(x_\lambda+2)^{-1} e_{\beta_{i,k}}\tensor e_{\beta_{i,j}}+
2 ((x_\lambda+2)(x_\lambda+y_\lambda+3))^{-1} e_{\beta_{i,i}}\tensor e_{\beta_{j,k}}\\&\qquad
-\frac{x_\lambda+1}{(x_\lambda+2)(x_\lambda+y_\lambda+3)}\, e_{\beta_{i,j}}\tensor e_{\beta_{i,k}}.
\end{split}
\end{equation}
It is now straightforward to show that~$\mathfrak R_\Psi(\lambda,\lambda+\eta)$ is spanned by the elements~\eqref{C160.i.4a} and~\eqref{C160.i.4b} and that~$\mathfrak R_\Psi(\lambda,\lambda+\eta)$ is generic.

Suppose now that~$0<t_{\lambda,\eta}<4$. By~\ref{C10}\cref{C10.c4} we have two possibilities. If~$\lambda(h_i)=0$ and~$i=j-1$
(hence~$x_\lambda=0$) we obtain from~\eqref{C160.0b} and~\eqref{C160.0d} that $\mathfrak R_\Psi(\lambda,\lambda+\eta)$ is spanned by the element~\eqref{C160.i.2}.
Finally, if $\lambda(h_j)=0$ and~$j=k-1$ (hence~$y_\lambda=0$), it follows from~\eqref{C160.0c} that $\Pi_\lambda(\beta_{i,j},\beta_{i,k})\in\bigwedge^2\lie n^+_\Psi$
hence the unique path~$(\lambda\from\lambda+\beta_{i,j}\from\lambda+\eta)$ is a relation.

We now prove~\eqref{C160.ii}.
Let $\eta=\beta_{i,j}+\beta_{j,k}=\beta_{i,k}+\beta_{j,j}$. The first three paths are rather easy and we obtain
\begin{align*}
&\Pi_\lambda(\beta_{j,k},\beta_{i,j})=e_{\beta_{i,j}}\tensor
e_{\beta_{j,k}},\\
&\Pi_\lambda(\beta_{i,k},\beta_{j,j})=
e_{\beta_{j,j}}\tensor e_{\beta_{i,k}}-(x_\lambda+1)^{-1}
e_{\beta_{i,j}}\tensor e_{\beta_{j,k}},\\
&\Pi_\lambda(\beta_{j,j},\beta_{i,k})=
e_{\beta_{i,k}}\tensor e_{\beta_{j,j}}-(y_\lambda+2)^{-1}
e_{\beta_{i,j}}\tensor e_{\beta_{j,k}}.
\end{align*}
The last path has some new features. We have
\begin{multline*}
\Pi_\lambda(\beta_{i,j},\beta_{j,k})=
e_{\beta_{j,k}}\tensor e_{\beta_{i,j}}+
e_{\beta_{j,j}}\tensor \bou_{\beta_{j,k},\beta_{j,j}}(\nu)
e_{\beta_{i,j}}\\+
e_{\beta_{i,k}}\tensor \bou_{\beta_{j,k},\beta_{i,k}}(\nu)
e_{\beta_{i,j}}+
e_{\beta_{i,j}}\tensor \bou_{\beta_{j,k},\beta_{i,j}}(\nu)
e_{\beta_{i,j}},
\end{multline*}
where~$\nu=\lambda+\beta_{i,j}$. Two terms are already familiar
\begin{align*}
&\bou_{\beta_{j,k},\beta_{j,j}}(\nu) e_{\beta_{i,j}}=-2 (\nu(\cal H_{j,k-1}))^{-1} e_{\beta_{i,k}},\\
&\bou_{\beta_{j,k},\beta_{i,k}}(\nu) e_{\beta_{i,j}}=-2(\nu(\cal H_{i,j-1}))^{-1}
e_{\beta_{j,j}}.
\end{align*}
Furthermore, we have, modulo~$\Ann_{ U(\lie g)}e_{\beta_{i,j}}$,
$$
\cal U_{i,j,j-1,k-1}(\nu)=\cal X^-_{i,k-1}(\nu)\prod_{t=i}^{j-1} (\nu(\cal H_{t,j-1})-\delta_{t,i})+
\cal X^-_{i,j-1}(\nu)\cal X^-_{j,k-1}(\nu)\prod_{t=i}^{j-1}\nu(\cal H_{t,k-1}).
$$
Using~\eqref{C160.10}, we obtain
\begin{multline*}
\cal U_{i,j,j-1,k-1}(\nu)=(-1)^{k-i}\prod_{t=i}^{j-1} (\nu(\cal H_{t,j-1})-\delta_{t,i})\Big(-\prod_{t=i+1}^{k-1} \nu(\cal H_{t,k-1}) f_{k-1}\cdots f_i\\+
\prod_{t=i+1}^{k-1} (\nu(\cal H_{t,k-1})+\delta_{t,j}) f_{j-1}\cdots f_i f_{k-1}\cdots f_j\Big)
\\+
(-1)^{k-i}\prod_{t=i+1}^{j-1}\nu(\cal H_{t,j-1})\prod_{t=j+1}^{k-1} \nu(\cal H_{t,k-1})
\prod_{t=i}^{j-1}\nu(\cal H_{t,k-1}) f_{j-1}\cdots f_i f_{k-1}\cdots f_j+\Ann_{ U(\lie g)}e_{\beta_{i,j}}.
\end{multline*}
Thus,
\begin{multline*}
\bou_{\beta_{j,k},\beta_{i,j}}(\nu)e_{\beta_{i,j}}\\
=\Big(\nu(\cal H_{i,j-1})\nu(\cal H_{i,k-1})\nu(\cal H_{j,k-1})\Big)^{-1}
\Big(\nu(\cal H_{i,k-1})-(\nu(\cal H_{i,j-1})-1)(\nu(\cal H_{j,k-1})-1)\Big) e_{\beta_{j,k}},
\end{multline*}
hence
\begin{multline*}
\Pi_{\lambda}(\beta_{i,j},\beta_{j,k})=
e_{\beta_{j,k}}\tensor
e_{\beta_{i,j}}
-2(y_\lambda+1)^{-1} e_{\beta_{j,j}}\tensor e_{\beta_{i,k}}-2x_\lambda^{-1} e_{\beta_{i,k}}\tensor e_{\beta_{j,j}}\\
+
\frac{(x_\lambda+y_\lambda+2-(x_\lambda-1)y_\lambda)}{x_\lambda(y_\lambda+1)(x_\lambda+y_\lambda+2)}\, e_{\beta_{i,j}}\tensor e_{\beta_{j,k}},
\end{multline*}
and we obtain the relations \eqref{C160.ii.4a} and~\eqref{C160.ii.4b}.
It is now straightforward to check that $\mathfrak R_\Psi(\lambda,\lambda+\eta)$
is generic if and only if~$x_\lambda\not=y_\lambda+1$.

Suppose now that~$0<t_{\lambda,\eta}<4$. Using~\ref{C10}\cref{C10.c4} and the above, we conclude that 
if~$t_{\lambda,\eta}=2$, the unique relation is~\eqref{C160.ii.2},
while in the case~$t_{\lambda,\eta}=1$, $x_\lambda=1$ and the unique path is not a relation.

In part~\eqref{C160.iii} we encounter some new features. We have
$\beta_{k,k}<\beta_{j,k}<\beta_{i,k}<\beta_{i,j}$
and 
\begin{align*}
&\Pi_{\lambda}(\beta_{k,k},\beta_{i,j})=e_{\beta_{i,j}}\tensor e_{\beta_{k,k}},\\
&\Pi_{\lambda}(\beta_{j,k},\beta_{i,k})=e_{\beta_{i,k}}\tensor e_{\beta_{j,k}}-
2y_\lambda^{-1} e_{\beta_{i,j}}\tensor e_{\beta_{k,k}},
\\
&\Pi_{\lambda}(\beta_{i,k},\beta_{j,k})=e_{\beta_{j,k}}\tensor
e_{\beta_{i,k}}-(x_\lambda+1)^{-1} e_{\beta_{i,k}}\tensor e_{\beta_{j,k}}
-\frac{2 x_\lambda}{(x_\lambda+1)(x_\lambda+y_\lambda+1)}\,
e_{\beta_{i,j}}\tensor e_{\beta_{k,k}}.
\end{align*}
The remaining path is rather interesting. Indeed, this turns out to be one of the only two cases when
the last term in~\eqref{C115.10} does not lie in the annihilator of the corresponding root
vector. As before, 
\begin{multline*}
\Pi_\lambda(\beta_{i,j},\beta_{k,k})=e_{\beta_{k,k}}\tensor e_{\beta_{i,j}}+
e_{\beta_{j,k}}\tensor \bou_{\beta_{k,k},\beta_{j,k}}(\mu) e_{\beta_{i,j}}\\
+e_{\beta_{i,k}}\tensor \bou_{\beta_{k,k},\beta_{i,k}}(\nu)
e_{\beta_{i,j}}+e_{\beta_{i,j}}\tensor \bou_{\beta_{k,k},\beta_{i,j}}(\mu) e_{\beta_{i,j}},
\end{multline*}
where~$\mu=\lambda+\beta_{i,j}$. Note that $\mu(\cal H_{r,k-1})=\lambda(\cal H_{r,k-1})+1>\lambda(h_{k-1})>2$, $r=i,j$.
We immediately get
$$
\bou_{\beta_{k,k},\beta_{j,k}}(\mu) e_{\beta_{i,j}}=-\mu(\cal H_{j,k-1})^{-1}e_{\beta_{i,k}}.
$$
Furthermore, by~\lemref{C115} we have
\begin{multline*}
\cal U_{i,k,k-1,k-1}(\mu)=
\cal X^-_{i,k-1}(\mu) \prod_{t=i}^{k-1} (\mu(\cal H_{t,k-1})-\delta_{t,i}-1)+\cal X^-_{i,k-1}(\mu-\varpi_{k-1})\prod_{t=i}^{k-1} \mu(\cal H_{t,k-1})
\\-\cal X^-_{i,j,k-1}(\mu-\varpi_{k-1})\cal X^-_{j,k-1}(\mu)\prod_{t=i}^{j-1}\mu(\cal H_{t,k-1})\prod_{t=j+1}^{k-1} (\mu(\cal H_{t,k-1})-1)
+\Ann_{ U(\lie g)}e_{\beta_{i,j}}.
\end{multline*}
Using~\eqref{C160.10} we obtain
\begin{multline*}
\bou_{\beta_{k,k},\beta_{i,k}}(\mu)
e_{\beta_{i,j}}=
-\frac12 (\mu(\cal H_{i,k-1})-1)^{-1}\Bigg(\frac{(\mu(\cal H_{i,k-1})-2)(\mu(\cal H_{j,k-1})-1)}{\mu(\cal H_{j,k-1})\mu(\cal H_{i,k-1})}
\\+\frac{\mu(\cal H_{j,k-1})-2}{\mu(\cal H_{j,k-1})-1}+\frac{1}{\mu(\cal H_{j,k-1})(\mu(\cal H_{j,k-1})-1)}\Bigg) e_{\beta_{j,k}}
=-\frac{\mu(\cal H_{j,k-1})-1}{\mu(\cal H_{i,k-1}) \mu(\cal H_{j,k-1})}\, e_{\beta_{j,k}}.
\end{multline*}
To compute the remaining term observe that by~\lemref{C115}
\begin{multline*}
(-1)^{i+j}\cal U_{i,j,k-1,k-1}(\mu)=\cal X^-_{j,k-1}(\mu-\varpi_{k-1})\cal X^-_{i,k-1} \prod_{t=i}^{j-1} \mu(\cal H_{t,k-1}-1-\delta_{t,i})\\+
\cal X^-_{i,k-1}(\mu-\varpi_{k-1})\cal X^-_{j,k-1}(\mu)\prod_{t=i}^{j-1} \mu(\cal H_{t,k-1})+\Ann_{ U(\lie g)}e_{\beta_{i,j}}\\
=\prod_{i\le t\le k-1,\,t\not=j} \mu(\cal H_{t,k-1}-1-\delta_{t,i}) f_{k-1}\cdots f_j \Big(\prod_{t=i+1}^{k-1} \mu(\cal H_{t,k-1}) f_{k-1}\cdots f_i\\-
\prod_{t=i+1}^{k-1}(\mu(\cal H_{t,k-1})+\delta_{t,j})f_{j-1}\cdots f_i f_{k-1}\cdots f_j\Big)\\+
\prod_{t=i+1}^{k-1}(\mu(\cal H_{t,k-1})-1)\prod_{i\le t\le k-1,t\not=j}\mu(\cal H_{t,k-1})
f_{k-1}\cdots f_i f_{k-1}\cdots f_j+\Ann_{ U(\lie g)}e_{\beta_{i,j}}.
\end{multline*}
Since~$\mu(h_i)=\lambda(h_i)+1>0$, $\mu(h_j)=\lambda(h_j)+1>0$, the monomials
$$
f_{k-2}\cdots f_i f_{k-1}\cdots f_j,\, f_{k-2}\cdots f_j f_{k-1}\cdots f_i
$$
are $\mu$-standard and hence the vectors
$$
u_1=f_{k-2}\cdots f_i f_{k-1}\cdots f_j v_\mu,\, u_2=f_{k-2}\cdots f_j f_{k-1}\cdots f_i v_\mu
$$
are linearly independent by~\corref{A45}. Clearly, $e_{k-1}^2 u_1=e_{k-1}^2 u_2=0$, while $h_{k-1} u_r=(\mu(h_{k-1})+2)u_r$, $r=1,2$.
Since~$\mu(h_{k-1})+2=\lambda(h_{k-1})+2\ge 4$, it follows from the standard $\lie{sl}_2$-theory that~$f_{k-1} u_1$, $f_{k-1} u_2$ are non-zero and linearly independent.
Therefore,
\begin{multline*}
\bou_{\beta_{k,k},\beta_{i,j}}(\mu) e_{\beta_{i,j}}=
(\mu(\cal H_{i,k-1})-1)^{-1} (\mu(\cal H_{j,k-1}))^{-1} \left(\frac{\mu(\cal H_{i,k-1})-2}{\mu(\cal H_{i,k-1})}+1\right) e_{\beta_{k,k}}
\\
=2 (\mu(\cal H_{j,k-1})\mu(\cal H_{i,k-1}))^{-1} e_{\beta_{k,k}}.
\end{multline*}
Thus,
\begin{multline*}
\Pi_\lambda(\beta_{i,j},\beta_{k,k})=
e_{\beta_{k,k}}\tensor e_{\beta_{i,j}}-
 (y_\lambda+1)^{-1} e_{\beta_{j,k}}\tensor e_{\beta_{i,k}}
\\
-((x_\lambda+y_\lambda+2)(y_\lambda+1))^{-1}\Big(y_\lambda  e_{\beta_{i,k}}\tensor  e_{\beta_{j,k}}-
2 e_{\beta_{i,j}}\tensor e_{\beta_{k,k}}\Big).
\end{multline*}
We immediately obtain the relations~\eqref{C160.iii.4a} and~\eqref{C160.iii.4b} 
(note that in this case~$y_\lambda>1$) and $\mathfrak R_\Psi(\lambda,\lambda+\eta)$ is easily seen to be generic.
Finally, suppose that~$0<t_{\lambda,\eta}<4$. Using~\ref{C10}\cref{C10.c4} we conclude that if~$t_{\lambda,\eta}=2$,
the unique relation is~\eqref{C160.iii.2},  while if~$t_{\lambda,\eta}=1$, the unique path is a relation.
\end{pf}
\subsection{}\label{C180}
Finally, we consider the case when~$m_\eta=6$, that is
$\eta=\beta_{i,j}+\beta_{k,l}=\beta_{i,k}+\beta_{j,l}=\beta_{i,l}+\beta_{j,k}\in\Psi+\Psi$, $i<j<k<l\in I$.
Let~$x_\lambda=\lambda(\cal H_{i,j-1})$, $y_\lambda=\lambda(\cal H_{j,k-1})$,
$z_\lambda=\lambda(\cal H_{k,l-1})$. Then~$\lambda(\cal H_{i,k-1})=x_\lambda+y_\lambda+1$, $\lambda(\cal H_{j,l-1})=y_\lambda+z_\lambda+1$ and
$\lambda(\cal H_{i,l-1})=x_\lambda+y_\lambda+z_\lambda+2$.
Note that if~$t_{\lambda,\eta}=6$, we have~$x_\lambda,y_\lambda,z_\lambda>0$.

All technical difficulties in computing the $\Pi_{\lambda}(\beta,\beta')$ which occur here have already been discussed and we omit the details.
Suppose first that~$t_{\lambda,\eta}=6$.
We have
\begin{align*}
&\Pi_{\lambda}(\beta_{k,l},\beta_{i,j})=e_{\beta_{i,j}}\tensor
e_{\beta_{k,l}},\\
&\Pi_{\lambda}(\beta_{j,l},\beta_{i,k})=e_{\beta_{i,k}}\tensor e_{\beta_{j,l}}-(y_\lambda+1)^{-1} e_{\beta_{i,j}}\tensor e_{\beta_{k,l}},\\
&\Pi_{\lambda}(\beta_{j,k},\beta_{i,l})=
e_{\beta_{i,l}}\tensor
e_{\beta_{j,k}}-(z_\lambda+1)^{-1} e_{\beta_{i,k}}\tensor e_{\beta_{j,l}}
-\frac{z_\lambda}{(z_\lambda+1)(y_\lambda+z_\lambda+2)} \,e_{\beta_{i,j}}\tensor e_{\beta_{k,l}},\\
&\Pi_{\lambda}(\beta_{i,l},\beta_{j,k})=
e_{\beta_{j,k}}\tensor e_{\beta_{i,l}}-
(x_\lambda +1)^{-1} e_{\beta_{i,k}}\tensor e_{\beta_{j,l}}
-\frac{x_\lambda }{(x_\lambda +1)(x_\lambda+y_\lambda+2)} e_{\beta_{i,j}}\tensor e_{\beta_{k,l}},\\
&\Pi_{\lambda}(\beta_{i,k},\beta_{j,l})=
e_{\beta_{j,l}}\tensor e_{\beta_{i,k}}-(z_\lambda+1)^{-1} e_{\beta_{j,k}}\tensor e_{\beta_{i,l}}-
(x_\lambda +1)^{-1} e_{\beta_{i,l}}\tensor e_{\beta_{j,k}}
\\&\phantom{\Pi_{\lambda}(\beta_{i,k}}+
((x_\lambda +1)(z_\lambda+1))^{-1} e_{\beta_{i,k}}\tensor e_{\beta_{j,l}}
-\frac{x_\lambda z_\lambda}{(x_\lambda +1)(z_\lambda+1)(x_\lambda+y_\lambda+z_\lambda+3)}\, e_{\beta_{i,j}}\tensor e_{\beta_{k,l}},\\
&\Pi_\lambda(\beta_{i,j},\beta_{k,l})=
e_{\beta_{k,l}}
\tensor e_{\beta_{i,j}}
-(y_\lambda+1)^{-1} e_{\beta_{j,l}}\tensor e_{\beta_{i,k}}
\\&\phantom{\Pi_\lambda(\beta_{i,j}}-
\frac{y_\lambda}{(y_\lambda+1)(y_\lambda+z_\lambda+2)}\, e_{\beta_{j,k}}\tensor
e_{\beta_{i,l}}
-\frac{y_\lambda}{(y_\lambda+1)(x_\lambda+y_\lambda+2)}\, e_{\beta_{i,l}}\tensor e_{\beta_{j,k}}\\
&\phantom{\Pi_\lambda(\beta_{i,j}}
-
\frac{x_\lambda+y_\lambda+z_\lambda+3+(x_\lambda+y_\lambda+1)(y_\lambda+1)(y_\lambda+z_\lambda+1)}{(y_\lambda+1)(x_\lambda+y_\lambda+2) (y_\lambda+z_\lambda+2)(x_\lambda+y_\lambda+z_\lambda+3)}\,
e_{\beta_{i,k}}\tensor e_{\beta_{j,l}}\\
&\phantom{\Pi_\lambda(\beta_{i,j}}+
\frac{(x_\lambda+y_\lambda+1)(y_\lambda+z_\lambda+1)+(y_\lambda+1)(x_\lambda+y_\lambda+z_\lambda+3)}{(y_\lambda+1)(x_\lambda+y_\lambda+2) (y_\lambda+z_\lambda+2)(x_\lambda+y_\lambda+z_\lambda+3)}\,
e_{\beta_{i,j}}\tensor e_{\beta_{k,l}}.
\end{align*}

Denote the paths in~${\Delta_\Psi}(\lambda,\lambda+\eta)$ by~$\bop_r$, $1\le r\le 6$, where the numbering corresponds to the order in which they
appear above. Clearly, none of these paths is a relation.
A direct computation shows that
$\mathfrak R_\Psi(\lambda,\lambda+\eta)$ is spanned by
\begin{align*}
\bor_1&=(x_\lambda+1) (y_\lambda+2) (x_\lambda-z_\lambda) (z_\lambda+1) \bop_1+y_\lambda
   (x_\lambda+y_\lambda+1) (x_\lambda-z_\lambda) (y_\lambda+z_\lambda+2)
   \bop_2\\&+z_\lambda (x_\lambda+1) (y_\lambda+1) (x_\lambda+y_\lambda+2)
   (y_\lambda+z_\lambda+1) \bop_3\\&-x_\lambda (y_\lambda+1)(z_\lambda+1)
   (x_\lambda+y_\lambda+1)  (y_\lambda+z_\lambda+2) \bop_4,
\\
\bor_2&=(x_\lambda+1) (z_\lambda+1)(x_\lambda+y_\lambda+2)  (x_\lambda+z_\lambda+2)
   (y_\lambda+z_\lambda+3) \bop_1\\&+y_\lambda(x_\lambda+1)(z_\lambda+2)
   (x_\lambda+y_\lambda+1)  (y_\lambda+z_\lambda+2)
   (x_\lambda+y_\lambda+z_\lambda+3) \bop_2\\&-(y_\lambda+1)
   (x_\lambda+y_\lambda+2) (x_\lambda+z_\lambda+2) (y_\lambda+z_\lambda+1)
   (x_\lambda+y_\lambda+z_\lambda+3) \bop_3\\&-x_\lambda (y_\lambda+1)(z_\lambda+1)
   (x_\lambda+y_\lambda+2)  (y_\lambda+z_\lambda+2)
   (x_\lambda+y_\lambda+z_\lambda+2) \bop_5,
\\
\intertext{and}
\bor_3&=(y_\lambda+2) (z_\lambda+1)(x_\lambda+y_\lambda+2)  (y_\lambda+z_\lambda+3)
   (x_\lambda+y_\lambda+z_\lambda+3) \bop_1\\&-(z_\lambda+2)(x_\lambda+y_\lambda+1)
    (y_\lambda+z_\lambda+2) (x_\lambda+2 y_\lambda+z_\lambda+4)
   \bop_2\\&-z_\lambda (y_\lambda+1) (x_\lambda+y_\lambda+z_\lambda+3) (x_\lambda+2
   y_\lambda+z_\lambda+4) \bop_3\\&-(y_\lambda+1)(z_\lambda+1) (x_\lambda+y_\lambda+1)
    (y_\lambda+z_\lambda+2) (x_\lambda+y_\lambda+z_\lambda+2)
   \bop_6.
\end{align*}
One can now check that~$\mathfrak R_\Psi(\lambda,\lambda+\eta)$ is generic unless $x_\lambda=z_\lambda$. In the latter case the
first relation reduces to~$\bop_3-\bop_4$. 

Finally, we list the relations in cases when~$0<t_{\lambda,\eta}<6$. 
By~\ref{C10}\cref{C10.c5}, we have three cases with~$t_{\lambda,\eta}=2$.
If~$y_\lambda=0$, or equivalently~$j=k-1$ and~$\lambda(h_j)=0$, while $x_\lambda,z_\lambda>0$, the relation is
$$
(z_\lambda+2) (x_\lambda+z_\lambda+4) \bop_2+(z_\lambda+1)
   (x_\lambda+z_\lambda+2) \bop_6.
$$
If~$x_\lambda=0$, $y_\lambda,z_\lambda>0$ the relation is
$$
(y_\lambda+1) (z_\lambda+2) (y_\lambda+z_\lambda+3) \bop_4+(y_\lambda+2)
   z_\lambda (y_\lambda+z_\lambda+2) \bop_5.
$$
If~$z_\lambda=0$, $x_\lambda,y_\lambda>0$ the relation is
$$
(x_\lambda+2) (y_\lambda+1) (x_\lambda+y_\lambda+3) \bop_3+x_\lambda
   (y_\lambda+2) (x_\lambda+y_\lambda+2) \bop_5.
$$
Thus, in all these cases the space~$\mathfrak R_\Psi(\lambda,\lambda+\eta)$ is generic. 
Finally, if~$x_\lambda=z_\lambda=0$, $y_\lambda>0$ the unique path~$\bop_5$ in~${\Delta_\Psi}(\lambda,\lambda+\eta)$ is not a relation.

Thus, we obtain the following
\begin{prop}
Let~$i<j<k<l\in I$.
Suppose that~$\eta=\beta_{i,j}+\beta_{k,l}=\beta_{i,k}+\beta_{j,l}=\beta_{i,l}+\beta_{j,k}\in\Psi+\Psi$. Then
$\dim\mathfrak R_\Psi(\lambda,\lambda+\eta)=\lfloor |\Delta_\Psi(\lambda,\lambda+\eta)/2|\rfloor$ and
$$
\cal N_\eta\subset P^+\cap \{ \xi\in \lie h^*\,:\, \xi(\cal H_{i,j-1}-\cal H_{k,l-1})=0\}
$$
and coincides with the latter set if~$\Psi$ is regular. 
\end{prop}

\subsection{}\label{concl}
Let~$\Psi=\Psi(i_1,\dots,i_k)$ be regular. It follows from Propositions~\ref{C120}, \ref{C140}, \ref{C160} and~\ref{C180} that the coefficients
in all relations in~$\mathfrak R_\Psi(\mu,\mu+\eta)$ depend on~$\mu(\cal H_{i_r,i_s-1})$, $1\le r<s\le n$. 

Let~$\lambda\in P^+$. By~\propref{C3}, $\Delta_\Psi[\lambda]$ is isomorphic to the quiver $\Xi_a(\boldsymbol{m})$, where
$a=\sum_{r=1}^k\lambda(h_{i_r})\pmod2$, $\boldsymbol{m}=(m_1,\dots,m_k)$, $m_r=\lambda(h_{i_r-1})+\lambda(h_{i_r})$. Let
$$
\zeta_r(\lambda)=\lambda(\cal H_{i_r+1,i_{r+1}-2})+2,\qquad 1\le r<k.
$$
Let~$(x_1,\dots,x_k)$ be the image of~$\mu\in\Delta_\Psi[\lambda]_0$ in~$\Xi_a(\boldsymbol{m})_0$. Then
$$
\mu(\cal H_{i_r,i_s-1})=x_r-x_s+\sum_{p=r+1}^s m_p+\zeta_{p-1}(\lambda)+r-s-1.
$$
Thus, the isomorphism of algebras~$\bt^{\lie g}_\Psi\to\Delta_\Psi$ gives rise to a family of relations on quivers $\Xi_a(\boldsymbol{m})$. The relations,
and in particular their genericity, depend on a family of positive integer parameters $\zeta_p(\lambda)$.
The resulting algebras
are Koszul and of global dimension at most $p(p+1)/2$, where~$p=\#\{j\,:\, m_j>0\}$. It is finite dimensional if and only if~$i_1>1$.
The explicit relations can be easily written down using the above Propositions.

\section*{List of notations}
\begin{tabular}{llllllll}
$I$&\ref{MR10}&
$\alpha_i$, $\varpi_i$&\ref{MR10}&$\varepsilon_i$, $\varphi_i$, $\varepsilon$, $\varphi$&\ref{MR10}
&$R$, $P$, $R^+$, $P^+$&\ref{MR10}\\ 
$\lie n^\pm_\Psi$, $\lie n^\pm$, $\lie b$&\ref{MR10}&$V(\lambda)$&\ref{MR30}&$\bv$, $\bv^{\circledast}$&\ref{MR30}&$\ba$, $\bt$, $\bs$, $\be$&\ref{MR30}\\
$1_\lambda\,:\,\lambda\in P^+$&\ref{MR30}&$\le_\Psi$, $\le$ &\ref{MR50}&
$d_\Psi(\lambda,\mu)$&\ref{MR50}&$\ba_\Psi^{\lie g}(F)$, $\ba_\Psi^{\lie g}$&\ref{MR50}\\$\le_\Psi\lambda$, $\lambda\le_\Psi$, $[\lambda,\mu]_\Psi$&\ref{MR50}&
$\Delta_0$, $\Delta_1$, $\bar\Delta$&\ref{PRE0}&$x^\pm$, $x\in\Delta_0$&\ref{PRE0}&$\bc\Delta$&\ref{PRE0}\\$\Delta_\Psi(F)$, $\Delta_\Psi$&\ref{MR70}&
$\mathfrak R_\Psi(\lambda,\lambda+\eta)$&\ref{MR75}&$m_\eta$, $t_{\lambda,\eta}$&\ref{M75}&$\cal N_\eta$&\ref{M75}\\$|\boldsymbol{x}|$, $\boldsymbol{x}\in\bz_+^r$&\ref{M95}&
$\boldsymbol{e}_i^{(r)}$&\ref{M95}&
$\Xi(\boldsymbol{m})$, $\Xi_a(\boldsymbol{m})$&\ref{M95}&$e_i$, $f_i$, $h_i$&\ref{PRE100}\\$v_\lambda$, $\xi_{-\lambda}$, $M^\lambda$&\ref{PRE100}&$\Pi_\lambda(\beta,\beta')$&\ref{SP30}&
$F_\beta(\lie h)$&\ref{COMP110}&$\pi_\lambda$, $\pi_{\lambda,\beta}$&\ref{COMP110}\\$\bou_{\beta,\gamma}$, $\bou_{\beta,\gamma}(\lambda)$&\ref{SP55}&$\alpha_{i,j}$&\ref{A48}&$\Sigma(i,j)$&\ref{A48}&$\bof_\sigma$&\ref{A48}\\
$\psi_\eta$&\ref{A50}&$\cal H_{i,j}$&\ref{A50}&$\cal X^\pm_{i,j,k}$, $\cal X^\pm_{i,j}$&\ref{A50}&
$\Gamma_a(\boldsymbol{m},\boldsymbol{n})$&\ref{A32}\\$e_{i,j}$&\ref{A90}&$\bou_{\alpha_{i,j},\alpha_{p,q}}$&\ref{A90}&$\cal Z_{\alpha_{i,j},\Psi}$&
\ref{A95}&$\beta_{i,j}$&\ref{C0}\\
$e_{\beta_{i,j}}$&\ref{C50}&$\cal U_{r,s,i,j}$&\ref{C60}&$\bar{\cal U}_{r,s,i,j}$&\ref{C115}&$\bou_{\beta_{i,j},\beta_{r,s}}$&\ref{C70}\\
\end{tabular}

\end{document}